\documentclass[11pt]{amsart}
\usepackage{amssymb,latexsym, amsmath, amsxtra}
\usepackage[dvips]{graphics}
\theoremstyle{plain}
\newtheorem{theorem}{Theorem}[section]

\newtheorem{conjecture}[theorem]{Conjecture}
\newtheorem{remark}[theorem]{Remark}

\theoremstyle{definition}

\theoremstyle{remark}

\numberwithin{equation}{section}

\begin{document}

\title{Applications of the $L$-functions ratios conjectures}

\abstract In upcoming papers by Conrey, Farmer and Zirnbauer there
appear conjectural formulas for averages, over a family,
 of ratios of products of shifted $L$-functions.
In this paper we will present various applications of these ratios
conjectures to a wide variety of problems that are of interest in
number theory, such as lower order terms in the zero statistics of
$L$-functions, mollified moments of $L$-functions and discrete
averages over zeros of the Riemann zeta function. In particular,
using the ratios conjectures we easily derive the answers to a
number of notoriously difficult computations.
\endabstract

\author {J.B. Conrey}
\address{American Institute of Mathematics,
360 Portage Ave, Palo Alto, CA 94306} \address{School of
Mathematics, University of Bristol, Bristol, BS8 1TW, United
Kingdom} \email{conrey@aimath.org}

\author{N.C. Snaith}
\address{School of Mathematics,
University of Bristol, Bristol, BS8 1TW, United Kingdom}
\email{N.C.Snaith@bris.ac.uk}

\thanks{
Research of the first author supported by the American Institute
of Mathematics and a Focused Research Group grant from the
National Science Foundation. The second author was  supported by a
Royal Society Dorothy Hodgkin and an EPSRC Advanced Research
Fellowship.  This paper was started when both authors were
visiting the Isaac Newton Institute for Mathematical Sciences for
the progamme ``Random matrix approaches  in number theory".
 }
\maketitle

\tableofcontents

\section{Introduction}
\label{sec:intro}

Applications of random matrix theory in number theory began with
Montgomery's pair correlation conjecture [36].  In this paper
Montgomery conjectured that, in the limit for large height up the
critical line, any local statistic of the zeros of the Riemann
zeta function is given by the corresponding statistic for
eigenvalues from the GUE ensemble of random matrix theory [34]. A
local statistic is one that involves only correlations between
zeros separated on a scale of a few mean spacings. Odlyzko checked
the statistics numerically for the pair correlation and the
nearest neighbour spacing distribution and found spectacular
agreement [37].  At leading order the zero statistics and
eigenvalues statistics are identical; asymptotically no factors of
an arithmetical nature appear. However, it is clear that
arithmetical contributions play a role in lower order terms, and
Bogomolny and Keating [3] identified these in the case of the pair
correlation function.

Katz and Sarnak [27] proposed that local statistics of zeros of
families of $L$-functions could be modelled by the eigenvalues of
matrices from the classical compact groups with Haar measure.  In
this way each family of $L$-functions is believed to have a
symmetry type: unitary, symplectic or orthogonal.  Iwaniec, Luo
and Sarnak [24] calculated the leading asymptotics for the
one-level densities (using test functions whose Fourier transforms
have limited support) for families of $L$-functions with each
symmetry type and found agreement with random matrix theory. Again
these leading terms had no arithmetic part.

More recently, random matrix theory has been applied to the
moments of $L$-functions averaged over a family.  These are
global, rather than local, statistics.  A characteristic feature
of a global statistic is that an arithmetic factor appears in the
leading order term.  In the original papers [29,30]  the leading
term was a product of the corresponding moment of a characteristic
polynomial from random matrix theory and a seemingly independent
Euler product.

Often in random matrix theory one can calculate such global
statistics exactly for any finite matrix size $N$.  In particular,
when evaluating moments of characteristic polynomials one obtains
an exact asymptotic expansion in $N$ as $N\to \infty$.  We now
understand conjecturally the analogue for moments of $L$-functions
and in particular how the arithmetic and random matrix factors
interact in the lower order terms.   For any family of
$L$-functions we can conjecture [10] an asymptotic expansion for
any moment which we believe is accurate essentially to the square
root of the size of the family.

A natural way to generalize these moment formulae is to consider
averages of ratios of products of $L$-functions or characteristic
polynomials. In two forthcoming papers [11,12] there appear
conjectural formulas for averages, over a family,  of ratios of
products of shifted $L$-functions. Those papers contain several
applications of these conjectures, as well as theorems proving the
 random matrix analogues of these conjectures. In [13] and
[5], different proofs of the random matrix theorems are given,
although not for the full range of the main parameter, the
dimension of the matrix.

The point is that these ratios conjectures are useful for
calculating both local and global statistics.  In fact quoting
from [4] ``The averages of products and ratios of characteristic
polynomials are more fundamental characteristics of random matrix
models than the correlation functions."  We would argue the same
can be said for $L$-functions.  From the ratios conjectures not
only can you obtain all $n$-level correlations, but also
essentially any local or global statistic.  An important feature
on the number theory side is that this includes all lower order
terms, in particular it shows the arithmetic contribution present
in local statistics.


In this paper we will present various applications of these ratios
conjectures.  In Section \ref{sect:ratios} we give the precise
statement and sketch the derivations of some examples of the
ratios conjecture for each of the three symmetry types: unitary,
symplectic and orthogonal. These examples, which have one or two
$L$-functions in the numerator and denominator, cover most of the
cases that we need in the applications in this paper, but the
conjectures are more general in that they can involve any number
of $L$-functions [11].  Theorem \ref{theo:sympderiv} and
\ref{theo:orthderiv} give auxiliary formulae useful in calculating
the most basic local statistic, the one-level density. In Section
\ref{sect:onelevel} we then show how the ratios conjecture can be
used to compute the one-level density of the simplest family of
$L$-functions with symplectic symmetry, namely Dirichlet
$L$-functions with real quadratic characters.  We state a similar
result for the orthogonal family associated with quadratic twists
of the Ramanujan $\tau$-function. In the following section we
consider lower order terms in the pair correlation of the zeros of
the Riemann zeta-function.   As mentioned above, Bogomolny and
Keating were the first to find these lower order terms; their
heuristic method involved a careful analysis of the
Hardy-Littlewood conjectures for prime pairs.  The strength of our
method is that it allows us to avoid such detailed considerations.

The next two sections consider averages of mollified
$L$-functions.   Mollifiers are used to obtain information about
small values of $L$-functions, in particular zeros.  Mollifiers
were first introduced in the context of the Riemann zeta function
to bound the number of zeros in a vertical strip to the right of
the half-line (i.e. zero density results).  Subsequently Selberg,
and then Levinson, obtained lower bounds for the proportion of
zeros satisfying the Riemann Hypothesis by mollifying zeta in the
neighbourhood of the critical line.  Recent uses have focused on
obtaining non-vanishing results at the central point for families
of $L$-functions.  All of these results involve complicated
analysis, for example Levinson's asymptotic evaluation of the
mollified second moment of zeta takes nearly fifty pages.  Before
embarking on such a calculation it would be useful to know ahead
of time what the answer is.  In Section \ref{sect:moll} we show
how to obtain these answers quickly.  For each of the families
that we've introduced we calculate the mollified second moment of
arbitrary linear combinations of derivatives and reveal the simple
structure of the result.  In all cases where these have been
rigourously calculated (only accomplished when the mollifier is
sufficiently short), these results are in agreement. In Section
\ref{sect:mollk} we show how to mollify any moment of the Riemann
zeta function and give detailed expressions in the case of the
fourth moment; none of these, apart from the second moment, have
been calculated without using the ratios conjecture.   It is
interesting to note that unlike other averages considered in this
paper, there does not seem to be a random matrix analogue of
mollifying as there is nothing that naturally corresponds to a
partial Dirichlet series.

Another kind of average which gives useful information about the
distribution of zeros is a discrete moment summing the zeta
function, or its derivatives, at or near the zeros. In Section
\ref{sect:discrete} we consider   moments of $|\zeta'(\rho)|$ and
$|\zeta(\rho+a)|$. Using the ratios conjecture we show how to
obtain all of the lower-order terms for these averages. While the
leading order terms had previously been conjectured or proved, it
was not known how to obtain these lower order terms.

In  Section \ref{sect:conn} we show how to use the ratios
conjecture to reproduce the asymptotic formulae used to obtain
non-vanishing results for various families. In addition we sketch
how one should go about proving that the proportion of
non-vanishing for the $k$th derivative $L^{(k)}(1/2,\chi)$
approaches 100\% as $k\to \infty$ for the family of all Dirichlet
$L$-functions; by contrast, in [35] a convincing argument is made
that one can't do better than 2/3 non-vanishing for
$\Lambda^{(k)}(1/2,\chi)$, where $\Lambda$ is the completed
$L$-function, without mollifying a higher power than the second.

In short, there are a number of difficult computations which the
ratios conjectures simplify significantly.  A few of these
computations have the property that they could be made into
theorems by proceeding alternatively; some are purely conjectural.
However, even for those that could be proved by other methods,
knowing the answer ahead of time is useful as a guide along the
way, a check at the end and even in deciding whether to commence
what could be a painful calculation.

Throughout this paper we assume the Riemann Hypothesis for all the
$L$-functions that arise.

\section{Ratios conjectures}\label{sect:ratios}

\subsection{A unitary example}

 An example of a basic conjecture for the zeta-function
follows. This was the example that Farmer first considered when
formulating his initial conjecture about averages of ratios of
zeta functions with shifts.
 With $s=1/2+it$, let
 \begin{eqnarray}
 \label{eq:Rzeta}
 R_\zeta(\alpha,\beta,\gamma,\delta):= \int_0^T
 \frac{\zeta(s+\alpha)\zeta(1-s+\beta)}{\zeta(s+\gamma)\zeta(1-s+\delta)}~dt.
\end{eqnarray}
Farmer [16] conjectured that for $\alpha,\beta,\gamma,\delta <<
1/\log T$,
\begin{eqnarray} \label{eq:farmer}
R_\zeta(\alpha,\beta,\gamma,\delta)\sim
T\frac{(\alpha+\delta)(\beta+\gamma)}{(\alpha+\beta)(\gamma+\delta)}-T^{1-\alpha-\beta}
\frac{(\delta-\beta)(\gamma-\alpha)}{(\alpha+\beta)(\gamma+\delta)},
\end{eqnarray}
as $T\rightarrow \infty$, provided that $\Re \gamma, \Re \delta
>0$. Our ratios conjecture gives us a recipe for computing a more
precise conjecture for $R_{\zeta}$.

 Briefly, we use the approximate functional
equation
\begin{equation}\label{eq:approxfunczeta}
\zeta(s)=\sum_{n\leq X} \frac{1}{n^s} + \chi(s)\sum_{n\leq Y}
\frac{1}{n^{1-s}}+ \text{remainder},
\end{equation}
where $s=\sigma+it$, $\chi(s)=2^s\pi^{s-1}
\sin(s\pi/2)\Gamma(1-s)$ and $XY=\tfrac{t}{2\pi}$,
 for the zeta functions in the
numerator and ordinary Dirichlet series expansions for those in
the denominator:
\begin{equation}
\frac{1}{\zeta(s)}= \sum_{n=1}^{\infty} \frac{\mu(n)}{n^s}.
\end{equation}
  We only use the pieces which have the same
number of $\chi(s)$ as $\chi(1-s)$ and
 we integrate term-by-term, retaining only the diagonal pieces.
We then complete all of the sums that we arrive at.

Thus, the term from the ``first'' part of the two approximate functional equations gives $T$ times
\begin{eqnarray}
 &&\sum_{hm=kn} \frac{\mu(h)\mu(k)}{m^{1/2+\alpha} n^{1/2+\beta} h^{1/2+\gamma} k^{1/2+\delta}}\nonumber\\
  &&\qquad\qquad=\prod_p \sum_{h+m=k+n} \frac{\mu(p^h)\mu(p^k)}{p^{(1/2+\alpha)m+(1/2+\beta)n +(1/2+\gamma)h
 +(1/2+\delta)k}}.
\end{eqnarray}
The only possibilities for $h$ and $k$ here are 0 and 1. Thus, we easily find that the right-hand sum above is
\begin{eqnarray}
=\frac{1}{\left(1-\frac{1}{p^{1+\alpha+\beta}}\right)}
\left(1-\frac{1}{p^{1+\beta+\gamma}}
-\frac{1}{p^{1+\alpha+\delta}}+\frac{1}{p^{1+\gamma+\delta}}\right);
\end{eqnarray}
 thus, the product over primes is
\begin{eqnarray}
\frac{\zeta(1+\alpha+\beta)\zeta(1+\gamma+\delta)}
{\zeta(1+\alpha+\delta)\zeta(1+\beta+\gamma)}A_\zeta(\alpha,
\beta, \gamma, \delta),
\end{eqnarray}
where
\begin{eqnarray}\label{eq:Azeta}
A_\zeta(\alpha, \beta, \gamma, \delta)=\prod_p
\frac{\left(1-\frac{1}{p^{1+\gamma+\delta}}\right)
\left(1-\frac{1}{p^{1+\beta+\gamma}}
-\frac{1}{p^{1+\alpha+\delta}}+\frac{1}{p^{1+\gamma+\delta}}\right)}
{\left(1-\frac{1}{p^{1+\beta+\gamma}}\right)
\left(1-\frac{1}{p^{1+\alpha+\delta}}\right)}.
\end{eqnarray}
The other term comes from the second piece of each approximate functional equation and is similar to the first piece except that $\alpha$ is replaced by $-\beta$ and $\beta$ is replaced by $-\alpha$. Also, because of the $\chi$-factors in the functional equation, we have an extra factor of
\begin{eqnarray}\label{eqn:chi}
\chi(s+\alpha) \chi(1-s+\beta) = \left(\frac {t}{2
\pi}\right)^{-\alpha-\beta}\left(1+O\left(\frac1{|t|}\right)\right).
\end{eqnarray}
Thus, the more precise ratios conjecture gives
\begin{conjecture} [{Conrey, Farmer and Zirnbauer [11]}] \label{conj:unit} With
constraints on $\alpha, \beta, \gamma$ and $\delta$ as described
below at (\ref{eq:constraints}), we have
\begin{eqnarray} \label{eqn:R}
 R_\zeta(\alpha,\beta,\gamma,\delta) \nonumber &=&\int_0^T
 \frac{\zeta(s+\alpha)\zeta(1-s+\beta)}{\zeta(s+\gamma)\zeta(1-s+\delta)}~dt\\
 &=&\int_0^T \left(\frac{\zeta(1+\alpha+\beta)\zeta(1+\gamma+\delta)}
{\zeta(1+\alpha+\delta)\zeta(1+\beta+\gamma)}A_\zeta(\alpha, \beta, \gamma, \delta) \right. \\
&& \quad \left. +\left(\frac {t}{2 \pi}\right)^{-\alpha-\beta}
\frac{\zeta(1-\alpha-\beta)\zeta(1+\gamma+\delta)}
{\zeta(1-\beta+\delta)\zeta(1-\alpha+\gamma)}A_\zeta(-\beta,
-\alpha,  \gamma, \delta) \right)~dt
+O\left(T^{1/2+\epsilon}\right),\nonumber
\end{eqnarray}
where $A_{\zeta}$ is defined at (\ref{eq:Azeta}).
\end{conjecture}

In the following sections we have similar conjecture for ratios of
$L$-functions averaged over various families. In these families
the $L$-functions are indexed by an integer $d$ and we consider
averages for $d<X$.  In all of these examples we constrain the
shifts as follows. For $\alpha$ a generic shift in the numerator
($\alpha$ and $\beta $ in the above example) and $\delta$ a
generic shift in the denominator, we require
\begin{subequations}\label{eq:constraints}
\begin{align}
 &-\frac{1}{4} < \Re \alpha< \frac{1}{4}\\
&\frac{1}{\log C} \ll \Re
\delta<\frac{1}{4}\label{eq:constraintsb}\\
&\Im \alpha,\Im \delta\ll_{\epsilon} C^{1-\epsilon}\quad \mbox{
(for every $\epsilon>0$)}
\end{align}
\end{subequations}
where $C=T$ in the above example and $C=X$ in the case of discrete
families of $L$-functions.  In conjectures that refer to these
conditions the error terms are believed to be uniform in the above
range of parameters.
\begin{remark}
\label{rem:constraints} Equation (\ref{eq:constraintsb}) can be
relaxed if for each shift in the denominator going to zero there
is a corresponding shift in the numerator going to zero at the
same rate.
\end{remark}
\begin{remark}
The bound of 1/4 on the absolute values of the real parts of the
shifts are to prevent divergence of  the Euler products that
appear in the ratios conjectures.
\end{remark}
\begin{remark}\label{rem:diff}
Because of the uniformity in the parameters
$\alpha,\beta,\gamma,\delta$ we can differentiate our conjectural
formulas with respect to these parameters and the results are
valid with the same range and error terms.
\end{remark}

For obtaining lower order terms in pair correlation in Section
\ref{sect:pair}, we need the following:
\begin{theorem} \label{theo:zetaderiv}Assuming Conjecture \ref{conj:unit}, we have
\begin{eqnarray}&&
\int_0^T
\frac{\zeta'}{\zeta}(s+\alpha)\frac{\zeta'}{\zeta}(1-s+\beta)~dt
=\int_0^T\left( \left(\frac{\zeta'}{\zeta}\right)'(1+\alpha+\beta)+\right. \\
&&\quad \left(\frac{t}{2\pi}\right)^ {-\alpha-\beta}
\zeta(1+\alpha+\beta)\zeta(1-\alpha-\beta)
\prod_p\frac{(1-\frac{1}{p^{1+\alpha+\beta}})
(1-\frac 2 p +\frac{1}{p^{1+\alpha+\beta}})}{(1-\frac 1 p )^2}\nonumber\\
&& \quad \left. -\sum_p \left(\frac{\log
p}{(p^{1+\alpha+\beta}-1)} \right)^2\right)~dt
+O(T^{1/2+\epsilon}),\nonumber
\end{eqnarray}
provided that $\frac{1}{\log T}\ll\Re \alpha, \Re \beta
<\frac{1}{4}$.
\end{theorem}

This theorem follows from (\ref{eqn:R}) by differentiating with
respect to $\alpha$ and $\beta$ and setting $\gamma=\alpha$ and
$\delta = \beta$.  To perform this calculation, it is helpful to
observe that $A(\alpha, \beta, \alpha, \beta)=1$. Also, when
differentiating the second term on the right side of (\ref{eqn:R})
it is useful to observe that for a function $f(z,w)$ which is
analytic at $(z,w)=(\alpha,\alpha)$,
\begin{equation}
\frac{d}{d\alpha}
\frac{f(\alpha,\gamma)}{\zeta(1-\alpha+\gamma)}\bigg|_{\gamma=\alpha}=-f(\alpha,\alpha).
\end{equation}

\subsection{Symplectic examples} As a second example we consider the family of Dirichlet $L$-functions $L(s,\chi_d)$ associated with real, even, Dirichlet characters $\chi_d$. Let
 \begin{eqnarray}
R_D(\alpha,\beta;\gamma,\delta):=\sum_{d\le X}
\frac{L(1/2+\alpha,\chi_d)L(1/2+\beta,\chi_d)}{L(1/2+\gamma,\chi_d)L(1/2+\delta,\chi_d)},
\end{eqnarray} with the usual conditions (\ref{eq:constraints}) on the shifts
$\alpha,\beta,\gamma$ and $\delta$.
Let us also consider the simpler example
\begin{eqnarray}
R_D(\alpha;\gamma):=\sum_{d\le X} \frac{L(1/2+\alpha,\chi_d)
}{L(1/2+\gamma,\chi_d) }. \end{eqnarray}  As part of our recipe,
we replace the $L(s,\chi_d)$ in the numerator by approximate
functional equation
\begin{equation}
L(\tfrac{1}{2}+\alpha,\chi_d)=\sum_{m<x}
\frac{\chi_d(m)}{m^{1/2+\alpha}}
+\left(\frac{d}{\pi}\right)^{-\alpha}
\frac{\Gamma(1/4-\alpha/2)}{\Gamma(1/4+\alpha/2)} \sum_{n<y}
\frac{\chi_d(n)}{n^{1/2-\alpha}} +\text{remainder},
\end{equation}
where $xy=d/(2\pi)$,
 and we replace the
$L(s,\chi_d)$ in the denominator by their infinite series:
\begin{equation}
\frac{1}{L(s,\chi_d)}=\sum_{h=1}^\infty \frac{\mu(h)\chi_d(h)}{h^s}.
\end{equation}

 We consider each of the $2^\lambda$ (if there are $\lambda$
factors in the numerator) pieces separately and average
term-by-term within those pieces. We only retain the terms where
we are averaging over squares; in other words we use the main part
of the formula
\begin{equation}\sum_{d\le X} \chi_d(n)=\bigg\{ {a(n) X^*+\mbox{ small if $n$ is a square} \atop \mbox{ small if $n$ is not a square}}
\end{equation}
where $X^*=\sum_{d\le X}1$ is the number of fundamental discriminants below $X$ and where
\begin{equation} a(n)=\prod_{p\mid n} \frac{p}{p+1}.
\end{equation}
After computing these `diagonal' terms, we complete the sums by
extending to infinity the ranges of the summation variables; we
identify these terms as ratios of products  of zeta functions
multiplied by absolutely convergent Euler products. The sum of
these expressions, one for each product of pieces of the
approximate functional equations, forms our conjectural answer.

Proceeding to details, let us first consider the simpler example
$R_{D}(\alpha;\gamma).$ We restrict attention to the `first' piece
of the approximate functional equation. Thus, we consider
\begin{eqnarray}
\sum_{d\le X} \sum_{h, m}\frac{\mu(h) \chi_d(hm)}{h^{1/2+\gamma}
m^{1/2+\alpha} }.
\end{eqnarray}
Retaining only the terms for which $hm$ is square, leads us to
\begin{eqnarray}
X^*\sum_{hm=\square} \frac{\mu(h) a(hm)}{h^{1/2+\gamma}
m^{1/2+\alpha} }.
\end{eqnarray}
We express this sum as an Euler product (to ``save'' variables we now replace $h$ by $p^h$ and $m$ by $p^m$):
\begin{eqnarray}
\prod_p\sum_{{h+m \atop \begin{tiny}
\mbox{even}\end{tiny}}}\frac{\mu(p^h)a(p^{h+m}) }{p^{h(1/2+\gamma)+
 m(1/2+\alpha) }}.
\end{eqnarray}
The effect of $\mu(p^h)$ is to limit the choices for $h$ to 0 or
1.   When $h=0$ we have
\begin{eqnarray} \sum_{{m \atop \begin{tiny}
\mbox{even}\end{tiny}}}\frac{a(p^m) }{p^{
 m(1/2+\alpha) }}=1+\frac{p}{p+1}\sum_{m=1}^\infty\frac{1}{p^{m(1+2\alpha)}}=1+\frac{p}{(p+1)}
\frac{1}{p^{1+2\alpha}}\frac{1}{(1-\frac{1}{p^{1+2\alpha}})},
\end{eqnarray}
and when $h=1$ there is a contribution of
\begin{eqnarray} \sum_{{m \atop \begin{tiny}
\mbox{odd}\end{tiny}}}-\frac{a(p^{m+1}) }{p^{1/2+\gamma}p^{
 m(1/2+\alpha) }}=-\frac{p}{(p+1)}\frac{1}{p^{1+\alpha+\gamma}}\frac{1}{(1-\frac{1}{p^{1+2\alpha}})} .\end{eqnarray}
Thus, the Euler product simplifies to
\begin{eqnarray}
\frac{\zeta(1+2\alpha)}{\zeta(1+\alpha+\gamma)}\prod_p\left(1-\frac{1}{p^{1+\alpha+\gamma}}
\right)^{-1}\left(1-\frac{1}{(p+1)p^{1+2\alpha}}-\frac{1}{(p+1)p^{\alpha+\gamma}}\right).
\end{eqnarray}
 The product over primes is absolutely convergent as long as $\Re \alpha, \Re \gamma >-1/4$.

The other piece can be determined by recalling the functional equation
\begin{eqnarray}
L(1/2+\alpha,\chi_d)=\left(\frac d \pi\right)^{-\alpha}
\label{eqn:functional}
\frac{\Gamma(1/4-\alpha/2)}{\Gamma(1/4+\alpha/2)}
L(1/2-\alpha,\chi_d).
\end{eqnarray}
Thus, in total we expect that
\begin{conjecture}[{Conrey, Farmer and Zirnbauer [11]}] \label{conj:symp} With
constraints on $\alpha$ and  $\gamma$ as described  at
(\ref{eq:constraints}), we have
\begin{eqnarray}
R_{D}(\alpha;\gamma)&=&\sum_{d\le X}  \nonumber \frac{L(1/2+\alpha,\chi_d)}{L(1/2+\gamma,\chi_d)}=\sum_{d\le X}\left(\frac{\zeta(1+2\alpha)}{\zeta(1+\alpha+\gamma)}A_D(\alpha;\gamma)\right.\\
&& \qquad \qquad + \left. \left(\frac d \pi\right)^{-\alpha}
\frac{\Gamma(1/4-\alpha/2)}{\Gamma(1/4+\alpha/2)}
\frac{\zeta(1-2\alpha)}{\zeta(1-\alpha+\gamma)}
A_D(-\alpha;\gamma)\right) +O(X^{1/2+\epsilon}),
\end{eqnarray}
where
\begin{eqnarray}
\label{eq:AD}
A_D(\alpha;\gamma)=\prod_p\left(1-\frac{1}{p^{1+\alpha+\gamma}}
\right)^{-1}\left(1-\frac{1}{(p+1)p^{1+2\alpha}}-\frac{1}{(p+1)p^{\alpha+\gamma}}\right).
\end{eqnarray}
\end{conjecture}

For applications to the one-level density in the next section, we
note that
\begin{eqnarray}
\sum_{d\le X}\frac{L'(1/2+r,\chi_d)}{L(1/2+r,\chi_d)}=
\frac{d}{d\alpha}R_D(\alpha;\gamma)\bigg|_{\alpha=\gamma=r}.
\end{eqnarray}
Now
\begin{eqnarray} \frac{d}{d\alpha}
 \frac{\zeta(1+2\alpha)}{\zeta(1+\alpha+\gamma)}
A_D(\alpha;\gamma)\bigg|_{\alpha=\gamma=r}=\frac{\zeta'(1+2r)}{\zeta(1+2r)}A_D(r;r)+A_D'(r;r)
\end{eqnarray}
and \begin{eqnarray}&&
 \frac{d}{d\alpha}
\left(\frac d \pi\right)^{-\alpha}
\frac{\Gamma(1/4-\alpha/2)}{\Gamma(1/4+\alpha/2)}
\frac{\zeta(1-2\alpha)}{\zeta(1-\alpha+\gamma)}
A_D(-\alpha;\gamma)\bigg|_{\alpha=\gamma=r}\\
&&\qquad =-\left(\frac d \pi\right)^{-r}
\frac{\Gamma(1/4-r/2)}{\Gamma(1/4+r/2)}
 \zeta(1-2r)  A_D(-r;r).\nonumber
\end{eqnarray}
Also, $A_D(r;r)=1,$
\begin{eqnarray}
A_D(-r;r)=
\prod_p\left(1-\frac{1}{(p+1)p^{1-2r}}-\frac{1}{p+1}\right)\left(1-\frac{1}{p}\right)^{-1},
\end{eqnarray}
and
\begin{eqnarray}
A'_D(r;r)=\sum_p \frac{\log p}{(p+1)(p^{1+2r}-1)}.
\end{eqnarray}

Thus, the ratios conjecture implies (see Remark \ref{rem:diff}):
\begin{theorem} \label{theo:sympderiv} Assuming Conjecture \ref{conj:symp},
$\frac{1}{\log X }\ll\Re r<\frac{1}{4}$ and $\Im r\ll_{\epsilon}
X^{1-\epsilon}$ we have
 \begin{eqnarray} \label{eqn:onelevel}
&&\sum_{d\le X}\frac{L'(1/2+r,\chi_d)}{L(1/2+r,\chi_d)}\nonumber \\
&&\qquad =
\sum_{d\le X} \bigg(\frac{\zeta'(1+2r)}{\zeta(1+2r)} +A_D'(r;r)-\left(\frac d \pi\right)^{-r}
\frac{\Gamma(1/4-r/2)}{\Gamma(1/4+r/2)}
 \zeta(1-2r)  A_D(-r;r)\bigg)\\
 &&\qquad \qquad +O(X^{1/2+\epsilon}),\nonumber
\end{eqnarray}
where $A_D(\alpha;\gamma)$ is defined in (\ref{eq:AD}).
\end{theorem}

 Now we look at the case of
two $L$-functions in the numerator and denominator.  Here we will
only work to keep the first main terms when the shifts $\alpha,
\beta, \gamma$ and $\delta$ are $\ll \tfrac{1}{\log X}$ and $X
\rightarrow \infty$. We consider, from the first part of the
functional equation for each of the $L$-functions,
\begin{eqnarray}
\sum_{d\le X}
\sum_{h,k,m,n}\frac{\mu(h)\mu(k)\chi_d(hkmn)}{h^{1/2+\gamma}k^{1/2+\delta}
m^{1/2+\alpha}n^{1/2+\beta}}.
\end{eqnarray}
Retaining only the terms for which $hkmn$ is square, leads us to
\begin{eqnarray} X^*\sum_{hkmn=\square}
\frac{\mu(h)\mu(k)a(hkmn)}{h^{1/2+\gamma}k^{1/2+\delta}
m^{1/2+\alpha}n^{1/2+\beta}}.
\end{eqnarray}
We express this sum as an Euler product (to ``save'' variables we now replace $h$ by $p^h$, etc.)
\begin{eqnarray}
\prod_p\sum_{{h+k+m+n \atop \begin{tiny}
\mbox{even}\end{tiny}}}\frac{\mu(p^h)\mu(p^k)a(p^{h+k+m+n})}{p^{h(1/2+\gamma)+
k(1/2+\delta)+m(1/2+\alpha)+n(1/2+\beta)}}.
\end{eqnarray}
We analyze the inner sum by dividing it into the four cases
according to $h=0,1$ and $k=0,1$; also it is helpful to note that
$\sum_{{m+n \atop \begin{tiny}
\mbox{even}\end{tiny}}}x^my^n=\frac{1+xy}{(1-x^2)(1-y^2)}$ and
$\sum_{{m+n \atop \begin{tiny}
\mbox{odd}\end{tiny}}}x^my^n=\frac{x+y}{(1-x^2)(1-y^2)}$. It is
more complicated to write down the exact formula for this,
complete with the arithmetic factor
$A_D(\alpha,\beta;\gamma,\delta)$. This factor is asymptotically 1
for small values of the parameters. Since we are interested in the
first main terms here we record that the relevant zeta factors in
the expression above are
\begin{eqnarray}
&&\frac{\zeta(1+2\alpha)\zeta(1+2\beta)\zeta(1+\alpha+\beta)\zeta(1+\gamma+\delta)}{\zeta(1+\alpha+\gamma)
\zeta(1+\alpha+\delta)\zeta(1+\beta+\gamma)\zeta(1+\beta+\delta)}\nonumber
\\
&&\qquad\qquad=\frac{(\alpha+\gamma)(\alpha+\delta)
(\beta+\gamma)(\beta+\delta)}{4\alpha\beta(\alpha+\beta)(\gamma+\delta)}+O(1/\log
X).
\end{eqnarray}
Thus we have, from the remaining parts of the functional equation:
\begin{conjecture}[{Conrey, Farmer and Zirnbauer [11]}]\label{conj:symp2} With
 $\alpha, \beta, \gamma, \delta \ll \tfrac{1}{\log X}$, we have
\begin{eqnarray} \label{eqn:R2symp} \nonumber
&&\frac{1}{X^*}R_D(\alpha,\beta;\gamma,\delta)=
\frac{(\alpha+\gamma)(\alpha+\delta)
(\beta+\gamma)(\beta+\delta)}{4\alpha\beta(\alpha+\beta)(\gamma+\delta)}-X^{-\alpha}
\frac{(-\alpha+\gamma)(-\alpha+\delta)
(\beta+\gamma)(\beta+\delta)}{4\alpha\beta(-\alpha+\beta)(\gamma+\delta)}\\
&&\qquad -X^{-\beta}\frac{(\alpha+\gamma)(\alpha+\delta)
(-\beta+\gamma)(-\beta+\delta)}{4\alpha\beta(\alpha-\beta)(\gamma+\delta)}
-X^{-\alpha-\beta}\frac{(-\alpha+\gamma)(-\alpha+\delta)
(-\beta+\gamma)(-\beta+\delta)}{4\alpha\beta(\alpha+\beta)(\gamma+\delta)}
\\
&&\qquad\qquad\qquad+O(1/\log X),\nonumber
\end{eqnarray}
as $X\rightarrow \infty$.
\end{conjecture}

\subsection{Orthogonal examples} As a third example, we consider the orthogonal
family of quadratic twists of the $L$-function $L_\Delta$ associated with the unique
weight 12 cusp form for the full modular group:
\begin{eqnarray}
L_\Delta(s,\chi_d)=\sum_{n=1}^\infty \frac{\chi_d(n)\tau^*(n)
}{n^s}=\prod _p \left(1-
\frac{\tau^*(p)\chi_d(p)}{p^s}+\frac{\chi_d(p^2)}{p^{2s}}\right)^{-1},
\end{eqnarray}
where $\tau^*(n)=\tau(n)/n^{11/2}$ and $\tau(n)$ is Ramanujan's tau-function.
For $d>0$, this has functional equation
\begin{eqnarray}
\xi_\Delta(s,\chi_d):=\left(\frac d{2\pi}\right)^s\Gamma(s+11/2)L_\Delta(s,\chi_d)=\xi_\Delta(1-s,\chi_d).
\end{eqnarray}
Let
\begin{eqnarray}
R_\Delta(\alpha;\gamma):= \sum_{d\le X} \frac{L_\Delta(1/2+\alpha,\chi_d)}{L_\Delta(1/2+\gamma,\chi_d)}
\end{eqnarray}
and let
\begin{eqnarray}
R_\Delta(\alpha,\beta;\gamma,\delta):= \sum_{d\le X}
\frac{L_\Delta(1/2+\alpha,\chi_d)L_\Delta(1/2+\beta,\chi_d)}{L_\Delta(1/2+\gamma,\chi_d)
L_\Delta(1/2+\delta,\chi_d)}.
\end{eqnarray}

As in the symplectic example we will calculate the full expression for $R_\Delta(\alpha;\gamma)$ and only the leading main terms for $R_\Delta(\alpha, \beta;\gamma, \delta)$.

Note that
\begin{eqnarray}
\frac{1}{L_\Delta(s,\chi_d)}=\prod_p\left(1-\frac{\tau^*(p)\chi_d(p)}{p^s}+\frac{\chi_d(p^2)}{p^{2s}}\right)=:
\sum_{n=1}^\infty\frac{\mu_\Delta(n)\chi_d(n)}{n^s}.
\end{eqnarray}
 To commence the calculation of $R_\Delta(\alpha;\gamma)$ we
replace each $L$-function in the numerator by the first half of
the approximate functional equation
\begin{equation}
L_{\Delta} (1/2+\alpha,\chi_d)=\sum_{m<x}
\frac{\chi_d(m)\tau^*(m)}{m^s}+\left(\frac{d}{2\pi}\right)^{-2\alpha}
\frac{\Gamma(6-\alpha)}{\Gamma(6+\alpha)} \sum_{n<y}
\frac{\chi_d(n)\tau^*(n)}{n^{1-s}} +\text{remainder},
\end{equation}
where $xy=d^2/(2\pi)$.
 We must then consider
\begin{eqnarray}
\sum_{d\le X}
\sum_{h,m}\frac{\mu_\Delta(h)\tau^*(m)\chi_d(hm)}{h^{1/2+\gamma}k^{1/2+\alpha}},
\end{eqnarray}
which leads to \begin{eqnarray}
&&X^*\sum_{hm=\square}\frac{\mu_\Delta(h)\tau^*(m)a(hm)}{h^{1/2+\gamma}k^{1/2+\alpha}}\\
&&\qquad\qquad=X^*\prod_p\bigg(1+\frac{p}{p+1}\sum_{h+m>0\atop
\begin{tiny}\mbox{even}\end{tiny}
}\frac{\mu_\Delta(p^h)\tau^*(p^m)}{p^{h(1/2+\gamma)+m(1/2+\alpha)}}\bigg).\nonumber
\end{eqnarray}
We note that $\mu_\Delta(p)=-\tau^*(p)$, $\mu_\Delta(p^2)=1$, and
$\mu_\Delta(p^m)=0$ for $m>2$, so that the product over primes
here is
\begin{eqnarray}\label{eqn:ep}
&&\qquad\prod_p\bigg(1+\frac{p}{p+1}\bigg(\sum_{m=1}^\infty\frac{\tau^*(p^{2m})}{p^{m(1+2\alpha)}}-
\frac{\tau^*(p)}{p^{1+\alpha+\gamma}}\sum_{m=0}^\infty\frac{\tau^*(p^{2m+1})}{p^{m(1+2\alpha)}}
+\frac{1}{p^{1+2\gamma}}\sum_{m=0}^\infty\frac{\tau^*(p^{2m})}{p^{m(1+2\alpha)}}\bigg)\bigg).
\end{eqnarray}
We note that
\begin{eqnarray}
\sum_{m=0}^\infty\tau^*(p^{2m})x^{2m}=\frac12\left\{\bigg(1-\tau^*(p)x+x^2\bigg)^{-1}
+\bigg(1+\tau^*(p)x+x^2\bigg)^{-1}\right\}
\end{eqnarray}
and
\begin{eqnarray}
\sum_{m=0}^\infty\tau^*(p^{2m+1})x^{2m+1}=\frac12\left\{\bigg(1-\tau^*(p)x+x^2\bigg)^{-1}
-\bigg(1+\tau^*(p)x+x^2\bigg)^{-1}\right\}.
\end{eqnarray}

 The ``polar'' part of the product (\ref{eqn:ep}) is
$\frac{\zeta(1+2\gamma)}{\zeta(1+\alpha+\gamma)};$  we can factor
these terms out and be left with a convergent Euler product.
However, we prefer at this point to factor out some other
$L$-functions present here with values near the 1-line and to be
left with an Euler product which is more rapidly convergent. To
this end, we recall the Rankin-Selberg convolution of $L_\Delta$
and the symmetric square $L$-function associated with $L_\Delta$.
We can write the Euler product for $L_\Delta$ as
\begin{eqnarray}
L_\Delta(s)=\prod_p\left(1-\frac{\alpha_p}{p^s}\right)^{-1}
\left(1-\frac{\overline{\alpha_p}}{p^s}\right)^{-1},
\end{eqnarray}
where $\alpha_p+\overline{\alpha_p}=\tau^*(p)$ and
$\alpha_p\overline{\alpha_p}=|\alpha_p|^2=1$. The Rankin-Selberg
$L$-function is
\begin{eqnarray}
L(\tau\otimes\tau,s)=\sum_{n=1}^\infty\frac{\tau^*(n)^2}{n^s}=\zeta(s)L_\Delta(\mbox{\rm{sym}}^2,s)\zeta(2s)^{-1},
\end{eqnarray}
where the symmetric square $L$-function is given by
\begin{eqnarray}
L_\Delta(\mbox{\rm{sym}}^2,s)=\prod_p\left(1-\frac{\alpha^2_p}{p^s}\right)^{-1}\left(1-\frac{1}{p^s}\right)^{-1}
\left(1-\frac{\overline{\alpha^2_p}}{p^s}\right)^{-1}
\end{eqnarray}
and is an entire function of $s$. As a Dirichlet series, we can
write
\begin{eqnarray}
L_\Delta(\mbox{\rm{sym}}^2,s)=\zeta(2s)^{-1}\sum_{n=1}^\infty
\frac{\tau^*(n^2)}{n^s}.
\end{eqnarray}
Thus, the product (\ref{eqn:ep}) can be expressed as
\begin{eqnarray}
\frac{\zeta(1+2\gamma)
L_\Delta(\mbox{\rm{sym}}^2,1+2\alpha)}{\zeta(1+\alpha+\gamma)L_\Delta(\mbox{\rm{sym}}^2,1+\alpha+\gamma)}B_\Delta(\alpha;\gamma),
\end{eqnarray}
where
\begin{eqnarray}
\label{eq:BDelta} B_\Delta(\alpha;\gamma)&=&\prod_p \nonumber
\label{eqn:B}
\bigg(1+\frac{p}{p+1}\bigg(\sum_{m=1}^\infty\frac{\tau^*(p^{2m})}{p^{m(1+2\alpha)}}-
\frac{\tau^*(p)}{p^{1+\alpha+\gamma}}\sum_{m=0}^\infty\frac{\tau^*(p^{2m+1})}{p^{m(1+2\alpha)}}
+\frac{1}{p^{1+2\gamma}}\sum_{m=0}^\infty\frac{\tau^*(p^{2m})}{p^{m(1+2\alpha)}}\bigg)\bigg)\\
&& \qquad \times
\frac{\left(1-\frac{\tau^*(p^2)}{p^{1+2\alpha}}+\frac{\tau^*(p^2)}{p^{2+4\alpha}}
-\frac{1}{p^{3+6\alpha}}\right)\left(1-\frac{1}{p^{1+2\gamma}}\right)}{
\left(1-\frac{\tau^*(p^2)}{p^{1+\alpha+\gamma}}+\frac{\tau^*(p^2)}{p^{2+2\alpha+2\gamma}}
-\frac{1}{p^{3+3\alpha+3\gamma}}\right)\left(1-\frac{1}{p^{1+\alpha+\gamma}}\right)}.
\end{eqnarray}
Note that $B(r;r)=1$; this follows from the fact that
$\tau^*(p^{2m+1})\tau^*(p)=\tau^*(p^{2m+2})+\tau^*(p^{2m}).$ Thus,
the ratios conjecture gives:
\begin{conjecture}[{Conrey, Farmer and Zirnbauer [11]}] \label{conj:orth} With
constraints on $\alpha$ and  $\gamma$ as described  at
(\ref{eq:constraints}), we have
\begin{eqnarray}&&\label{eqn:ratorth}
R_\Delta(\alpha;\gamma)=\sum_{d\le
X}\frac{L_\Delta(1/2+\alpha,\chi_d)}{L_\Delta(1/2+\gamma,\chi_d)}
=\sum_{d\le X} \bigg( \frac{\zeta(1+2\gamma)
L_\Delta(\mbox{\rm{sym}}^2,1+2\alpha)}{\zeta(1+\alpha+\gamma)L_\Delta(\mbox{\rm{sym}}^2,1+\alpha+\gamma)}B_\Delta(\alpha;\gamma)\nonumber\\
&&\qquad  + \left(\frac d{2\pi}\right)^{-2\alpha}
\frac{\Gamma(6-\alpha)}{\Gamma(6+\alpha)}\frac{\zeta(1+2\gamma)
L_\Delta(\mbox{\rm{sym}}^2,1-2\alpha)}{\zeta(1-\alpha+\gamma)L_\Delta(\mbox{\rm{sym}}^2,1-\alpha+\gamma)}B_\Delta(-\alpha;\gamma)
\bigg) +O(X^{1/2+\epsilon}),
\end{eqnarray}
where $B_{\Delta}(\alpha,\gamma)$ is defined in (\ref{eq:BDelta}).
\end{conjecture}

For application to the one-level density, we note that
\begin{eqnarray}
\sum_{d\le X}
\frac{L_\Delta'(1/2+r,\chi_d)}{L_\Delta(1/2+r,\chi_d)}
=\frac{d}{d\alpha}R_\Delta(\alpha;\gamma)\bigg|_{\alpha=\gamma=r}.
\end{eqnarray}
Now
\begin{eqnarray}&& \frac{d}{d\alpha}
\frac{\zeta(1+2\gamma)
L_\Delta(\mbox{\rm{sym}}^2,1+2\alpha)}{\zeta(1+\alpha+\gamma)L_\Delta(\mbox{\rm{sym}}^2,1+\alpha+\gamma)}B_\Delta(\alpha;\gamma)
\bigg|_{\alpha=\gamma=r}\nonumber
\\&&\qquad\qquad=-\frac{\zeta'(1+2r)}{\zeta(1+2r)}
+\frac{L_\Delta'(\mbox{\rm{sym}}^2,1+2r)}{L_\Delta(\mbox{\rm{sym}}^2,1+2
r)} +B_\Delta'(r;r)
\end{eqnarray}
 and
\begin{eqnarray}&&
 \frac{d}{d\alpha}
\left(\frac d{2\pi}\right)^{-2\alpha}
\frac{\Gamma(6-\alpha)}{\Gamma(6+\alpha)}\frac{\zeta(1+2\gamma)
L_\Delta(\mbox{\rm{sym}}^2,1-2\alpha)}{\zeta(1-\alpha+\gamma)L_\Delta(\mbox{\rm{sym}}^2,1-\alpha+\gamma)}B_\Delta(-\alpha;\gamma)
  \bigg|_{\alpha=\gamma=r}\\
&&\qquad =-\left(\frac d{2\pi}\right)^{-2r}
\frac{\Gamma(6-r)}{\Gamma(6+r)}\frac{\zeta(1+2r)
L_\Delta(\mbox{\rm{sym}}^2,1-2r)}{ L_\Delta(\mbox{\rm{sym}}^2,1
)}B_\Delta(-r;r).\nonumber
 \end{eqnarray}
Thus, the ratios conjecture implies:
\begin{theorem} \label{theo:orthderiv}
Assuming Conjecture \ref{conj:orth}, if $\frac{1}{\log X} \ll \Re
r<\frac{1}{4}$ and $\Im r\ll_{\epsilon} X^{1-\epsilon}$, then
 \begin{eqnarray} \label{eqn:orthlev1}
&&\sum_{d\le X} \nonumber \frac{L_\Delta'(1/2+r,\chi_d)}{L_\Delta(1/2+r,\chi_d)}
= \sum_{d\le X}\bigg(-\frac{\zeta'(1+2r)}{\zeta(1+2r)} +\frac{L_\Delta'(\mbox{\rm{sym}}^2,1+2r)}{L_\Delta(\mbox{\rm{sym}}^2,1+2 r)} +B_\Delta'(r;r)\\
&& \qquad -\left(\frac d{2\pi}\right)^{-2r}
\frac{\Gamma(6-r)}{\Gamma(6+r)}\frac{\zeta(1+2r)
L_\Delta(\mbox{\rm{sym}}^2,1-2r)}{ L_\Delta(\mbox{\rm{sym}}^2,1
)}B_\Delta(-r;r)\bigg)+O(X^{1/2+\epsilon})
\end{eqnarray}
where $B_{\Delta}(\alpha,\gamma)$ is defined in (\ref{eq:BDelta}).
\end{theorem}

We now determine the main terms when
$\alpha,\beta,\gamma,\delta\ll \tfrac{1}{\log X}$ and $X\to
\infty$ for the average over this family of the ratio
$R_\Delta(\alpha,\beta;\gamma,\delta)$ of two $L$-functions over
two $L$-functions. We are quickly led to consider
\begin{eqnarray}
\sum_{hkmn=\square}\frac{\mu_\Delta(h)\mu_\Delta(k)
\tau^*(m)\tau^*(n)a(hkmn)}
{h^{1/2+\gamma}k^{1/2+\delta}m^{1/2+\alpha}n^{1/2+\beta}}.
\end{eqnarray}
When we go to Euler products, we find that this expression
evaluates to
\begin{eqnarray}
\frac{\zeta(1+\alpha+\beta)\zeta(1+2\gamma)\zeta(1+\gamma+\delta)\zeta(1+2\delta)}
{\zeta(1+\alpha+\gamma)\zeta(1+\alpha+\delta)\zeta(1+\beta+\gamma)\zeta(1+\beta+\delta)}
A_\Delta(\alpha,\beta;\gamma,\delta),
\end{eqnarray}
where $A$ is analytic if the real parts of
$\alpha,\beta,\gamma,\delta$ are smaller than 1/4 in absolute
value; moreover $A_{\Delta}(0,0;0,0)=1$. Thus, this part is
\begin{eqnarray}
=
\frac{(\alpha+\gamma)(\alpha+\delta)(\beta+\gamma)(\beta+\delta)}
{(\alpha+\beta)(2\gamma)(\gamma+\delta)(2\delta)}+O(1/\log X) .
\end{eqnarray}
Taking the symmetric  sum of four of these terms, arising from the
product of the approximate functional equations of the two
$L$-functions in the numerator, we find that
\begin{conjecture}[{Conrey, Farmer and Zirbauer [11]}] \label{conj:orth2} With
 $\alpha, \beta, \gamma, \delta \ll \tfrac{1}{\log X}$, we have
\begin{eqnarray} \label{eqn:ratorth2}
\frac{1}{X^*}R_\Delta(\alpha,\beta;\gamma,\delta) &=&
\frac{1}{X^*}\sum_{d\le X}
\frac{L_\Delta(1/2+\alpha,\chi_d)L_\Delta(1/2+\beta,\chi_d)}{L_\Delta(1/2+\gamma,\chi_d)
L_\Delta(1/2+\delta,\chi_d)}\\&=&\nonumber
\frac{(\alpha+\gamma)(\alpha+\delta)(\beta+\gamma)(\beta+\delta)}{(\alpha+\beta)
(2\gamma)(\gamma+\delta)(2\delta)}
\\&& \qquad + X^{-2\alpha}
\frac{(-\alpha+\gamma)(-\alpha+\delta)(\beta+\gamma)(\beta+\delta)}
{(-\alpha+\beta)(2\gamma)(\gamma+\delta)(2\delta)}\nonumber\\
&& \qquad \qquad +X^{-2\beta}\nonumber
\frac{(\alpha+\gamma)(\alpha+\delta)(-\beta+\gamma)(-\beta+\delta)}{(\alpha-\beta)(2\gamma)(\gamma+\delta)(2\delta)}
\\ && \qquad \qquad \qquad
- X^{-2\alpha-2\beta}\nonumber
\frac{(-\alpha+\gamma)(-\alpha+\delta)(-\beta+\gamma)(-\beta+\delta)}{(\alpha+\beta)(2\gamma)
(\gamma+\delta)(2\delta)}\nonumber \\
&& \qquad\qquad\qquad\qquad+O(1/\log X),\nonumber
\end{eqnarray}
as $X\to \infty$.
\end{conjecture}

\section{One-level density}\label{sect:onelevel}
In this section we use the ratios conjecture to compute the
one-level density function for zeros of quadratic Dirichlet
$L$-functions, complete with lower order terms.  \"Ozl\"uk and
Snyder [38] have proven such results (assuming the generalized
Riemann Hypothesis) for test functions $f$ for which the support
of $\hat{f}$ is limited. The ratios conjectures imply a result
consistent with [38] but with no constraint on the support of the
Fourier transform of the test function.

For simplicity, we assume that:
\begin{eqnarray} \label{eq:fconditions}
&&f(z) \text{ is holomorphic throughout the strip }  |\Im z| <2,\nonumber\\
&& \text{is real on the real line and even,}\\
&& \text{and that } f(x)\ll 1/(1+x^2) \text{ as } x\to
\infty.\nonumber
\end{eqnarray}

 We consider
\begin{eqnarray}
S_1(f):=\sum_{d\le X} \sum_{\gamma_d} f(\gamma_d),
\end{eqnarray}
where $\gamma_d$ denotes the ordinate of a generic zero of $L(s,\chi_d)$ on the half-line (we are assuming that all of the complex zeros are on the 1/2-line).

We have
\begin{eqnarray}
S_1(f)=\sum_{d\le X} \frac{1}{2\pi i}
\left(\int_{(c)}-\int_{(1-c)}\right)
\frac{L'(s,\chi_d)}{L(s,\chi_d)} f(-i(s-1/2))~ds,
\end{eqnarray}
where  $(c)$ denotes a vertical line from $c-i \infty$ to
$c+i\infty$ and $3/4>c>1/2+\tfrac{1}{\log X}$. The integral on the
$c$-line is
\begin{eqnarray}
\label{eq:density1} \frac{1}{2\pi  }\int_{-\infty}^\infty
f(t-i(c-1/2)) \sum_{d\le X}
\frac{L'(1/2+(c-1/2+it),\chi_d)}{L(1/2+(c-1/2+it),\chi_d)} ~dt.
\end{eqnarray}
It follows by the Riemann Hypothesis that on the path of
integration $(c)$
\begin{equation}
\label{eq:logderivLbound} \frac{L'(s,\chi_d)}{L(s,\chi_d)}\ll
\log^2(|s|d).
\end{equation}
For $|t|>X^{1-\epsilon}$ we estimate the integral using
(\ref{eq:logderivLbound}) and (\ref{eq:fconditions}) and the
result is $\ll X^{\epsilon}$.  By the ratios conjecture
(\ref{eqn:onelevel}), if $|t|<X^{1-\epsilon}$ the sum over $d$ in
(\ref{eq:density1}) is:
\begin{eqnarray}\label{eq:density2}&&
 \sum_{d\le X} \bigg(\frac{\zeta'(1+2r)}{\zeta(1+2r)} +A_D'(r;r)-\left(\frac d \pi\right)^{-r}
\frac{\Gamma(1/4-r/2)}{\Gamma(1/4+r/2)}
 \zeta(1-2r)  A_D(-r;r)\bigg)\bigg|_{r=c-1/2+it}\\
 &&\qquad \qquad +O(X^{1/2+\epsilon}).\nonumber
 \end{eqnarray}
Since the quantity in (\ref{eq:density2}) is $\ll X^{1+\epsilon}$
for $|t|<X^{1-\epsilon}$ and $f(t)\ll \tfrac{1}{t^2}$ we can
extend the integration in $t$ to infinity.
 Finally, since the integrand is regular at $r=0$, we can move the path of
integration to $c=1/2$ and so obtain
\begin{eqnarray}
&&\frac{1}{2\pi}\int_{-\infty}^\infty f(t) \sum_{d\le X} \bigg(\frac{\zeta'(1+2it)}{\zeta(1+2it)} +A_D'(it;it)-\left(\frac d \pi\right)^{-it}
\frac{\Gamma(1/4-it/2)}{\Gamma(1/4+it/2)}
 \zeta(1-2it)  A_D(-it;it)\bigg) ~dt\\
 &&\qquad \qquad +O(X^{1/2+\epsilon}).\nonumber
\end{eqnarray}
For the integral on the $1-c$ line, we change variables, letting  $s\to 1-s$, and we use the functional equation
\begin{eqnarray}
\frac{L'(1-s,\chi_d)}{L(1-s,\chi_d)}=\frac{X'(s,\chi_d)}{X(s,\chi_d)}-
\frac{L'(s,\chi_d)}{L(s,\chi_d)}
\end{eqnarray}
where
\begin{eqnarray}
\frac{X'(s,\chi_d)}{X(s,\chi_d)}=-\log \frac{d}{\pi} -\frac12 \frac {\Gamma'}{\Gamma}\left(\frac{1-s}{2}\right)-\frac12 \frac {\Gamma'}{\Gamma}\left(\frac{ s}{2}\right)
\end{eqnarray}
The contribution from the $L'/L$ term is now exactly as before, since $f$ is even.  Thus, we obtain
\begin{theorem} \label{theo:1level} Assuming Conjecture \ref{conj:symp} and $f$ satisfying (\ref{eq:fconditions}), we have
\begin{eqnarray}
&&\sum_{d\le X} \sum_{\gamma_d} f(\gamma_d)= \frac{1}{2\pi} \int_{-\infty}^\infty f(t)\sum_{d\le X}  \bigg( \log \frac{d}{\pi} +\frac12 \frac {\Gamma'}{\Gamma}(1/4+it/2)+\frac12 \frac {\Gamma'}{\Gamma}(1/4-it/2)+\\
&& \qquad \qquad 2\bigg(\frac{\zeta'(1+2it)}{\zeta(1+2it)} +A_D'(it;it)-\left(\frac d \pi\right)^{-it}
\frac{\Gamma(1/4-it/2)}{\Gamma(1/4+it/2)}
 \zeta(1-2it)  A_D(-it;it)\bigg)\bigg) ~dt\nonumber\\
 &&\qquad \qquad \qquad +O(X^{1/2+\epsilon}),\nonumber
\end{eqnarray}
where
\begin{eqnarray}
A_D(-r;r)=
\prod_p\left(1-\frac{1}{(p+1)p^{1-2r}}-\frac{1}{p+1}\right)\left(1-\frac{1}{p}\right)^{-1},
\end{eqnarray}
and
\begin{eqnarray}
A'_D(r;r)=\sum_p \frac{\log p}{(p+1)(p^{1+2r}-1)}.
\end{eqnarray}
\end{theorem}

The low-lying zeros of this family of $L$-functions are expected
to display the same statistics as the eigenvalues of the matrices
from $USp(2N)$ chosen with respect to Haar measure.  Thus in the
large $X$ limit, the one level density of the scaled zeros will
have the form, as proved by \"Ozl\"uk   and Snyder [38],
\begin{eqnarray}
&&\lim_{X\rightarrow\infty} \frac{1}{X^*} \sum_{d\leq X}
\sum_{\gamma_d} f\left(\gamma_d \frac{\log \frac{d}{\pi}} {2\pi}
\right) = \int_{-\infty}^{\infty} f(x) \bigg(1-\frac{\sin (2\pi
x)}{2\pi x}\bigg) dx,
\end{eqnarray}
where $X^*$ is the number of terms in the sum (and is proportional
to $X$).

Defining $f(t)=g(\tfrac{t\log X}{2\pi})$ and scaling the variable
$t$ from Theorem \ref{theo:1level} as $\tau=\frac{t\log X}{2\pi}$
\begin{eqnarray}
&&\sum_{d\leq X} \sum _{\gamma_d}g\left(\frac{\gamma_d\log
X}{2\pi}\right)=\frac{1}{\log X} \int_{-\infty}^\infty
g(\tau)\sum_{d\le X} \bigg( \log \frac{d}{\pi} +\frac12 \frac
{\Gamma'}{\Gamma}\left(1/4+\frac{i\pi\tau}{\log
X}\right)\\
&&\qquad+\frac12 \frac {\Gamma'}{\Gamma}\left(1/4-\frac{i\pi \tau
}{\log X}\right)+2\bigg(\frac{\zeta'(1+\frac{4i\pi\tau}{\log
X})}{\zeta(1+\frac{4i\pi\tau}{\log X})}+A_D'\left(\frac{2\pi
i\tau}{\log X};\frac{2\pi i\tau}{\log
X}\right)\nonumber\\
&& \qquad   -e^{-\frac{2\pi i\tau}{\log X}\log\frac{d}{\pi}}
\frac{\Gamma(1/4-\frac{i\pi \tau}{\log X})}{\Gamma(1/4+\frac{i\pi
\tau}{\log X})}
 \zeta\left(1-\frac{4\pi i\tau}{\log X}\right)  A_D\left(-\frac{2\pi i\tau}{\log X};\frac{2\pi i\tau}{\log X}\right)
 \bigg)\bigg) d\tau\nonumber\\
 &&\qquad \qquad \qquad +O(X^{1/2+\epsilon}).\nonumber
 \end{eqnarray}
 For large $X$ only the $\log\frac{d}{\pi}$ term, the $\frac{\zeta'}{\zeta}$ term and the final
 term in the integral contribute, yielding the asymptotic
 \begin{eqnarray}
&&\sum_{d\leq X} \sum _{\gamma_d}g\left(\frac{\gamma_d\log
X}{2\pi}\right)\sim\frac{1}{\log X} \int_{-\infty}^\infty
g(\tau)\bigg( X^* \log X-X^*\frac{\log X}{2\pi
i\tau}+X^*\frac{e^{-2\pi i\tau}}{2\pi i\tau}\log X \bigg)d\tau.
\end{eqnarray}
However, since $g$ is an even function,  the middle term above
drops out and the last term can be duplicated with a change of
sign of $\tau$, leaving
\begin{eqnarray}
&&\lim_{X\rightarrow\infty} \frac{1}{X^*} \sum_{d\leq X}
\sum_{\gamma_d} g\left(\gamma_d \frac{\log \frac{d}{\pi}} {2\pi}
\right) = \int_{-\infty}^{\infty} g(\tau)\bigg(1+\frac{e^{-2\pi
i\tau}}{4\pi i\tau}+\frac{e^{2\pi i\tau}}{-4\pi i\tau}\bigg)d\tau,
\end{eqnarray}
and resulting in exactly the answer expected.

In much the same way as for Theorem \ref{theo:1level} we can
compute the lower order terms in the one-level density for the
zeros of the functions from the orthogonal family
$L_\Delta(s,\chi_d)$ by using (\ref{eqn:orthlev1}).
\begin{theorem}
Assuming Conjecture \ref{conj:orth} and with $f$ satisfying
(\ref{eq:fconditions}), we have
\begin{eqnarray}
&&\sum_{d\le X} \sum_{\gamma_{\Delta,d}} f(\gamma_{\Delta,d})= \frac{1}{2\pi}
\int_{-\infty}^\infty f(t)\sum_{d\le X} \bigg(
 2\log \tfrac{d}{2\pi} +\frac{\Gamma'}{\Gamma}(6+it)+\frac{\Gamma'}{\Gamma}(6-it)\\
&&\qquad +2\bigg(-\frac{\zeta'(1+2it)}{\zeta(1+2it)}
+\frac{L_\Delta'(\mbox{\rm{sym}}^2,1+2it)}
{L_\Delta(\mbox{\rm{sym}}^2,1+2 it)} +B_\Delta'(it;it)\nonumber\\
&& \qquad -\left(\frac d{2\pi}\right)^{-2it}
\frac{\Gamma(6-it)}{\Gamma(6+it)}\frac{\zeta(1+2it)
L_\Delta(\mbox{\rm{sym}}^2,1-2it)}{ L_\Delta(\mbox{\rm{sym}}^2,1 )}B_\Delta(-it;it)\bigg)~dt\nonumber\\
 &&\qquad \qquad \qquad +O(X^{1/2+\epsilon})\nonumber
\end{eqnarray}
where $B_{\Delta}$ is defined in (\ref{eqn:B}).
\end{theorem}

 In the same way as above, the main terms here give the
one-level density of eigenvalues of matrices from the group
$SO(2N)$, which in the limit of large $N$ is $1+\tfrac{\sin 2\pi
x}{2\pi x}$.

\section{Pair-correlation}
\label{sect:pair} We show how to use the ratios conjecture to
compute the pair-correlation of the zeros of the Riemann
zeta-function, originally conjectured by Montgomery [36], together
with lower order (arithmetic) terms that have been found
heuristically  by Bogolmony and Keating [3] (see [28] for a more
expository description, [2] for numerical calculation of these
lower order terms and [19] for related rigorous results). When
Farmer formulated his original ratio conjecture (\ref{eq:farmer})
he observed in [16] that it implied the leading order terms of
Montgomery's pair correlation conjecture.  Farmer's method is
completely different from what we present below.

We want to evaluate the sum
\begin{eqnarray}
S(f)=\sum_{0<\gamma, \gamma'<T} f(\gamma-\gamma')
\end{eqnarray}
for a test function $f$ satisfying (\ref{eq:fconditions}).
 We rewrite the sum in question in terms of contour
integrals. Let $1/2+\tfrac{1}{\log T}< a < b<3/4$ and let
$\mathcal C_1$ be the positively oriented rectangular contour with
corners $a, a+iT, 1-a+i T, 1-a$ and let $\mathcal C_2$ be the
rectangular contour with corners $b, b+iT, 1-b+i T, 1-b$. Then
\begin{eqnarray}
S(f)=\frac{1}{(2\pi i)^2} \int_{\mathcal C_1}\int _{\mathcal C_2} \frac{\zeta'}{\zeta}(z)\frac{\zeta'}{\zeta}(w) f(-i(z-w)) ~dw ~dz;
\end{eqnarray}
the point, of course, is that the poles inside the contours are
simple poles with residue 1 at the zeros $z=1/2+i\gamma$ and
$w=1/2+i\gamma'$ of the zeta-function. The integrals along the
horizontal sides are small and may be ignored. Thus, we consider 4
double integrals. We consider each of the 4 double integrals
separately; call them $I_1,\dots, I_4$, where $I_1$ has vertical
parts $a$ and $b$, $I_2$ has vertical parts $1-a$ and $1-b$, $I_3$
has vertical parts $a$ and $1-b$ and $I_4$ has vertical parts
$1-a$ and $b$.

It is easy to see using the Riemann Hypothesis that
$I_1=O(T^\epsilon)$ just by moving the contours to the right of 1
and integrating term-by-term.

For $I_2$, we use the functional equation
$\frac{\zeta'}{\zeta}(s)=\frac{\chi'}{\chi}(s)-\frac{\zeta'}{\zeta}(1-s)$
for $s=w$ and $s=z$  and find similarly that
\begin{eqnarray}
I_2&=&\frac{1}{(2\pi)^2}
\int_{0}^T\int_{0}^T\frac{\chi'}{\chi}(1/2+iu)\frac{\chi'}{\chi}(1/2+iv)f(u-v)
~du ~dv+O(T^{\epsilon}). \end{eqnarray} Using the fact that
\begin{eqnarray}
\frac{\chi'}{\chi}(1/2+it)=-\log\frac{|t|}{2\pi}\left(1+O\left(\frac
1 {|t|}\right)\right)
\end{eqnarray}
and that $f$ is even, we see, after the substitution $u=v+\eta$, that
\begin{eqnarray}
I_2&=&\frac{2}{(2\pi)^2}\int_{0}^T \int_v^T \log\frac{v}{2\pi}\log \frac{u}{2\pi}f(u-v)~du~dv+O(T^\epsilon)\\
&=&\nonumber \frac{2}{(2\pi)^2}\int_0^T f(\eta) \int_{0}^{T-\eta}
\log\frac{v}{2\pi}\log \frac{v+\eta}{2\pi} ~dv
~d\eta+O(T^{\epsilon}).
\end{eqnarray}
Recall that $f$ satisfies
\begin{equation}
f(x)\ll\frac{1}{1+x^2}
\end{equation}
for real $x$.  Letting $v\to vT$ in the inner integral above, we have
\begin{eqnarray}
I_2&=& \frac{2}{(2\pi)^2}T \int_0^T f(\eta) \int_0^{1-\frac{\eta}{T}}
 \log\frac{vT}{2\pi}\log \frac{vT+\eta}{2\pi} ~dv ~d\eta+O(T^{\epsilon}).
\end{eqnarray}
We may extend the upper limit  of the inner integral to $v=1$,
introducing an error term of size $\ll \int_0^T \eta f(\eta)
\log^2T d\eta\ll \log^3 T$. We can also replace $\log (vT+\eta)$
by $\log vT$ with the same error term. Thus,
\begin{eqnarray}
I_2&=& \frac{2}{(2\pi)^2}T \int_0^T f(\eta) \int_0^{1 }
 \log^2\frac{vT}{2\pi}  ~dv ~d\eta+O(T^{\epsilon})\\
&=&\nonumber
 \frac{1}{(2\pi)^2}\int_{-T}^Tf(\eta) \int_0^T \log^2 \frac{v}{2\pi} ~dv~d\eta +O(T^\epsilon).
\end{eqnarray}

Next we consider $I_3$. Letting $z=w+\eta$, it is
\begin{eqnarray}
I_3&=&\frac{-1}{(2\pi i)^2}\int_{1-b}^{1-b+iT}\int_{a}^{a+iT}
\frac{\zeta'}{\zeta}(w)\frac{\zeta'}{\zeta}(z)f(-i(z-w))~dw ~dz\\
&=& \frac{-1}{(2\pi )^2i}\int_{1-a-b-iT}^{1-a-b+iT} f(-i\eta)
\int_{T_1}^{T_2} \frac{\zeta'}{\zeta}(a+i
t)\frac{\zeta'}{\zeta}(a+it +\eta)~dt ~d\eta,\nonumber
\end{eqnarray}
where $T_1=\max\{0,-\Im\eta\}$ and $T_2=\min\{T-\Im\eta,T\}.$
We use the functional equation
\begin{eqnarray}\frac{\zeta'}{\zeta}(a+\eta+it)=\frac{\chi'}{\chi}(a+\eta+it)
-\frac{\zeta'}{\zeta}(1-a-\eta-it).
\end{eqnarray}
The term with the $\chi'/\chi$ is small as is seen by moving the contour to the right. Thus, we see that
\begin{eqnarray}
\qquad I_3&=& \frac{1}{(2\pi )^2i}\int_{1-a-b-iT}^{1-a-b+iT} f(-i\eta) \int_{T_1}^{T_2} \frac{\zeta'}{\zeta}(a+i t)\frac{\zeta'}{\zeta}(1-a-it -\eta)~dt ~d\eta+O(T^{\epsilon})\\
&=& \nonumber\frac{1}{(2\pi )^2i}\int_{1-a-b-iT}^{1-a-b+iT} f(-i\eta) \int_{T_1}^{T_2} \frac{\zeta'}{\zeta}(s+(a-1/2))\frac{\zeta'}{\zeta}(1-s+(1/2-a -\eta))~dt ~d\eta\\
&&\qquad +O(T^{\epsilon}),\nonumber
\end{eqnarray}
where $s=1/2+it$.  By Theorem \ref{theo:zetaderiv}, we have
\begin{eqnarray}
I_3 &=& \frac{1}{(2\pi )^2i}\int_{1-a-b-iT}^{1-a-b+iT} f(-i\eta)
\int_{T_1}^{T_2}
\left( \left(\frac{\zeta'}{\zeta}\right)'(1-\eta)+\right. \\
&&\quad \left(\frac{t}{2\pi}\right)^
{\eta} \zeta(1-\eta)\zeta(1+\eta)
\prod_p\frac{(1-\frac{1}{p^{1-\eta}})
(1-\frac 2 p +\frac{1}{p^{1-\eta}})}{(1-\frac 1 p )^2}\nonumber\\
&& \qquad \left. -\sum_p \left(\frac{\log p}{(p^{1-\eta}-1)}\right)^2\right)~dt~d\eta
 +O(T^{1/2+\epsilon}).\nonumber
\end{eqnarray}
Let $\delta=a+b-1$ and let $g(-\eta,t)$ be the integrand in the
second integral above. We can extend the range of the inner
integration, much as we did for the $I_2$ integral to the interval
$[0,T]$ with an error term of size $\ll T^\epsilon \int_\eta
|\eta| |f(\eta)|d\eta \ll T^\epsilon$. Thus, we obtain
\begin{eqnarray}
I_3=\frac{1}{(2\pi)^2i}\int_0^T\int_{-\delta-iT}^{-\delta+iT}f(-i\eta)g(-\eta,t)~d\eta
~dt+O(T^{1/2+\epsilon}).
\end{eqnarray}

Now we consider $I_4$.  Again letting $z=w+\eta$, we have
\begin{eqnarray}
I_4&=&\frac{1}{(2\pi i)^2}\int_{1-a}^{1-a+iT}\int_b^{b+iT}
\frac{\zeta'}{\zeta}(w)\frac{\zeta'}{\zeta}(z)f(-i(z-w))~dz ~dw\\
&=& \frac{1}{(2\pi )^2i}\int_{a+b-1-iT}^{a+b-1+iT} f(-i\eta)
\int_{T_1}^{T_2} \frac{\zeta'}{\zeta}(1-a+i
t)\frac{\zeta'}{\zeta}(1-a+it +\eta)~dt ~d\eta.\nonumber
\end{eqnarray}
We use the functional equation
\begin{equation}\frac{\zeta'}{\zeta}(1-a+it)=\frac{\chi'}{\chi}(1-a+it)-\frac{\zeta'}{\zeta}(a-it).
\end{equation}
Again, the contribution of the $\chi'/\chi$ term is negligible. Thus,
\begin{eqnarray}
\qquad I_4&=&
 \frac{1}{(2\pi )^2i}\int_{a+b-1-iT}^{a+b-1+iT} f(-i\eta) \int_{T_1}^{T_2} \frac{\zeta'}{\zeta}(a-i t)\frac{\zeta'}{\zeta}(1-a+it +\eta)~dt ~d\eta+O(T^{\epsilon})\\
&=& \frac{1}{(2\pi )^2i}\int_{a+b-1-iT}^{a+b-1+iT}
f(-i\eta)\int_{T_1}^{T_2}
  \frac{\zeta'}{\zeta}(1-s+(a-1/2))\frac{\zeta'}{\zeta}(s+(1/2-a +\eta))~dt ~d\eta \nonumber\\
&&\qquad +O(T^{\epsilon}).\nonumber
\end{eqnarray}
Now, by Theorem \ref{theo:zetaderiv},
\begin{eqnarray}
I_4&=&
 \frac{1}{(2\pi )^2i}\int_{a+b-1-iT}^{a+b-1+iT} f(-i\eta) \int_{T_1}^{T_2}
\left( \left(\frac{\zeta'}{\zeta}\right)'(1+\eta)+\right. \\
&&\quad \left(\frac{t}{2\pi}\right)^
{-\eta} \zeta(1-\eta)\zeta(1+\eta)
\prod_p\frac{(1-\frac{1}{p^{1+\eta}})
(1-\frac 2 p +\frac{1}{p^{1+\eta}})}{(1-\frac 1 p )^2}\nonumber\\
&& \qquad \left. -\sum_p \left(\frac{\log
p}{(p^{1+\eta}-1)}\right)^2\right)~dt~d\eta
 +O(T^{1/2+\epsilon}).\nonumber
\end{eqnarray}
Using the notation introduced after the calculation of $I_3$, and again extending the range of the integration in the inner integral, we can write the expression for $I_4$ as
\begin{eqnarray}
I_4=\frac{1}{(2\pi)^2i}\int_0^T\int_{\delta-i
T}^{\delta+iT}f(-i\eta)g(\eta,t)~d\eta ~dt+O(T^{1/2+\epsilon}).
\end{eqnarray}
Combining this with what we found for $I_3$ we have, after a change of variables,
\begin{eqnarray}
I_3+I_4=\frac{2}{(2\pi)^2i}\int_0^T\int_{\delta-iT}^{\delta+iT}f(i\eta)g(\eta,t)~d\eta
~dt+O(T^{1/2+\epsilon}).
\end{eqnarray}

Now let
\begin{eqnarray}
A(\eta)=\prod_p\frac{(1-\frac{1}{p^{1+\eta}})
(1-\frac 2 p +\frac{1}{p^{1+\eta}})}{(1-\frac 1 p )^2}
\end{eqnarray}
and
\begin{eqnarray}
B(\eta)=\sum_p \left(\frac{\log p}{(p^{1+\eta}-1)}\right)^2
\end{eqnarray}
so that
\begin{eqnarray}
g(\eta,t)=\left(\frac{\zeta'}{\zeta}\right)'(1+\eta)+
 \left(\frac{t}{2\pi}\right)^
{-\eta} \zeta(1-\eta)\zeta(1+\eta)A(\eta)-B(\eta).
\end{eqnarray}
Near 0, we see that (note that $A'(0)=0$),
\begin{eqnarray}
g(\eta,t)=\frac{ \log\frac{t}{2\pi}}{\eta}+O(1).
\end{eqnarray}
We move the path of integration in $\eta$ to   the imaginary axis
from $-T$ to $T$ with a principal value as we pass through 0; the contribution from $1/2$  of the residue from the pole of $g$ at $\eta=0$   is
\begin{equation}
\pi \int_0^T f(0) \log \frac{t}{2 \pi} ~dt.
\end{equation}

Combining our expressions for $I_1,\dots, I_4$, and changing $\eta$ into $ir$  we have
\begin{theorem} \label{theo:2point} Assuming  Conjecture
\ref{conj:unit}, and with $f$ satisfying (\ref{eq:fconditions}),
we have
\begin{eqnarray}
 &&\sum_{\gamma, \gamma'\le T} f(\gamma-\gamma') =\frac{1}{(2\pi)^2}\int_0^T\bigg(2\pi f(0) \log \frac{t}{2 \pi}
+ \int_{-T}^T f(r) \bigg( \log^2 \frac{t}{2\pi} +2\bigg(\left(\frac{\zeta'}{\zeta}\right)'(1+ir)\\
 && \qquad \qquad
+ \left(\frac{t}{2\pi}\right)^ {-ir}
\zeta(1-ir)\zeta(1+ir)A(ir)-B(ir)\bigg)\bigg) ~dr \bigg)~dt
+O(T^{1/2+\epsilon});\nonumber
\end{eqnarray}
here the integral is to be regarded as a principal value near $r=0$,
\begin{eqnarray}
A(\eta)=\prod_p\frac{(1-\frac{1}{p^{1+\eta}})
(1-\frac 2 p +\frac{1}{p^{1+\eta}})}{(1-\frac 1 p )^2},
\end{eqnarray}
and
\begin{eqnarray}
B(\eta)=\sum_p \left(\frac{\log p}{(p^{1+\eta}-1)}\right)^2.
\end{eqnarray}
\end{theorem}

We believe that this formula, originally found by Bogomolny and
Keating [3], is very accurate, indeed, down to a square root error
term. It includes all of the lower order terms that arise from
arithmetical considerations and should include all of the
fluctuations found in any of the extensive numerical experiments
that have been done. We have not scaled any of the terms here so
that terms of different scales are shown all at once.

To see the leading order term from Montgomery's pair-correlation
conjecture, let $L=\log \tfrac{T}{2\pi}$, $g(x
\frac{L}{2\pi})=f(x)$, and scale the variable $r$ in the inner
integral in Theorem \ref{theo:2point} as $y=r\frac{L}{2\pi}$:
\begin{eqnarray}
 &&\sum_{\gamma, \gamma'\le T} g((\gamma-\gamma')\tfrac{L}{2\pi})\\
  &&\quad=\frac{1}{(2\pi)^2}
 \int_0^T\bigg(2\pi g(0) \log \frac{t}{2 \pi}
+ \frac{2\pi}{L}\int_{-T\frac{L}{2\pi}}^{T\frac{L}{2\pi}}
 g(y) \bigg( \log^2 \frac{t}{2\pi} +2\bigg(\left(\frac{\zeta'}{\zeta}\right)'(1+\tfrac{2\pi i y}{L}
 )\nonumber\\
 &&\quad\quad
+ e^ {-2\pi i y \frac{\log(t/2\pi)}{L}} \zeta(1-\tfrac{2\pi i
y}{L})\zeta(1+\tfrac{2\pi i y}{L})A(\tfrac{2\pi i
y}{L})-B(\tfrac{2\pi i y}{L})\bigg)\bigg) ~dy \bigg)~dt
+O(T^{1/2+\epsilon}).\nonumber
\end{eqnarray}
For large $T$, only the $\log^2 \frac{t}{2\pi}$ and the two terms
containing zeta functions contribute, so we have the asymptotic
\begin{eqnarray}
 &&\sum_{\gamma, \gamma'\le T} g((\gamma-\gamma')\tfrac{L}{2\pi}) \sim \frac{1}{(2\pi)^2}
 \int_0^T\bigg(2\pi g(0) \log \frac{t}{2 \pi}\\
 &&\qquad\qquad
+ \frac{2\pi}{L}\int_{-T\frac{L}{2\pi}}^{T\frac{L}{2\pi}}
 g(y) \bigg( \log^2 \frac{t}{2\pi} -\frac{L^2}{2\pi^2 y^2}
+ e^ {-2\pi i y \frac{\log(t/2\pi)}{L}} \left(
\frac{L^2}{2\pi^2y^2}\right)\bigg)\bigg) ~dy \bigg)~dt.\nonumber
\end{eqnarray}
Integrating over $t$, we find that
\begin{eqnarray}
 \sum_{\gamma, \gamma'\le T} g((\gamma-\gamma')\tfrac{L}{2\pi})
 &\sim&
 \frac{T}{2\pi}\log\tfrac{T}{2\pi} \left(g(0)
+ \int_{-\infty}^{\infty}
 g(y) \left( 1 -\frac{1}{2\pi^2 y^2}
+    \frac{\cos(2\pi y)}{2\pi^2y^2}\right) ~dy \right)\\
&=&\frac{T}{2\pi}\log\tfrac{T}{2\pi} \left(g(0) +
\int_{-\infty}^{\infty}
 g(y) \bigg( 1 -\bigg(\frac{\sin \pi y}{\pi y}\bigg)^2\bigg) ~dy
 \right).\nonumber
\end{eqnarray}
The expression $1 -(\sin^2 \pi y)/(\pi y)^2$ is exactly the
limiting two-point correlation function predicted by Montgomery
[36].

\section{Mollifying second moments} \label{sect:moll}
The technique of mollifying is used for computing information
about zeros in families of $L$-functions, for example for
obtaining lower bounds for the proportion of zeros on the critical
line or for showing that not many $L$-functions in a family vanish
at the central point. The general set up is that we have a family
of $L$-functions to average over. Before performing the average we
multiply by a Dirichlet polynomial whose coefficients arise from
the inverses of the members of the family, multiplied by a
smoothing function. We will compute one example arising from each
of the three basic symmetry types.

As we discussed in the introduction, mollifier calculations are in
general quite complicated. The ratios conjectures give a
relatively easy way to obtain the relevant asymptotic formula.
Thus, they can serve as a guide as to whether to embark on a
calculation and a check as to whether a calculation is correct.
They also provide evidence that mean-value formulas which can be
proven for short mollifiers remain correct for long mollifiers.
So, these calculations are valuable even though we assume RH.

\subsection{A Unitary example} We start with the Riemann zeta-function in $t$-aspect as a prototype of a unitary family.
So, let
\begin{eqnarray}
M(s,P)=\sum_{n\le y} \frac{\mu(n) P\left(\frac{\log y/n}{\log
y}\right)}{n^s}\end{eqnarray} where $\mu(n)$ is the M\"{o}bius
function,
\begin{equation} \frac{1}{\zeta(s)}=\sum_{n=1}^\infty
\frac{\mu(n)}{n^s},
\end{equation}
and $P$ is a polynomial satisfying $P(0)=0$.  Also,
\begin{equation} y=T^\theta
\end{equation}
where, classically, the following results have been proven for
$\theta< 1/2$, and, with a more modern treatment, for $\theta<
4/7$ [7].   Conjecturally, the asymptotic formula we obtain should
be valid for any fixed $\theta$, no matter how large. We want to
consider
\begin{eqnarray}
I=\int_0^T |\zeta(1/2+it)|^2 |M(1/2+i t,P)|^2 ~dt,
\end{eqnarray}
and more generally
\begin{eqnarray}
I_\zeta(\alpha,\beta,P_1,P_2)=\int_0^T \zeta(s+\alpha)
\zeta(1-s+\beta) M(s,P_1)M(1-s,P_2) ~dt,
\end{eqnarray}
where $s=1/2+it$.
Also, it is useful to discuss the scaled and differentiated form of this quantity, namely,
\begin{eqnarray}
I_\zeta(Q_1,Q_2,P_1,P_2):=  Q_1\bigg(\frac {-1}{\log
T}\frac{d}{d\alpha}\bigg)Q_2 \bigg(\frac {-1}{\log
T}\frac{d}{d\beta}\bigg)I_\zeta(\alpha,\beta,P_1,P_2)\bigg|_{\alpha=\beta=0},
\end{eqnarray}
for polynomials $Q_1$ and $Q_2$.

 To relate this to our ratios
conjecture we note that by Perron's formula
\begin{eqnarray}
\frac{1}{2 \pi i}\int_{(c)}  x^z \frac{dz}{z^{m+1}}  = \bigg\{ {\frac{\log^m x}{m!} \mbox{ if $x>1$} \atop 0 \mbox{ if $0<x< 1$}}
\end{eqnarray}
where $c>0$. Therefore, if $P(x)=\sum_{m\geq 1} p_m x^m$, then
\begin{eqnarray}
M(s,P) =\sum_{m\geq 1} \frac{p_m m!}{\log^my} \frac{1}{2\pi i}
\int_{(c)} \frac{y^z}{z^{m+1}} \frac{1}{\zeta(s+z)} ~ dz.
\end{eqnarray}
This expression leads us to
\begin{eqnarray}
I_\zeta(\alpha,\beta,P_1,P_2)=\sum_{m,n} \frac{p_{1,m} m! p_{2,n}
n! }{\log^{m+n}y} \frac{1}{(2\pi i)^2} \int_{(c_1)} \int_{(c_2)}
\frac{y^{w+z}}{w^{m+1}z^{n+1} } R_\zeta(\alpha,\beta,w,z) ~dw ~dz,
\end{eqnarray}
 where $c_1=c_2=1/\log y$, and $R_{\zeta}$ is defined at (\ref{eq:Rzeta}). Using
 the
 ratios conjecture \ref{conj:unit}, we see that the double integral above is equal
 to
\begin{eqnarray}
&&\frac{1}{(2\pi i)^2} \int_{(c_1)} \int_{(c_2)} \frac{y^{w+z}}{w^{m+1}z^{n+1} }
\int_0^T \left(\frac{\zeta(1+\alpha+\beta)\zeta(1+w+z)}
{\zeta(1+\alpha+z)\zeta(1+\beta+w)}A_\zeta(\alpha, \beta, w, z) \right.
\\
&& \qquad   \left. +\left(\frac {t}{2 \pi}\right)^{-\alpha-\beta}
\frac{\zeta(1-\alpha-\beta)\zeta(1+w+z)}
{\zeta(1-\beta+z)\zeta(1-\alpha+w)}A_\zeta(-\beta, -\alpha,  w,z)
\right)~dt ~dw ~dz +O(T^{1/2+\epsilon})\nonumber
\end{eqnarray}

From this formula we could work out a precise conjecture with all
lower order terms included. However, we are mainly interested in
the leading order term when $\alpha, \beta \approx 1/\log T$.  The
leading order terms come from the residues of the poles in $w$ and
$z$ at zero; to obtain these we use arguments similar to the proof
of the Prime Number Theorem to move the paths of integration
slightly to the left of zero, allowing us to replace the contours
of $z$ and $w$ with circles of radius $1/\log T$ and $2/\log T$
respectively.  The error term is then certainly $1/\log T$ smaller
than the main term.  Also we use $A = 1 +O(1/\log T)$ and
$\zeta(1+x)= 1/x+O(1)$ for small $x$ and large $T$. Then we have
\begin{eqnarray}
I_\zeta(\alpha,\beta,P_1,P_2)&=&\sum_{m,n} \frac{p_{1,m} m!
p_{2,n} n! }
{\log^{m+n}y}\frac{1}{(2\pi i)^2} \oint \oint \frac{y^{w+z}}{w^{m+1}z^{n+1} }\nonumber \\
&&\quad\times \int_0^T \left(\frac{ (\alpha+z) (\beta+w)}{
(\alpha+\beta) (w+z)}
  +\left(\frac {t}{2 \pi}\right)^{-\alpha-\beta}
\frac{ (-\beta+z) (-\alpha+w)}{ (-\alpha-\beta) (w+z)}
  \right)~dt ~dw ~dz \nonumber \\
  && \qquad\qquad+ O(T/\log T).
\end{eqnarray}
It is convenient to write, for $\Re (w+z)>0$,
\begin{equation}
\frac{y^{w+z}}{w+z}=\int_0^y u^{w+z}\frac{du}{u}
\end{equation}
so that the above becomes
\begin{eqnarray}
I_\zeta(\alpha,\beta,P_1,P_2)&=&\frac{1}{\alpha+\beta}\sum_{m,n}
\frac{p_{1,m} m! p_{2,n} n! }{\log^{m+n}y}\int_0^T
\int_1^y\frac{1}{(2\pi i)^2} \oint
\oint \frac{u^{w+z}}{w^{m+1}z^{n+1} } \nonumber\\
&&\quad
  \times\left(  (\alpha+z) (\beta+w)
  -\left(\frac {t}{2 \pi}\right)^{-\alpha-\beta}
  (-\beta+z) (-\alpha+w)
  \right)  ~dw ~dz~\frac{du}{u} ~ dt\nonumber \\
  &&  \qquad\qquad+O(T/\log T);
\end{eqnarray}
 note that the integration in $u$  is for $u\ge 1$ since for $u<1$ the integrals in $z$ and $w$ are 0.

Now
\begin{eqnarray}\label{eqn:poly1}
\sum_m \frac{p_{1,m}m!}{\log^my} \frac{1}{2\pi
i}\oint\frac{u^w}{w^{m+1}} ~dw=P_1\left(\frac{\log u}{\log
y}\right)
\end{eqnarray}
and
\begin{eqnarray} \label{eqn:poly2}
\sum_m \frac{p_{1,m}m!}{\log^my} \frac{1}{2\pi
i}\oint\frac{u^w}{w^{m}} ~dw=\frac{1}{\log y}P_1'\left(\frac{\log
u}{\log y}\right).
\end{eqnarray}
Therefore,
\begin{eqnarray}
I_\zeta(\alpha,\beta,P_1,P_2)&=&\frac{T}{\alpha+\beta} \int_1^y
  \left(  \bigg(\alpha+\frac{d}{dz}\bigg) \bigg(\beta+\frac{d}{dw}\bigg)
  -T^{-\alpha-\beta}
  \bigg(-\beta+\frac{d}{dz}\bigg) \bigg(-\alpha+\frac{d}{dw}\bigg)
  \right) \nonumber\\
  &&\qquad \times P_1\bigg(\frac{w+\log u}{\log y} \bigg)P_2\bigg(\frac{z+\log u}{\log y}
  \bigg)\bigg|_{w=z=0 }~\frac{du}{u} +O(T/\log T).
\end{eqnarray}
Letting $u=y^r$,  we deduce that
\begin{eqnarray}
I_\zeta(\alpha,\beta,P_1,P_2)&=&\frac{T\log y}{\alpha+\beta}
\left( \bigg(\alpha+\frac{d}{dz}\bigg)
\bigg(\beta+\frac{d}{dw}\bigg)
  -T^{-\alpha-\beta}
  \bigg(-\beta+\frac{d}{dz}\bigg) \bigg(-\alpha+\frac{d}{dw}\bigg)
  \right) \nonumber\\
  &&\qquad\times\int_0^1  P_1\bigg(\frac{w}{\log y}+r\bigg)P_2\bigg(\frac{z}{\log y}+r\bigg)
  ~dr\bigg|_{w=z=0}+O(T/\log T).
\end{eqnarray}
It is useful to rewrite the main term of this as
\begin{eqnarray}
&& \frac{T\log y(1-T^{-\alpha-\beta})}{\alpha+\beta}
   \bigg(-\beta+\frac{d}{dz}\bigg) \bigg(-\alpha+\frac{d}{dw}\bigg)
   \int_0^1 P_1\bigg(\frac{w}{\log y}+r\bigg)P_2\bigg(\frac{z}{\log y}+r\bigg)~dr\bigg|_{w=z=0}\nonumber\\
&&\qquad + \frac{T\log
y}{\alpha+\beta}(\alpha+\beta)\bigg(\frac{d}{dw}+\frac{d}{dz}\bigg)
\int_0^1  P_1\bigg(\frac{w}{\log y}+r\bigg)P_2\bigg(\frac{z}{\log
y}+r\bigg)~dr\bigg|_{w=z=0}.
\end{eqnarray}
The second term here is $ =TP_1(1)P_2(1)$. For the first term, we
write
\begin{equation}
\frac{1-T^{-\alpha-\beta}}{\alpha+\beta}=\log T\int_0^1 T^{-u(\alpha+\beta)}~du
\end{equation}
and note, for example,  that $\log y(-\alpha+d/dw)P_1(w/\log
y+r)|_{w=0}=(d/dw)y^{-\alpha w}P_1(w +r) |_{w=0}$.  Finally,
recalling that $y=T^{\theta}$, we have
\begin{eqnarray} \label{eqn:mollunit}
&& I_{\zeta}(\alpha,\beta,P_1,P_2)=
TP_1(1)P_2(1)\\
&&\quad+\frac{T}{\theta} \frac{d}{dw}\frac{d}{dz} y^{-\alpha
w-\beta z }\int_0^1\int_0^1 T^{-(\alpha+\beta)u}
P_1(w+r)P_2(z+r)~dr~du\bigg|_{w=z=0}+O(T/\log T).\nonumber
\end{eqnarray}
This formula appears in [7], page 11. To compute
$I_\zeta(Q_1,Q_2,P_1,P_2)$ we observe, for example, that
\begin{eqnarray}
Q_1\bigg(\frac{-1}{\log T}\frac{d}{d\alpha}\bigg) y^{-\alpha w}T^{-\alpha u}\bigg|_{\alpha=0}=Q_1(w\theta+u).
\end{eqnarray}
Thus, we have
\begin{theorem} Let $P_1,P_2,Q_1$ and $Q_2$ be polynomials, with $P_1(0)=P_2(0)=0$.
Assuming the ratios conjecture \ref{conj:unit}, for any fixed
$\theta>0$, we have (using $s=1/2+it$)
 \begin{eqnarray}
&&  \quad\frac{1}{T} \int_0^T Q_1\bigg(\frac{-1}{\log T} \frac d
{d\alpha}\bigg) Q_2\bigg(\frac{-1}{\log T} \frac d
{d\beta}\bigg)\zeta(s+\alpha)\zeta(1-s+\beta)
M(s,P_1)M(1-s,P_2)~dt\bigg|_{\alpha=\beta=0}
\\&&\quad\qquad=
P_1(1)P_2(1)Q_1(0)Q_2(0) \nonumber \\
&&\qquad\qquad\quad+ \frac{d}{dw}\frac{d}{dz}
\frac{1}{\theta}\int_0^1\int_0^1 P_1(w+r)P_2(z+r)Q_1(w\theta+u)
Q_2(z\theta+u)~dr~du\bigg|_{w=z=0 }\nonumber\\
&&\qquad\qquad \qquad\qquad\qquad+O(1/\log T)\nonumber \\
&&\quad\qquad=
P_1(1)P_2(1)Q_1(0)Q_2(0)\nonumber\\
&&\qquad\qquad\quad+\frac{1}{\theta}\int_0^1\int_0^1\big(P_1'(r)Q_1(u)+\theta
P_1(r)Q_1'(u)\big)\big(P_2'(r)Q_2(u)+\theta
P_2(r)Q_2'(u)\big)~dr~du\nonumber \\
&&\qquad\qquad \qquad\qquad\qquad+O(1/\log T).\nonumber
\end{eqnarray}
\end{theorem}
As remarked earlier, if $\theta< 4/7$, then this is a theorem of
[7] which generalizes work of Levinson [33].  Farmer [16] was the
first to propose that this formula should hold for any fixed value
of $\theta>0$; he calls this the ``long mollifiers" conjecture.
Other examples of mollifying a second moment in a unitary family
are in [17], [25], and [35].

\subsection{A symplectic example}
We consider mollifying in the family of $L$-functions $L(s,\chi_d)$ associated with real Dirichlet characters.
Let
\begin{eqnarray}
M(\chi_d,P)=\sum_{n\le y} \frac{\mu(n)\chi_d(n) P\left(\frac{\log
\frac y n}{\log y}\right)}{n^{1/2}},
\end{eqnarray}
where $P$ is a polynomial satisfying $P(0)=P'(0)=0$  and
$y=X^\theta$.   Consider the second mollified moment
\begin{eqnarray}
{\mathcal M}(\alpha,\beta,P_1,P_2)=\sum_{d\le
X}L(1/2+\alpha,\chi_d)L(1/2+\beta,\chi_d)
M(\chi_d,P_1)M(\chi_d,P_2).
\end{eqnarray}
As in our previous example, we can express
\begin{eqnarray}
M(\chi_d,P)=\sum_n \frac{p_n n!}{\log^ny}\frac{1}{2\pi i}
\int_{(c)} \frac{y^w}{L(1/2+w,\chi_d)w^{n+1}} ~dw,
\end{eqnarray}
where the $p_n$ are the coefficients of the polynomial $P$.  So, letting $p_{m,1}$ and $p_{n,2}$ be the coefficients of $P_1$ and $P_2$ we have
\begin{eqnarray}\label{eqn:M}
&&{\mathcal M}(\alpha,\beta,P_1,P_2)\nonumber\\
&&\qquad  =\sum_{m,n}\frac{p_{m,1}m! p_{n,2}
n!}{\log^{m+n}y}\frac{1}{(2\pi
i)^2}\int_{(c_1)}\int_{(c_2)}\frac{y^{w+z}}{w^{m+1}z^{n+1}}\sum_{d\le
X}\frac{L(1/2+\alpha,\chi_d)L(1/2+\beta,\chi_d)}{L(1/2+w,\chi_d)L(1/2+z,\chi_d)}~dw~dz.
\end{eqnarray}
For the sum over $d$ we substitute from (\ref{eqn:R2symp}); we find that
\begin{eqnarray}
&&\frac{1}{X^*}{\mathcal M}(\alpha,\beta,P_1,P_2)\sim
\sum_{m,n}\frac{p_{m,1}m! p_{n,2} n!}{\log^{m+n}y}\frac{1}{(2\pi i)^2}\int_{(c_1)}\int_{(c_2)}\frac{y^{w+z}}{w^{m+1}z^{n+1}}\\
&&\qquad \bigg(
\frac{(\alpha+w)(\alpha+z)
(\beta+w)(\beta+z)}{4\alpha\beta(\alpha+\beta)(w+z)}-X^{-\alpha}
\frac{(-\alpha+w)(-\alpha+z)
(\beta+w)(\beta+z)}{4\alpha\beta(-\alpha+\beta)(w+z)}\nonumber\\
&&\qquad \qquad -X^{-\beta}\frac{(\alpha+w)(\alpha+z)
(-\beta+w)(-\beta+z)}{4\alpha\beta(\alpha-\beta)(w+z)}\nonumber\\
&&\qquad \qquad \qquad
-X^{-\alpha-\beta}\frac{(-\alpha+w)(-\alpha+z)
(-\beta+w)(-\beta+z)}{4\alpha\beta(\alpha+\beta)(w+z)} \bigg)
~dw~dz. \nonumber\end{eqnarray} For simplicity from now on we
write asymptotic formulas but, as in the previous section, they
could all be replaced by equality with an error term that is one
log smaller than the main term.

As before, we replace $\frac{y^{w+z}}{w+z}$ by $\int_1^y
u^{w+z}\frac{du}{u}$.  Then the poles are all at $w=0$ and $z=0$
and only the numerators in the last set of brackets depend on $w$
and $z$.  Removing the factor $(w+z)$ from the denominator, we
expand this bracket into an expression that is a polynomial of
total degree 4 in $w$ and $z$ with maximum degree 2 in each
variable:
\begin{eqnarray}
&&\frac{1}{4\alpha\beta}\left(\frac{1-X^{-\alpha-\beta}}{\alpha+\beta}+X^{-\alpha}
\frac{1-X^{\alpha-\beta}}{\alpha-\beta}\right)w^2 z^2 +\frac{(1-X^{-\alpha})
(1-X^{-\beta})}{4\alpha\beta}(w^2 z+w z^2)\nonumber\\
&&\quad +
\left(\frac{1-X^{-\alpha-\beta}}{4(\alpha+\beta)}-X^{-\alpha}
\frac{1-X^{\alpha-\beta}}{4(\alpha-\beta)}\right)(w^2+z^2)\nonumber\\
&&\qquad
+\left(\frac{(\alpha+\beta)(1-X^{-\alpha-\beta})}{4\alpha\beta}+
\frac{(\alpha-\beta)(X^{-\alpha}-X^{-\beta})}{4\alpha\beta}\right) w
z +\frac 1 4 (1+X^{-\alpha})(1+X^{-\beta})(w+z)\nonumber\\
&&\qquad \qquad +
\frac{\alpha\beta(1-X^{-\alpha-\beta})}{4(\alpha+\beta)}
+X^{-\alpha}\frac{\alpha\beta(1-X^{\alpha-\beta})}{4(\alpha-\beta)}.
\end{eqnarray}

Using the analogue of formulas (\ref{eqn:poly1}) and (\ref{eqn:poly2}), we see that
we now should replace $w^2 z^2$ in this expression by
\begin{eqnarray}\label{eqn:analogue}
\frac{1}{\log^4 y}\int_1^y P_1''\left(\frac{\log u}{\log y} \right) P_2''\left(\frac{\log u}{\log y} \right) ~\frac{du}{u}
=\frac{1}{\log^3 y}\int_0^1 P_1''(r)P_2''(r)~dr.
\end{eqnarray}
Likewise, $w^2z+w z^2$ should be replaced by
\begin{eqnarray}
\frac{1}{\log^2 y} \int_0^1\big(P_1''(r)P_2'(r)+P_1'(r)P_2''(r)\big)~dr ,
\end{eqnarray}
$w^2+z^2$ by
\begin{eqnarray}
\frac{1}{\log y}\int_0^1\big(P_1''(r)P_2(r)+P_1(r)P_2''(r)\big)~dr ,
\end{eqnarray}
$w z$ by
\begin{eqnarray}
\frac{ \int_0^1 P_1'(r)P_2'(r)~dr }{\log y},
\end{eqnarray}
$w+z$ by
\begin{eqnarray}
\int_0^1 \big(P_1'(r)P_2(r)+P_1(r)P_2'(r)\big)~dr ,
\end{eqnarray}
and the constant term by
\begin{eqnarray}
\log y \int_0^1 P_1(r)P_2(r)~dr.
\end{eqnarray}
In this way, we find that
 \begin{eqnarray} \label{eqn:ff}
&&\frac{4}{X^*}{\mathcal M}(\alpha,\beta,P_1,P_2)\sim
\frac{1}{\alpha\beta}\left(\frac{1-X^{-\alpha-\beta}}{\alpha+\beta}+X^{-\alpha}
\frac{1-X^{\alpha-\beta}}{\alpha-\beta}\right)\frac{\int_0^1 P_1''(r)P_2''(r)~dr}{\log^3y} \\
&&\quad +\frac{(1-X^{-\alpha})(1-X^{-\beta})}{\alpha\beta}\frac{
\int_0^1\big(P_1''(r)P_2'(r)+P_1'(r)P_2''(r)\big)~dr}{\log^2y}
\nonumber \\ &&\qquad+
\left(\frac{1-X^{-\alpha-\beta}}{(\alpha+\beta)}-X^{-\alpha}
\nonumber \frac{1-X^{\alpha-\beta}}{(\alpha-\beta)}\right)
\frac{\int_0^1\big(P_1''(r)P_2(r)+P_1(r)P_2''(r)\big)~dr}{\log y}\\
&&\qquad
+\left( (1+X^{-\alpha})\frac{1-X^{-\beta}}{\beta}+ \nonumber (1+X^{-\beta})\frac{1-X^{-\alpha}}{\alpha } \right) \frac{ \int_0^1 P_1'(r)P_2'(r)~dr }{\log y}\\
&&\nonumber\qquad \qquad + (1+X^{-\alpha})(1+X^{-\beta}) \int_0^1 \big(P_1'(r)P_2(r)+P_1(r)P_2'(r)\big)~dr\\
&&\qquad \qquad +
\left(\frac{\alpha\beta(1-X^{-\alpha-\beta})}{\alpha+\beta}\nonumber
+X^{-\alpha}\frac{\alpha\beta(1-X^{\alpha-\beta})}{\alpha-\beta}\right)\log
y \int_0^1 P_1(r)P_2(r)~dr.
\end{eqnarray}
This gives our final formula for ${\mathcal M}(\alpha,\beta,P_1,P_2)$.

If, instead, we consider the mollified second moment of $\xi(1/2,\chi_d)$, we can put our answer into a more symmetric form. Recall that by the functional equation (\ref{eqn:functional}) we have
\begin{eqnarray}
\xi(1/2+\alpha,\chi_d):=\left(\frac{d}{\pi}\right)^{\alpha/2}\Gamma\left(\frac 14 +\frac \alpha 2\right)L(1/2+\alpha,\chi_d)=\xi(1/2-\alpha,\chi_d).
\end{eqnarray}
Therefore,  if we multiply (\ref{eqn:ff}) by $X^{(\alpha+\beta)/2}$ we will obtain the asymptotic formula for the mollified second moment of $\xi$:
\begin{eqnarray}
\mathcal{N}(\alpha,\beta,P_1,P_2):=\sum_{d\le X} \xi(1/2+\alpha,\chi_d)\xi(1/2+\beta,\chi_d)M(\chi_d,P_1)M(\chi_d,P_2).
\end{eqnarray}
We have
\begin{eqnarray} \label{eqn:fxif}
&&\frac{4}{X^*}{\mathcal N}(\alpha,\beta,P_1,P_2)\sim \nonumber
\frac{1}{\alpha\beta}\left(\frac{X^{\frac{\alpha+\beta}{2}}-X^{\frac{-(\alpha+\beta)}2}}
{\alpha+\beta}-\frac{X^{\frac{\alpha-\beta}{2}}-X^{\frac{\beta-\alpha}{2}}}
 {\alpha-\beta}\right)\frac{\int_0^1 P_1''(r)P_2''(r)~dr}{\log^3y} \\
&&\quad +\frac{(X^{\frac \alpha 2}-X^{\frac{-\alpha}2})(X^{\frac \beta 2}-X^{\frac{-\beta}{2}})}{\alpha\beta}\frac{\int_0^1\big(P_1''(r)P_2'(r)+P_1'(r)P_2''(r)\big)~dr}{\log^2y} \nonumber \\ &&\qquad+ \left(\frac{X^{\frac{\alpha+\beta}{2}}-X^{\frac{-(\alpha+\beta)}2}}{(\alpha+\beta)}+ \nonumber
\frac{X^{\frac{\alpha-\beta}{2}}-X^{\frac{\beta-\alpha}{2}}}
{(\alpha-\beta)}\right)\frac{\int_0^1\big(P_1''(r)P_2(r)+P_1(r)P_2''(r)\big)~dr}{\log y}\\
&&\qquad
+\left( (X^{\frac \alpha 2}+X^{\frac{-\alpha}2})\frac{X^{\frac \beta 2}-X^{\frac{-\beta}{2}}}{\beta}+   (X^{\frac \beta 2}+X^{\frac{-\beta}{2}})\frac{X^{\frac \alpha 2}-X^{\frac{-\alpha}2}}{\alpha } \right) \frac{ \int_0^1 P_1'(r)P_2'(r)~dr }{\log y}\\
&&\qquad \qquad  + (X^{\frac \alpha 2}+X^{\frac{-\alpha}2})(X^{\frac \beta 2}+X^{\frac{-\beta}{2}}) \int_0^1 \big(P_1'(r)P_2(r)+P_1(r)P_2'(r)\big)~dr\nonumber \\
&&\qquad \qquad +
\left(\frac{\alpha\beta(X^{\frac{\alpha+\beta}{2}}-X^{\frac{-(\alpha+\beta)}2})}
{\alpha+\beta}\nonumber - \frac{\alpha\beta(X^{\frac{\alpha-\beta
}{2}} -X^{\frac{\beta-\alpha }{2}})}{\alpha-\beta}\right)\log y
\int_0^1 P_1(r)P_2(r)~dr.
\end{eqnarray}
We introduce a scaling, writing $\alpha=2a /\log X$ and $\beta=2b/
\log X$. Then it is not difficult, remembering that
$y=X^{\theta}$, to see that the above can be rewritten as
\begin{eqnarray} \label{eqn:fxif1}
&&\frac{4}{X^*}{\mathcal N}(\alpha,\beta,P_1,P_2)\sim \nonumber
\frac 1 {2\theta^3}\int_0^1\frac{\sinh   au }{a }
\frac{\sinh   bu }{b }~du
 \int_0^1 P_1''(r)P_2''(r)~dr  \\
&&\quad + \frac 1 {\theta^2} \frac{\sinh  a  }{a} \frac
{\sinh b}{b}\int_0^1\big(P_1''(r)P_2'(r)+P_1'(r)P_2''(r)\big)~dr \nonumber \\ &&\qquad+ \nonumber
\frac{2}{\theta}\int_0^1 \cosh au \cosh bu~du
\int_0^1\big(P_1''(r)P_2(r)+P_1(r)P_2''(r)\big)~dr \\
&&\qquad
+\frac{2}{\theta}\left(   \frac{\cosh a\sinh b}{b}+
 \frac{\cosh b\sinh a}{a} \right)   \int_0^1 P_1'(r)P_2'(r)~dr \\
&&\qquad \qquad  + 4\cosh a\cosh b  \int_0^1 \big(P_1'(r)P_2(r)+P_1(r)P_2'(r)\big)~dr\nonumber \\
&&\qquad \qquad + 8\theta \nonumber a b \int_0^1 \sinh au \sinh
bu~du\int_0^1 P_1(r)P_2(r)~dr.
\end{eqnarray}

We now apply
$Q_1\left(\frac{d}{da}\right)Q_2\left(\frac{d}{db}\right)$ to this
expression to obtain
\begin{eqnarray}
\mathcal{N}(Q_1,Q_2,P_1,P_2):=
Q_1\left(\frac{d}{da}\right)Q_2\left(\frac{d}{db}\right)
\mathcal{N}\left(\frac {2a}{\log X},\frac {2b}{\log
X},P_1,P_2\right)\bigg|_{a=b=0}.
\end{eqnarray}
We may assume that $Q_1$ and $Q_2$ are even functions, since
for an odd number $r$ we have $\xi^{(r)}(1/2,\chi_d)=0$. To perform this calculation, we observe, for example, that
\begin{eqnarray}
&&Q_1\left(\frac{d}{da}\right)Q_2\left(\frac{d}{db}\right)\int_0^1\frac{\sinh
au  }{a }
\frac{\sinh bu}{b }~du\bigg|_{a=b=0}\\
&&\quad =
Q_1\left(\frac{d}{da}\right)Q_2\left(\frac{d}{db}\right)\int_0^1
\int_0^u \cosh  at_1  ~dt_1\int_0^u \cosh  bt_2 ~dt_2~du\bigg|_{a=b=0}\nonumber\\
&&\quad = \frac 1 4 \int_0^1 \int_0^u (Q_1(t_1)+Q_1(-t_1))~dt_1\int_0^u
(Q_2(t_2)+Q_2(-t_2))~dt_2 ~du\nonumber\\
&&\quad =  \int_0^1\tilde{Q_1}(u)\tilde{Q_2}(u)~du,\nonumber
\end{eqnarray}
where we have used the notation
\begin{eqnarray}
\tilde{Q}(u)=\int_0^u Q(t)~dt.
\end{eqnarray}
By similar, but easier, calculations we find that
\begin{eqnarray}
&&\nonumber\frac 4 {X^*}\mathcal{N}(Q_1,Q_2,P_1,P_2)\sim \frac 1
{2\theta^3}\int_0^1\tilde{Q_1}(u)\tilde{Q_2}(u)~du
 \int_0^1 P_1''P_2''(r)~dr  \\
&&\quad + \frac 1 {\theta^2}\tilde{Q_1}(1)\tilde{Q_2}(1)\int_0^1\big(P_1''(r)P_2'(r)+P_1'(r)P_2''(r)\big)~dr \nonumber \\ &&\qquad+ \nonumber
\frac{2}{\theta}\int_0^1 Q_1(u)Q_2(u)~du
\int_0^1\big(P_1''(r)P_2(r)+P_1(r)P_2''(r)\big)~dr \\
&&\qquad
+\frac{2}{\theta}\left(  Q_1(1)\tilde{Q_2}(1)+\tilde{Q_1}(1)Q_2(1) \right)   \int_0^1 P_1'(r)P_2'(r)~dr \\
&&\qquad \qquad  + 4Q_1(1)Q_2(1)  \int_0^1
\big(P_1'(r)P_2(r)+P_1(r)P_2'(r) \big)~dr\nonumber \\&&\qquad
\qquad + 8\theta \nonumber \int_0^1 Q_1'(u)Q_2'(u)~du\int_0^1
P_1(r)P_2(r)~dr.
\end{eqnarray}
The right hand side here can be written in a more compact form as
\begin{eqnarray}
&&\frac{1}{2\theta}\int_0^1\int_0^1 \left(\frac 1 \theta P_1''(r)\tilde{Q_1}(u)-4\theta P_1(r)Q_1'(u) \right)\left( \frac{1}{\theta} P_2''(r)\tilde{Q_2}(u)-4\theta P_2(r) Q_2'(u) \right)~du ~dr \\
&&\qquad + \left(\frac 1\theta P_1'(1)\tilde{Q_1}(1)+2
P_1(1)Q_1(1)\right) \left(\frac 1\theta P_2'(1)\tilde{Q_2}(1)+2
P_2(1)Q_2(1)\right).\nonumber
\end{eqnarray}
To verify this assertion we need to use identities which follow from integration-by-parts, such as
\begin{eqnarray}
\int_0^1 \left(P_1'(r)P_2(r)+P_1(r)P_2'(r)\right)~dr=P_1(1)P_2(1),
\end{eqnarray}
\begin{eqnarray}
\int_0^1 \left(P_1''(r)P_2'(r)+P_1'(r)P_2''(r)\right)~dr=P_1'(1)P_2'(1),
\end{eqnarray}
\begin{eqnarray}
\int_0^1 P_1''(r)P_2(r)~dr=P_1'(1)P_2(1)-\int_0^1 P_1'(r)P_2'(r)~dr,
\end{eqnarray}
and
\begin{eqnarray}
\int_0^1\tilde{Q_1}(u)Q_2'(u)~du = \tilde{Q_1}(1)Q_2(1)-\int_0^1
Q_1(u)Q_2(u)~du.
\end{eqnarray}
In the last equation note that we have used $\tilde{Q}(0)=0$.
\begin{theorem} \label{theo:mollsymplec}
Assuming Conjecture \ref{conj:symp2}, we have for even polynomials
$Q_1$ and $Q_2$, and $P_1$ and $P_2$ polynomials satisfying
$P_1(0)=P_1'(0)=P_2(0)=P_2'(0)=0$, and $y=X^\theta$ with any
$\theta>0$,
\begin{eqnarray}
&&Q_1\left(\frac{2}{\log
X}\frac{d}{d\alpha}\right)Q_2\left(\frac{2} {\log
X}\frac{d}{d\beta}\right)\sum_{d\le X}
\xi(1/2+\alpha,\chi_d)\xi(1/2+\beta,\chi_d)
M(\chi_d,P_1)M(\chi_d,P_2)\bigg|_{\alpha=\beta=0}\nonumber\\
&&\quad = X^*\bigg(\frac{1}{8\theta}\int_0^1\int_0^1 \left(\frac 1
\theta P_1''(r) \tilde{Q_1}(u)-4\theta P_1(r)Q_1'(u) \right)\left(
\frac{1}{\theta} P_2''(r)
\tilde{Q_2}(u)-4\theta P_2(r) Q_2'(u) \right)~du ~dr\nonumber \\
&&\qquad \qquad + \frac{1}{4}\left(\frac 1\theta
P_1'(1)\tilde{Q_1}(1)+2 P_1(1)Q_1(1)\right) \left(\frac 1\theta
P_2'(1)\tilde{Q_2}(1)+2 P_2(1)Q_2(1)\right)+O(1/\log X)\bigg).
\end{eqnarray}
\end{theorem}

Examples of second moment mollifying in a symplectic family occur
in [40] and [15].

\subsection{An orthogonal example}
Here we compute
\begin{eqnarray}
\quad\mathcal M_\Delta(\alpha,\beta;P_1,P_2):=\sum_{d\le X}
L_\Delta(1/2+\alpha,\chi_d)L_\Delta(1/2+\beta,\chi_d)M_\Delta(\chi_d,P_1)M_\Delta(\chi_d,P_2),
\end{eqnarray}
where
\begin{eqnarray}
M_\Delta(\chi_d,P):=\sum_{m\le y} \frac{\mu_\Delta(m)\chi_d(m)
P\left(\frac {\log \frac y m}{\log y}\right)}{m^{1/2}}.
\end{eqnarray}
As in equation (\ref{eqn:M}), we have
\begin{eqnarray}\label{eqn:MD}
&&{\mathcal M}_\Delta(\alpha,\beta,P_1,P_2)\nonumber\\
&&\qquad  =\sum_{m,n}\frac{p_{m,1}m! p_{n,2}
n!}{\log^{m+n}y}\frac{1}{(2\pi
i)^2}\int_{(c_1)}\int_{(c_2)}\frac{y^{w+z}}{w^{m+1}z^{n+1}}\sum_{d\le
X}\frac{L_\Delta(1/2+\alpha,\chi_d)L_\Delta(1/2+\beta,\chi_d)}
{L_\Delta(1/2+w,\chi_d)L_\Delta(1/2+z,\chi_d)}~dw~dz.
\end{eqnarray}
Using (\ref{eqn:ratorth2}) leads to
\begin{eqnarray}
\mathcal{M}_\Delta(\alpha,\beta,P_1,P_2)&\sim&
\frac{1}{X^*}\sum_{m,n}\frac{p_{m,1}m! p_{n,2} n!}{\log^{m+n}y}\frac{1}{(2\pi
i)^2}\int_{(c_1)}\int_{(c_2)}\frac{y^{w+z}}{w^{m+1}z^{n+1}4wz(w+z)}\nonumber\\
&&\bigg(
\frac{(\alpha+w)(\alpha+z)(\beta+w)(\beta+z)}{(\alpha+\beta) }
\nonumber\\&& \qquad + X^{-2\alpha}
\frac{(-\alpha+w)(-\alpha+z)(\beta+w)(\beta+z)}{(-\alpha+\beta) }\\
&& \qquad \qquad +X^{-2\beta}\nonumber
\frac{(\alpha+w)(\alpha+z)(-\beta+w)(-\beta+z)}{(\alpha-\beta) }
\\ && \qquad \qquad \qquad
- X^{-2\alpha-2\beta}\nonumber
\frac{(-\alpha+w)(-\alpha+z)(-\beta+w)(-\beta+z)}{(\alpha+\beta) }
\bigg) ~dw ~dz.
\end{eqnarray}
We expand the brackets into powers of $w$ and $z$ yielding
\begin{eqnarray}
&&\left(\frac{1-X^{-2\alpha-2\beta}}{\alpha+\beta}-\frac{X^{-2\alpha}-X^{-2\beta}}{\alpha-\beta}
\right)w^2 z^2\\
&&\qquad
+(1+X^{-2\alpha})(1+X^{-2\beta}) (w^2 z+w z^2)\nonumber\\
&& \qquad +\bigg(\frac{\alpha \beta(1-X^{-2
\alpha-2\beta})}{\alpha+\beta}  +\frac{\alpha \beta
(X^{-2\alpha}-X^{-2\beta})}{\alpha-\beta}\bigg)(w^2+z^2)\nonumber\\
&&\qquad \qquad  +\bigg((\alpha+\beta)(1-X^{-2\alpha-2\beta})-
(\alpha-\beta)(X^{-2\alpha}-X^{-2\beta})\bigg) w z\nonumber\\
&&\qquad \qquad \qquad + \alpha \beta(1-X^{-2 \alpha})(1-X^{-2 \beta}) (w+z)\nonumber\\
&&\qquad \qquad \qquad \qquad + \bigg(\frac{\alpha^2
\beta^2(1-X^{-2 \alpha-2\beta})}{\alpha+\beta}  -\frac{\alpha^2
\beta^2 (X^{-2\alpha}-X^{-2\beta})}{\alpha-\beta}\bigg).\nonumber
\end{eqnarray}
As we did in the other cases, we replace $\frac{y^{w+z}}{w+z}$ by
$\int_1^y u^{w+z}\frac{du}{u}$. In a similar manner to
(\ref{eqn:analogue}), we evaluate the sums over $m$ and $n$ using
\begin{eqnarray}
&&\sum_{m,n}\frac{p_{m,1}m! p_{n,2} n!}{\log^{m+n}y}\frac{1}{(2\pi
i)^2}\int_{(c_1)}\int_{(c_2)}\int_1^y
\frac{u^{w+z}}{w^{m+1}z^{n+1}wz}~\frac{du}{u} ~dw ~dz\\
&&\qquad\qquad\ =\log^2 y\int_1^y \tilde{P_1}\left(\frac{\log
u}{\log y}\right) \tilde{P_2}\left(\frac{\log u}{\log y}\right)
~\frac{du}{u},\nonumber
\end{eqnarray}
\begin{eqnarray}
&&\sum_{m,n}\frac{p_{m,1}m! p_{n,2} n!}{\log^{m+n}y}\frac{1}{(2\pi
i)^2}\int_{(c_1)}\int_{(c_2)}\int_1^y
\frac{u^{w+z}}{w^{m+1}z^{n+1}}~\frac{du}{u} ~dw ~dz \\
&&\qquad\qquad=\int_1^y P_1\left(\frac{\log u}{\log y}\right)
P_2\left(\frac{\log u}{\log y}\right) ~\frac{du}{u},\nonumber
\end{eqnarray}
and, in general,
\begin{eqnarray}\label{eqn:poly}
&&\sum_{m,n}\frac{p_{m,1}m! p_{n,2} n!}{\log^{m+n}y}\frac{1}{(2\pi
i)^2}\int_{(c_1)}\int_{(c_2)}\int_1^y \frac{u^{w+z}w^a z^b }{w^{m+1}z^{n+1}wz}~\frac{du}{u} ~dw ~dz \\
 && \qquad \qquad =
(\log y)^{2-a-b}\int_1^y P_1^{(a-1)}\left(\frac{\log u}{\log \nonumber y}\right) P_2^{(b-1)}\left(\frac{\log u}{\log y}\right) ~\frac{du}{u}\\
&& \qquad \qquad  =(\log y)^{3-a-b} \int_0^1
P_1^{(a-1)}(t)P_2^{(b-1)}(t)~dt, \nonumber
\end{eqnarray}
where $P^{(a)}$ means the $a$th derivative of $P$; if $a<0$ then it means the $(-a)$th integral of $P$, so that, for example, $P^{(-1)}=\tilde{P}$.
Inputting this into our expression for $\mathcal M_\Delta$ leads to
\begin{eqnarray}
&&\frac{4}{X^*}\mathcal{M}_\Delta(\alpha,\beta,P_1,P_2)
\nonumber\\
&& \quad \sim \frac{1}
{\log y}\left(\frac{1-X^{-2\alpha-2\beta}}{\alpha+\beta}-\frac{X^{-2\alpha}-X^{-2\beta}}
{\alpha-\beta}\right)\int_0^1 P_1'(t)P_2'(t)~dt\nonumber\\
&&\quad \quad
+(1+X^{-2\alpha})(1+X^{-2\beta}) \int_0^1\left(P_1(t)P_2'(t)+P_1'(t)P_2(t)\right)~dt\nonumber\\
&& \qquad \quad +\log y \bigg(\frac{\alpha \beta(1-X^{-2
\alpha-2\beta})}{\alpha+\beta}  +\frac{\alpha \beta
(X^{-2\alpha}-X^{-2\beta})}{\alpha-\beta}\bigg)\int_0^1\left(P_1'(t)\tilde{P_2}(t)
+\tilde{P_1}(t)P_2'(t)\right)~dt\\
&&\qquad \qquad  +\log y\bigg((\alpha+\beta)(1-X^{-2\alpha-2\beta})-
(\alpha-\beta)(X^{-2\alpha}-X^{-2\beta})\bigg)\int_0^1P_1(t)P_2(t)~dt\nonumber\\
&&\qquad \qquad \qquad + \log^2 y~\alpha \beta(1-X^{-2 \alpha})(1-X^{-2 \beta})
 \int_0^1\left(\tilde{P_1}(t)P_2(t)+P_1(t)\tilde{P_2}(t)\right)~dt\nonumber\\
&&\qquad \qquad \qquad \qquad + \log^3 y ~\bigg(\frac{\alpha^2
\beta^2(1-X^{-2 \alpha-2\beta})}{\alpha+\beta}  -\frac{\alpha^2
\beta^2
(X^{-2\alpha}-X^{-2\beta})}{\alpha-\beta}\bigg)\int_0^1\tilde{P_1}(t)\tilde{P_2}(t)~dt.\nonumber
\end{eqnarray}
We want to compare mollifying in an orthogonal family with that in a symplectic family. To this end, we consider, as we did for the symplectic family, mollifying the xi-functions. In this situation it just means multiplying the above result by
$X^{\alpha+\beta}$. This gives
\begin{eqnarray}&&
\frac{4}{X^*}\sum_{d\le X} \xi_\Delta(1/2+\alpha,\chi_d)\xi_\Delta(1/2+\beta,\chi_d)
M_\Delta(\chi_d,P_1)M_\Delta(\chi_d,P_2)\nonumber\\
&&\quad \sim \frac{1} {\log
y}\left(\frac{X^{\alpha+\beta}-X^{-\alpha-\beta}}{\alpha+\beta}
+\frac{X^{\alpha-\beta}-X^{\beta-\alpha}}{\alpha-\beta}\right)\int_0^1 P_1'(t)P_2'(t)~dt\nonumber\\
&&\quad \quad
+(X^\alpha+X^{-\alpha})(X^\beta+X^{-\beta}) \int_0^1\left(P_1(t)P_2'(t)+P_1'(t)P_2(t)\right)~dt\nonumber\\
&& \qquad \quad +\log y \bigg(\frac{\alpha
\beta(X^{\alpha+\beta}-X^{- \alpha-\beta})}{\alpha+\beta}
-\frac{\alpha \beta
(X^{\alpha-\beta}-X^{\beta-\alpha})}{\alpha-\beta}\bigg)\int_0^1\left(P_1'(t)\tilde{P_2}(t)+\tilde{P_1}(t)P_2'(t)\right)~dt\\
&&\qquad \qquad  +\log y\bigg((\alpha+\beta)(X^{\alpha+\beta}-X^{-\alpha-\beta})+
(\alpha-\beta)(X^{\alpha-\beta}-X^{\beta-\alpha})\bigg)\int_0^1P_1(t)P_2(t)~dt\nonumber\\
&&\qquad \qquad \qquad + \log^2 y~\alpha \beta(X^\alpha-X^{-
\alpha})
(X^\beta-X^{- \beta}) \int_0^1\left(\tilde{P_1}(t)P_2(t)+P_1(t)\tilde{P_2}(t)\right)~dt\nonumber\\
&&\qquad \qquad \qquad \qquad + \log^3 y ~\bigg(\frac{\alpha^2
\beta^2(X^{\alpha+\beta}-X^{- \alpha-\beta})}{\alpha+\beta}
+\frac{\alpha^2 \beta^2
(X^{\alpha-\beta}-X^{\beta-\alpha})}{\alpha-\beta}\bigg)\int_0^1\tilde{P_1}(t)\tilde{P_2}(t)~dt.\nonumber
\end{eqnarray}
If we now scale, letting $\alpha=a/\log X$ and $\beta=b/\log X$,
and continuing our $y=X^\theta$ convention, then we can rewrite
the above as
\begin{eqnarray}&&
\frac{4}{X^*}\sum_{d\le X} \xi_\Delta\left(\frac 12 +\frac{a}{\log X},\chi_d\right)
\xi_\Delta\left(\frac 12 +\frac{b}{\log X},\chi_d\right)M_\Delta(\chi_d,P_1)M_\Delta(\chi_d,P_2)\nonumber\\
&&\quad \sim \frac{2}{\theta}
 \left( \frac{\sinh(a+b)}{a+b}+\frac{\sinh(a-b)}{a-b}\right) \int_0^1 P_1'(t)P_2'(t)~dt\nonumber\\
&&\quad \quad
+4\cosh a\cosh b \int_0^1\left(P_1(t)P_2'(t)+P_1'(t)P_2(t)\right)~dt\nonumber\\
&& \qquad \quad +\theta \bigg(\frac{2a b \sinh(a+b)}{a+b}
-\frac{2a b
\sinh (a-b)}{a-b}\bigg)\int_0^1 \left(P_1'(t)\tilde{P_2}(t)+\tilde{P_1}(t)P_2'(t)\right)~dt\\
&&\qquad \qquad  +\theta\bigg(2(a+b)\sinh(a+b)+
2(a-b)\sinh(a-b)\bigg)\int_0^1P_1(t)P_2(t)~dt\nonumber\\
&&\qquad \qquad \qquad + \theta^2~4ab\sinh a\sinh b \int_0^1\left(\tilde{P_1}(t)P_2(t)+P_1(t)
\tilde{P_2}(t)\right)~dt\nonumber\\
&&\qquad \qquad \qquad \qquad + \theta^3
~\bigg(\frac{2a^2b^2\sinh(a+b)}{a+b} +\frac{2a^2b^2
\sinh(a-b)}{a-b}\bigg)\int_0^1\tilde{P_1}(t)\tilde{P_2}(t)~dt.\nonumber
\end{eqnarray}
We now apply $Q_1(\frac d{da})Q_2(\frac d {db})|_{a=b=0}$ to both
sides of this expression; we assume that $Q_1$ and $Q_2$ are even.
We  use the notation  $\mathcal M_\Delta(Q_1,Q_2,P_1,P_2)$ as in
the symplectic example. Thus,
\begin{eqnarray}&&
\frac{4}{X^*}\mathcal{M}_\Delta(Q_1,Q_2,P_1,P_2)\nonumber\\
&&\quad \sim \frac{4}{\theta}
 \int_0^1Q_1(u)Q_2(u)~du \int_0^1 P_1'(t)P_2'(t)~dt\nonumber\\
&&\quad \quad
+4Q_1(1)Q_2(1) \int_0^1\left(P_1(t)P_2'(t)+P_1'(t)P_2(t)\right)~dt\nonumber\\
&& \qquad \quad +4\theta\int_0^1 Q_1'(u)Q_2'(u)~du \int_0^1 \left(P_1'(t)\tilde{P_2}(t)+\tilde{P_1}(t)P_2'(t)\right)~dt\\
&&\qquad \qquad  +4\theta(Q_1(1)Q_2'(1)+Q_1'(1)Q_2(1))\int_0^1P_1(t)P_2(t)~dt\nonumber\\
&&\qquad \qquad \qquad + 4\theta^2~Q_1'(1)Q_2'(1) \int_0^1\left(\tilde{P_1}(t)P_2(t)+
P_1(t)\tilde{P_2}(t)\right)~dt\nonumber\\
&&\qquad \qquad \qquad \qquad + 4\theta^3 ~\int_0^1
Q_1''(u)Q_2''(u)~du
\int_0^1\tilde{P_1}(t)\tilde{P_2}(t)~dt.\nonumber
\end{eqnarray}
This expression can be simplified to obtain
\begin{theorem} \label{theo:mollorthog}
Assuming Conjecture \ref{conj:orth2}, with even polynomials $Q_1$,
$Q_2$, and polynomials $P_1$, $P_2$,  satisfying
$P_1(0)=P_2(0)=0$, and using $y=X^\theta$, we have for arbitrary
$\theta$,
 \begin{eqnarray}
&&\frac{1}{X^*}Q_1\left(\frac{1}{\log X} \frac{d}{d\alpha}\right)
Q_2\left(\frac{1}{\log X} \frac{d}{d\beta}\right)\\
&&\qquad \times\sum_{d\le X}
\xi_\Delta(1/2+\alpha,\chi_d)\xi_\Delta(1/2+\beta,\chi_d)
M_\Delta(\chi_d, P_1)M_\Delta(\chi_d,P_2) \bigg|_{\alpha=\beta=0}\nonumber\\
&& \quad \qquad=
  \frac{1}{\theta}\int_0^1\int_0^1 \bigg(P_1'(t)Q_1(u)-\theta^2\tilde{P_1}(t)Q_1''(u)\bigg)
  \bigg(P_2'(t)Q_2(u)-\theta^2\tilde{P_2}(t)Q_2''(u)\bigg)~dt ~du \nonumber\\
&&\qquad \qquad
+\bigg(P_1(1)Q_1(1)+\theta\tilde{P_1}(1)Q_1'(1)\bigg)
\bigg(P_2(1)Q_2(1)+\theta\tilde {P_2}(1)Q_2'(1)\bigg)\nonumber\\
&&\qquad\qquad + \theta\bigg( Q_1'(0)Q_2(0)\int_0^1
\tilde{P}_1(t)P'_2(t) dt + Q_1(0)Q'_2(0) \int_0^1
P'_1(t)\tilde{P}_2(t)dt \bigg)+O(1/\log X).\nonumber
\end{eqnarray}
\end{theorem}

Examples of second moment mollifying in an orthogonal family occur
in [26] and [31]; in [32] a fourth moment mollification is
performed.


\section{Mollifying the $k$th moment of
$\zeta(s)$}\label{sect:mollk}
 Chris Hughes has unpublished notes
giving an asymptotic  formula for
\begin{eqnarray}\int_0^T|\zeta(1/2+it)|^4 |A(1/2+it)|^2~dt
\end{eqnarray}
where
\begin{eqnarray} A(s)=\sum_{n\le y}\frac{a_n}{n^s}
\end{eqnarray}
is an arbitrary Dirichlet polynomial and where $y=T^\theta$ with
$\theta < 5/27$. For applications to zeros of $\zeta(s)$ it would
extremely useful to specialize this formula to the case that
$A(s)=M(s)$ is a mollifying polynomial, but this would still
involve a lot of work.  Via ratios we produce a conjectural
formula which can serve as a check against the more complicated
rigorous proof via Hughes' formula.  There are (at least) two
obvious choices for a mollifying polynomial $M(s)$. One is $
M(s)=M_1(s,P)^2 $ where
\begin{eqnarray}
\label{eq:M1} M_1(s,P)=\sum_{n\le y}\frac{\mu(n)P\left(\frac {\log
\frac y n}{\log y}\right)}{n^s}
\end{eqnarray}
with $y=T^{\theta}$.  The other is $M(s)=M_2(s,P)$ with
\begin{eqnarray}
\label{eq:M2} M_2(s,P)= \sum_{n\le y}\frac{\mu_2(n)P\left(\frac
{\log \frac y n}{\log y}\right)}{n^s}
\end{eqnarray}
where $y=T^{\theta}$ and $\mu_2$ is the coefficient in the
generating function for $1/\zeta(s)^2$.

Here we will compute what the ratios conjecture tells us about the
asymptotics for the $k$th mollified moments  in the   case where
we mollify with $M_k(s,P)$, where $P(x)=\sum_mp_mx^m$ is a
polynomial satisfying
\begin{eqnarray}
P(0)=P'(0)=\cdots=P^{(k^2-1)}(0)=0.\end{eqnarray} These conditions
on $P(x)$ ensure that we have a smooth cut-off at $n=y$.  It is
only in the course of the calculation that we see why we need
$k^2-1$ derivatives to be zero.

We note that
\begin{eqnarray}
M_k(s,P)= \sum_{n\le y} \frac{\mu_k(n) P \left(\frac{\log
\frac{y}{n}}{\log y}\right) }{ n^s}=\sum_m \frac{p_m m!}{(\log
y)^m}\frac{1}{2\pi i} \int_{(c)} \frac{y^w}{\zeta^k(s+w) w^{m+1}}
~dw,
\end{eqnarray}
where $\mu_k$ is the coefficient in the generating function for
$1/\zeta(s)^k$, $y=T^{\theta}$ and $c>0$.

Thus, using $s=1/2+it$,
\begin{eqnarray}
&&{\mathcal M}_k(\alpha,\beta):=\frac{1}{T}\int_0^T
\zeta(s+\alpha_1)\cdots\zeta(s+\alpha_k)\zeta(1-s-\beta_1)\cdots\\
&&\qquad\qquad\qquad\qquad\times
\zeta(1-s-\beta_k) M_k(s,Q)M_k(1-s,P) ~dt\nonumber\\
&&\qquad = \sum_{m,n} \frac{q_{m}m!
p_{n}n!}{\log^{m+n}y}\frac{1}{(2\pi i)^2}\int_{(c_1)}
\int_{(c_2)}\frac{y^{w+z}}{w^{m+1}z^{n+1}}\nonumber\\
&& \qquad \qquad \times \frac{1}{T}\int_0^T
\frac{\zeta(s+\alpha_1)\cdots\zeta(s+\alpha_k)\zeta(1-s-\beta_1)\cdots
\zeta(1-s-\beta_k) }{\zeta(s+w)^k\zeta(1-s+z)^k}~dt ~dw
~dz\nonumber
\end{eqnarray}
Using the contour integral form of  the ratios conjectures (see,
for example, [9, Lemma 2.1 or 10, Lemma 2.5.1]; the sum of
residues of this integral equals the $\binom{2k}{k}$ terms in the
ratio conjectures as we have previously been writing them, for
example (\ref{eqn:R})), the integral over $t$ is asymptotic to
\begin{eqnarray}\label{eq:kratio} &&
\frac{T^{\tfrac{1}{2} (-\sum_{j=1}^k (\alpha_j-\beta_j))}}
{(w+z)^{k^2}} \;\frac{(-1)^{\frac{2k(2k-1)}{2}}}{(2\pi i)^{2k} k!
k!} \oint\cdots\oint \frac{G(v_1,\ldots,v_{2k})
\Delta(v_1,\ldots,v_{2k})^2}
{\prod_{i=1}^{2k}\prod_{j=1}^{k}(v_i-\alpha_j)(v_i-\beta_j)}
dv_1\cdots dv_{2k}
\end{eqnarray}
where
\begin{eqnarray}
G(v_1,\ldots,v_k,v_{k+1},\ldots,v_{2k})=\frac{\prod_{j=1}^k
(v_j+z)^k\;\;\prod_{j=1}^k(w-v_{j+k})^k
}{D(v_{1+k},\ldots,v_{2k};v_1,\ldots,v_k)}T^{\frac{1}{2}
\sum_{j=1}^k(v_j-v_{j+k})}.
\end{eqnarray}
Here
$D(v_{k+1},\ldots,v_{2k};v_1,\ldots,v_k)=\prod_{j=1}^k\prod_{i=1}^k
(v_j-v_{i+k})$.

Noting the identity
\begin{eqnarray}\label{eqn:id}
\frac{y^{w+z}}{(w+z)^A}=\frac{1}{(A-1)!}\int_0^y u^{w+z}
\left(\log \frac y u\right)^{A-1} \frac {du}{u},
\end{eqnarray}
with  $A=k^2$, we have
\begin{eqnarray}
&&{\mathcal M}_k(\alpha,\beta)\sim\sum_{m,n} \frac{q_{m}m!
p_{n}n!}{\log^{m+n}y}\frac{1}{(2\pi i)^2}\int_{(c_1)} \int_{(c_2)}
\frac{T^{\tfrac{1}{2} (-\sum_{j=1}^k (\alpha_j-\beta_j))}}
{(k^2-1)!} \frac{(-1)^{\frac{2k(2k-1)}{2}}}{(2\pi i)^{2k} k! k!}\nonumber\\
&&\qquad\times \int_0^y\oint\cdots\oint \frac{G(v_1,\ldots,v_{2k})
\Delta(v_1,\ldots,v_{2k})^2}
{\prod_{i=1}^{2k}\prod_{j=1}^{k}(v_i-\alpha_j)(v_i-\beta_j)}
dv_1\cdots dv_{2k}\; u^{w+z} (\log \frac{y}{u})^{k^2-1}
\frac{du}{u} \frac{dw\;dz}{w^{m+1}z^{n+1}}.
\end{eqnarray}

So, focusing on just the integrals over $u$, $w$ and $z$,
\begin{eqnarray}
&&\int_{(c_1)} \int_{(c_2)}\int_0^y \prod_{j=1}^k
(v_j+z)^k\;\;\prod_{j=1}^k(w-v_{j+k})^k  \frac{u^{w+z}(\log
\frac{y}{u})^{k^2-1}} {w^{m+1}z^{n+1}} \frac{du}{u} dw dz\nonumber\\
&&\sim\int_{(c_1)} \int_{(c_2)} \log^{k^2}y \int_0^1
e^{\eta(w+z)\log y} (1-\eta)^{k^2-1} d\eta\frac{\prod_{j=1}^k
(v_j+z)^k\;\;\prod_{j=1}^k(w-v_{j+k})^k }{w^{m+1}z^{n+1}}dw dz,
\end{eqnarray}
where the substitution was $\frac{\log u}{\log y}=\eta$ and we
note that the part of the integral with $u<1$ will not contribute
since for these values of $u$ we can move the path of integration
in $w$ and $z$ as far to the right as we like.  Now let
$y=T^{\theta}$, $\alpha_j=a_j/\log T=a_j\theta/\log y$, and
similarly for $b_j$ and $\beta_j$, as well as making the
replacements $w\rightarrow w/\log y$, $z\rightarrow z/\log y$ and
$v_j\rightarrow v_j/ \log y$.  So we end up with
\begin{eqnarray}
&&{\mathcal M}_k(a/\log T,b/\log T)\sim\sum_{m,n} q_{m}m!
p_{n}n!\frac{1}{(2\pi i)^2} \frac{e^{\tfrac{1}{2} (-\sum_{j=1}^k
(a_j-b_j))}}
{(k^2-1)!} \frac{(-1)^{\frac{2k(2k-1)}{2}}}{(2\pi i)^{2k} k! k!}\nonumber\\
&&\qquad \times\oint\cdots\oint \int_{(c_1)} \int_{(c_2)}\int_0^1
e^{\eta(w+z)}(1-\eta)^{k^2-1} d\eta \frac{\prod_{j=1}^k
(v_j+z)^k\;\;\prod_{j=1}^k(w-v_{j+k})^k
}{D(v_{k+1},\ldots,v_{2k};v_1,\ldots,v_k)w^{m+1}z^{n+1}} dwdz\\
&&\qquad\times\frac{e^{\frac{1}{2\theta}(\sum_{j=1}^k(v_j-v_{j+k}))}\Delta(v_1,\ldots,v_{2k})^2}
{\prod_{i=1}^{2k}\prod_{j=1}^{k}(v_i-\theta a_j)(v_i-\theta b_j)}
dv_1\cdots dv_{2k}.\nonumber
\end{eqnarray}

Now we note that
\begin{eqnarray}
&&\prod_{j=1}^k
(v_j+z)^k\;\;\prod_{j=1}^k(w-v_{j+k})^k\\
&&\qquad=\frac{d^k}{du_1^k}\cdots \frac{d^k}{du_{2k}^k}
e^{u_1(v_1+z)+\cdots
+u_k(v_k+z)+u_{k+1}(w-v_{k+1})+\cdots+u_{2k}(w-v_{2k})}\big|_{u_1=\cdots=u_{2k}=0},\nonumber
\end{eqnarray}
and that
\begin{equation}
\frac{1}{2\pi i} \oint \frac{e^{aw}}{w^{m+1}} dw=\frac{a^m}{m!},
\end{equation}
and use these to write
\begin{eqnarray}
&&{\mathcal M}_k(a/\log T,b/\log T)\sim\frac{e^{\tfrac{1}{2}
(-\sum_{j=1}^k (a_j-b_j))}} {(k^2-1)!}
\frac{(-1)^{\frac{2k(2k-1)}{2}}}{(2\pi i)^{2k} k!
k!}\;\frac{d^k}{du_1^k}\cdots \frac{d^k}{du_{2k}^k} \int_0^1
(1-\eta)^{k^2-1}
\\&&\quad\times Q(\eta+u_{k+1}+\cdots+u_{2k})
P(\eta+u_1+\cdots+ u_k) \;\oint\cdots\oint e^{\frac{1}{2\theta}(\sum_{j=1}^k(v_j-v_{j+k}))}\nonumber\\
&&\quad\times \frac{e^{u_1v_1+\cdots+u_kv_k-u_{k+1}v_{k+1}-\cdots
- u_{2k}v_{2k}}\Delta(v_1,\ldots,v_{2k})^2}
{D(v_{k+1},\ldots,v_{2k};v_1,\ldots,v_k)\prod_{i=1}^{2k}\prod_{j=1}^{k}(v_i-\theta
a_j)(v_i-\theta b_j)} dv_1\cdots
dv_{2k}\;d\eta\big|_{u_1=\cdots=u_{2k}=0}\nonumber
\end{eqnarray}

Now we concentrate on the contour integral over the $v_j$
variables:
\begin{eqnarray}
&&I_v(u_1,\ldots,u_{2k})\nonumber \\
&&\qquad:=\frac{1}{(2\pi i)^{2k} k! k!}\oint\cdots\oint
e^{(\tfrac{1}{2\theta}+u_1)v_1+(\tfrac{1}{2\theta} +u_2)v_2
+\cdots +
(\tfrac{1}{2\theta}+u_k)v_k-(\tfrac{1}{2\theta}+u_{k+1})v_{k+1}-\cdots
-(\tfrac{1}{2\theta}+u_{2k})v_{2k}}\nonumber\\
&&\qquad\qquad\times\frac{\Delta(v_1,\ldots,v_{2k})
\Delta(v_1,\ldots,v_k)\Delta(v_{k+1},\ldots,v_{2k})}{D(\theta
a_1,\ldots,\theta a_k;v_1,\ldots,v_{2k}) D(\theta
b_1,\ldots,\theta b_k;v_1,\ldots,v_{2k})} dv_1\cdots dv_{2k}.
\end{eqnarray}
Expanding the determinants
$\Delta(z_1,\ldots,z_k)=\det[z_j^{m-1}]_{1\leq j,m\leq k}$, we
obtain
\begin{eqnarray}
I_v&=&\frac{1}{(2\pi i)^{2k} k! k!}\oint\cdots\oint
\frac{e^{(\tfrac{1}{2\theta}+u_1)v_1+(\tfrac{1}{2\theta} +u_2)v_2
+\cdots +
(\tfrac{1}{2\theta}+u_k)v_k-(\tfrac{1}{2\theta}+u_{k+1})v_{k+1}-\cdots
-(\tfrac{1}{2\theta}+u_{2k})v_{2k}}}{D(\theta a_1,\ldots,\theta
a_k;v_1,\ldots,v_{2k}) D(\theta
b_1,\ldots,\theta b_k;v_1,\ldots,v_{2k})}\nonumber\\
&&\qquad\times\left( \sum_S {\rm sgn}(S) v_1^{S_0}v_2^{S_1} \cdots
v_k^{S_{k-1}} v_{k+1}^{S_k} \cdots v_{2k}^{S_{2k-1}} \right)
\left( \sum_Q {\rm sgn}(Q) v_1^{Q_0}\cdots v_k^{Q_{k-1}}
\right)\nonumber \\& &\ \ \ \ \times\left( \sum_R {\rm sgn}(R)
v_{k+1}^{R_0}\cdots v_{2k}^{R_{k-1}} \right)  dv_1\cdots dv_{2k}.
\end{eqnarray}
Here $Q$ and $R$ are permutations of $\{0,1,\ldots,k-1\}$ and $S$
is a permutation of $\{0,1,\ldots,2k-1\}$.

Since the integrand is symmetric amongst $v_1,\ldots,v_k$ and also
amongst $v_{k+1},\ldots,v_{2k}$, in each term of the sum over $Q$
we permute the variables $v_1,\ldots,v_k$ so that $v_j$ appears
with the exponent $j-1$, for $j=1,\ldots, k$.  In the sum over $S$
the effect is to redefine the permutations, and the additional
sign involved with this exactly cancels ${\rm sgn(Q)}$.  We do the
same with the sum over $R$, and as a result we are left with
$k!^2$ copies of the sum over the permutation $S$:
\begin{eqnarray}
I_v&=&\frac{1}{(2\pi i)^{2k} }\oint\cdots\oint
\frac{e^{(\tfrac{1}{2\theta}+u_1)v_1+(\tfrac{1}{2\theta} +u_2)v_2
+\cdots +
(\tfrac{1}{2\theta}+u_k)v_k-(\tfrac{1}{2\theta}+u_{k+1})v_{k+1}-\cdots
-(\tfrac{1}{2\theta}+u_{2k})v_{2k}}}{D(\theta a_1,\ldots,\theta
a_k;v_1,\ldots,v_{2k}) D(\theta
b_1,\ldots,\theta b_k;v_1,\ldots,v_{2k})}\nonumber\\
&&\times\sum_S {\rm sgn}(S) v_1^{S_0} v_2^{S_1+1} \cdots
v_k^{S_{k-1}+(k-1)} v_{k+1}^{S_k} v_{k+2}^{S_{k+1}+1} \cdots
v_{2k}^{S_{2k-1}+(k-1)} dv_1\cdots dv_{2k}.
\end{eqnarray}

This can then be written as the following determinant, where we
have written out the $j$th row, with $j=1, \ldots,2k$.
\begin{eqnarray}
I_v&=& \det \bigg\{\frac{1}{2\pi i}
\oint\frac{e^{(\tfrac{1}{2\theta}+u_1)v_1} v_1^{j-1}}{
\prod_{i=1}^k(v_1-\theta a_i)(v_1-\theta b_i)}dv_1,\frac{1}{2\pi
i} \oint\frac{e^{(\tfrac{1}{2\theta}+u_2)v_2} v_2^{j}}{
\prod_{i=1}^k(v_2-\theta a_i)(v_2-\theta b_i)}dv_2,\ldots,\nonumber\\
&&\qquad\qquad\qquad\qquad\ldots,\frac{1}{2\pi i}
\oint\frac{e^{(\tfrac{1}{2\theta}+u_k)v_k} v_k^{j+k-2}}{
\prod_{i=1}^k(v_k-\theta a_i)(v_k-\theta b_i)}dv_k,\\
&&\frac{1}{2\pi i}
\oint\frac{e^{-(\tfrac{1}{2\theta}+u_{k+1})v_{k+1}}
v_{k+1}^{j-1}}{ \prod_{i=1}^k(v_{k+1}-\theta a_i)(v_{k+1}-\theta
b_i)}dv_{k+1}, \frac{1}{2\pi i}
\oint\frac{e^{-(\tfrac{1}{2\theta}+u_{k+2})v_{k+2}} v_{k+2}^{j}}{
\prod_{i=1}^k(v_{k+2}-\theta a_i)(v_{k+2}-\theta b_i)}dv_{k+2},\ldots,\nonumber\\
&&\qquad\qquad\qquad\qquad \ldots,\frac{1}{2\pi i}
\oint\frac{e^{-(\tfrac{1}{2\theta}+u_{2k})v_{2k}} v_{2k}^{j+k-2}}{
\prod_{i=1}^k(v_{2k}-\theta a_i)(v_{2k}-\theta
b_i)}dv_{2k}\bigg\}.\nonumber
\end{eqnarray}

Setting the $a$'s and $b$'s equal to zero and performing the
integration, we have (for integer $n$)
\begin{equation}
\frac{1}{2\pi i} \oint e^{bv} v^n dv = \left\{\begin{array}{cc}0&
n\geq 0\\\frac{b^{-n-1} }{(-n-1)!}&n<0\end{array}\right.,
\end{equation}
giving us
\begin{eqnarray}\label{eq:Iv}
\qquad\quad I_v(\vec{u})&=&\det \left( \begin{array}{cccccc}
\frac{(\frac{1}{2\theta} +u_1)^{2k-1}}{(2k-1)!} &\cdots
&\frac{(\frac{1}{2\theta} +u_k)^{k}}{k!}
&\frac{-(\frac{1}{2\theta} +u_{k+1})^{2k-1}}{(2k-1)!}
&\cdots&\frac{(-1)^k(\frac{1}{2\theta} +u_{2k})^{k}}{k!}
\\\frac{(\frac{1}{2\theta} +u_1)^{2k-2}}{(2k-2)!} &\cdots
&\frac{(\frac{1}{2\theta} +u_k)^{k-1}}{(k-1)!}
&\frac{(\frac{1}{2\theta} +u_{k+1})^{2k-2}}{(2k-2)!}
&\cdots&\frac{(-1)^{k-1}(\frac{1}{2\theta} +u_{2k})^{k-1}}{(k-1)!}
\\\vdots&\ddots&\vdots&\vdots&\ddots&\vdots\\
\frac{(\frac{1}{2\theta} +u_1)^0}{0!}
&\cdots&\frac{(\frac{1}{2\theta} +u_k)^{1-k}}{(1-k)!}
&\frac{(\frac{1}{2\theta} +u_{k+1})^{0}}{0!} &\cdots
&\frac{(-1)^{k-1}(\frac{1}{2\theta} +u_{2k})^{1-k}}{(1-k)!}
\end{array}\right).
\end{eqnarray}
Note that many of the lower matrix entries are zero as we apply
the convention that $\tfrac{1}{n!}=0$ for integer $n<0$.

So, the mollified $k$th moment with the $\alpha$'s and $\beta$'s
set to zero now has the form
\begin{eqnarray}
&&{\mathcal M}_k(0,0)\sim \frac{(-1)^{\frac{2k(2k-1)}{2}}}
{(k^2-1)!} \;\frac{d^k}{du_1^k}\cdots \frac{d^k}{du_{2k}^k}
\int_0^1 (1-\eta)^{k^2-1}  I_v(\vec{u})\\ && \qquad \times
Q(\eta+u_{k+1}+\cdots+u_{2k}) P(\eta+u_1+\cdots+ u_k) \;
d\eta\big|_{u_1=\cdots=u_{2k}=0}.\nonumber
\end{eqnarray}

The differentiation with respect to the $u$ variables is not
difficult, but we will now restrict ourselves to the case $k=2$
where we can write the result fairly concisely.  In this case we
have
\begin{eqnarray}
&&{\mathcal M}_2(0,0)\sim \frac{1} {3!} \;\frac{d^2}{du_1^2}\cdots
\frac{d^2}{du_{4}^2} \int_0^1 (1-\eta)^{3}
\\&&\qquad\times Q(\eta+u_{3}+u_{4})
P(\eta+u_1+u_2)d\eta \;\nonumber\\
&&\qquad\times\left.  \det \left( \begin{array}{cccc}
\frac{(\frac{1}{2\theta} +u_1)^{3}}{3!}  &\frac{(\frac{1}{2\theta}
+u_2)^{2}}{2!} &\frac{-(\frac{1}{2\theta} +u_{3})^{3}}{3!}
&\frac{(\frac{1}{2\theta} +u_{4})^{2}}{2!}
\\\frac{(\frac{1}{2\theta} +u_1)^{2}}{2!}
&(\frac{1}{2\theta} +u_2) &\frac{(\frac{1}{2\theta}
+u_{3})^{2}}{2!} &-(\frac{1}{2\theta} +u_{4})
\\  (\frac{1}{2\theta}+u_1)&1&-(\frac{1}{2\theta}+u_3)&1\\
1 &0 &1  &0
\end{array}\right) \right|_{u_1=\cdots=u_{4}=0}.\nonumber
\end{eqnarray}
Performing the differentiation leads to
\begin{eqnarray}
&&{\mathcal M}_2(0,0)\sim \int_0^1\frac{(1-\eta)^3}{6}
 \bigg(
Q(\eta)P^{(4)}(\eta)+ Q^{(4)}(\eta)P(\eta)
 +4(Q^{(3)}(\eta)P'(\eta)+Q'(\eta)P^{(3)}(\eta)) \nonumber \\
 &&\qquad +6Q''(\eta)P''(\eta)+\frac{2}{\theta}(Q^{(4)}(\eta)P'(\eta)+Q'(\eta)P^{(4)}(\eta))
 +
\frac{8}{\theta}(Q^{(3)}(\eta)P''(\eta)+Q''(\eta)P^{(3)}(\eta))
 \\
&&\qquad
+\frac{2}{\theta^2}(Q^{(4)}(\eta)P''(\eta)+Q''(\eta)P^{(4)}(\eta)
 )+  \frac{4}{\theta^2}Q^{(3)}(\eta)P^{(3)}(\eta) \nonumber\\
 &&\qquad +\frac{2}{3\theta^3} (Q^{(4)}(\eta)P^{(3)}(\eta)+Q^{(3)}(\eta)P^{(4)}(\eta))
 +\frac{1}{12\theta^4} Q^{(4)}(\eta)P^{(4)}(\eta)\bigg)~d\eta .\nonumber \end{eqnarray}
 Integration by parts gives
\begin{theorem}
\label{theo:4zetamoll}
 Assuming the ratios conjecture as indicated in (\ref{eq:kratio}), if $Q$
and $P$ are polynomials which vanish at 0 and whose first three
derivatives vanish at 0, then for any $\theta>0$ we have
\begin{eqnarray}\label{eq:notdavid}
&&\frac{1}{T}\int_0^T|\zeta(1/2+it)|^4 M_2(1/2+it,Q)M_2(1/2-it,P)~dt  \nonumber\\
&&\quad = P(1)Q(1) +\frac{1}{\theta}\int_0^1\frac{(1-\eta)^3}{6}
 \bigg(
2(Q^{(4)}(\eta)P'(\eta)+Q'(\eta)P^{(4)}(\eta))
 \\
&&\qquad +
8(Q^{(3)}(\eta)P''(\eta)+Q''(\eta)P^{(3)}(\eta))
+\frac{2}{\theta}(Q^{(4)}(\eta)P''(\eta)+Q''(\eta)P^{(4)}(\eta)
 )+  \frac{4}{\theta}Q^{(3)}(\eta)P^{(3)}(\eta) \nonumber\\
 &&\qquad +\frac{2}{3\theta^2} (Q^{(4)}(\eta)P^{(3)}(\eta)+Q^{(3)}(\eta)P^{(4)}(\eta))
 +\frac{1}{12\theta^3} Q^{(4)}(\eta)P^{(4)}(\eta)\bigg)~d\eta +O(1/\log T).\nonumber
\end{eqnarray}
\end{theorem}
\begin{remark}
While we don't know which are the minimizing polynomials,
 with the choice
that $P(x)= Q(x)=x^4$, the right side of (\ref{eq:notdavid}) is
equal to
\begin{equation}
1 +\frac{208}{35\theta} + \frac{48}{5\theta^2} +
\frac{32}{5\theta^3} + \frac{2}{ \theta^4}.
\end{equation}
\end{remark}

By a similar calculation one can show that with polynomials $P_i$
and $Q_j$ satisfying $P_i(0)=P_i'(0)=Q_j(0)=Q_j'(0)=0$, we have
\begin{eqnarray}\label{eq:david}&&\frac{1}{T}\int_0^T
|\zeta(1/2+it)|^4 M_1(s,P_1)M_1(s,P_2)M_1(1-s,Q_1)M_1(1-s,Q_2)
~dt\nonumber \\
&& \qquad = \frac{1} {16} \;\frac{d}{du_1}\cdots \frac{d}{du_{4}}
\frac{d}{dU_1}\cdots \frac{d}{dU_{4}} \iiiint _{\mathcal{R}}
P_1(\tfrac{\eta_1}{2}+\tfrac{\eta_2}{2}+u_3+u_4)
P_2(\tfrac{\eta_3}{2}+\tfrac{\eta_4}{2}
+ U_3+U_4)\nonumber \\
&&\qquad\qquad\qquad\qquad \times
Q_1(\tfrac{\eta_1}{2}+\tfrac{\eta_3}{2}+u_1+u_2)
Q_2(\tfrac{\eta_2}{2}+\tfrac{\eta_4}{2}+U_1+U_2)d\eta_1\cdots d\eta_4 \;\\
&&\qquad\qquad\qquad\qquad\times  \left.
I_v(u_1+U_1,u_2+U_2,u_3+U_3,u_4+U_4)
\right|_{u_1=\cdots=u_4=U_1=\cdots=U_{4}=0}\nonumber
\\
&&\qquad \qquad \qquad \qquad \qquad +O(1/\log
T),\nonumber\end{eqnarray} where $\mathcal{R}$ is the subset of
$[-1,1]^4$ for which $\eta_1+\eta_2\geq0$, $\eta_3+\eta_4\geq0$,
 $\eta_1+\eta_3\geq0$ and $\eta_2+\eta_4\geq0$.

\begin{remark}
While we don't know which are the minimizing polynomials,
 with the choice
that $P_1(x)=P_2(x)=Q_1(x)=Q_2(x)=x^2$, the right side of
(\ref{eq:david}) is equal to
\begin{equation}
1 +\frac{68}{21\theta} + \frac{10}{3\theta^2} +
\frac{64}{45\theta^3} + \frac{2}{9 \theta^4}.
\end{equation}
\end{remark}

\subsection{A third power mollification}
In the work of Hughes alluded to earlier it is likely that in
applications to zeros of $\zeta(s)$ on the critical line the
moment
\begin{eqnarray}
\int_0^T |\zeta(1/2+it)|^2
|M_1(1/2+it,P_1)+\zeta(1/2+it)M_2(1/2+it,P_2)|^2~dt
\end{eqnarray}
will need to be evaluated.  Here, as at (\ref{eq:M1}) and
(\ref{eq:M2}), $M_1$ is a mollifier with arithmetic coefficients
$\mu(m)$ smoothed by $P_1$ and
\begin{eqnarray}
M_2(s,P_2):= \sum_{m\le y_2} \frac{\mu_2(m)P_2\left(\frac{\frac
{y_2}{m}}{\log y_2}\right)}{m^s},
\end{eqnarray}
where the $\mu_2$ are the coefficients of $1/\zeta(s)^2$. In order
to evaluate this integral we can use the results
(\ref{eqn:mollunit}) and Theorem \ref{theo:4zetamoll}; in addition
to these we need to evaluate
\begin{eqnarray} \label{eq:I3}
I_3(\alpha,\beta,\gamma;P_1,P_2):=\int_0^T
\zeta(s+\alpha)\zeta(s+\beta)\zeta(1-s+\gamma)
M_1(1-s,P_1)M_2(s,P_2)~dt.
\end{eqnarray}
Proceeding as usual,
\begin{eqnarray}
I_3 &= &\sum_{m,n}\frac{p_{1,m}m!p_{2,n}n!}{\log^m y_1 \log ^n
y_2} \frac{1}{(2\pi i)^2}
\int_{(c_1)}\int_{(c_2)}\frac{y^{w+z}}{w^{m+1}z^{n+1}} \\
&&\qquad\qquad\times \int_0^T
\frac{\zeta(s+\alpha)\zeta(s+\beta)\zeta(1-s+\gamma)}{\zeta(s+w)^2\zeta(1-s+z)}~dt~dw~dz
\nonumber
\end{eqnarray}
 By the ratios
conjecture,
\begin{eqnarray}
&&I_3\sim T\sum_{m,n}\frac{p_{1,m}m!p_{2,n}n!}{\log^m y_1 \log ^n
y_2} \frac{1}{(2\pi i)^2}
\int_{(c_1)}\int_{(c_2)}\frac{y_1^wy_2^z}{w^{m+1}z^{n+1}(w+z)^2}
\bigg(\frac{(\alpha+z)(\beta+z)(\gamma+w)^2}{(\alpha+\gamma)(\beta+\gamma) }\\
&&\qquad -T^{-\alpha-\gamma}
\frac{(-\gamma+z)(\beta+z)(-\alpha+w)^2}{(\alpha+\gamma)(\beta-\alpha)
} -T^{-\beta-\gamma}
\frac{(\alpha+z)(-\gamma+z)(-\beta+w)^2}{(\alpha-\beta)(\beta+\gamma)
}\bigg)~dw~dz\nonumber
\end{eqnarray}
Taking the limit (Mathematica can be helpful here) as $\alpha,
\beta, \gamma\to 0$, we have
\begin{eqnarray}
&&I_3\sim T\sum_{m,n}\frac{p_{1,m}m!p_{2,n}n!}{\log^m y_1 \log ^n
y_2} \frac{1}{(2\pi i)^2}
\int_{(c_1)}\int_{(c_2)}\frac{y_1^zy_2^w}{z^{m+1}w^{n+1}(w+z)^2}\\
&&\qquad  \big((w+z)^2+w^2 zL+2w z^2 L +w^2z^2 L^2
/2\big)~dw~dz.\nonumber
\end{eqnarray}
The contribution of the $(w+z)^2$ term at the beginning of the
brackets above is $P_1(1)P_2(1)$, essentially using
(\ref{eqn:poly}). For the rest of the terms, we write, using
(\ref{eqn:id}),
\begin{eqnarray}
\frac{y_1^zy_2^w}{(w+z)^2}=
\frac{y_2^{w+z}(y_1/y_2)^z}{(w+z)^2}=(y_1/y_2)^z\int_0^{y_2}u^{w+z}
\log \frac {y_2}{ u} \frac{du}{u}.
\end{eqnarray}
The integral over $w$ is 0 unless $u\ge 1$; the integral over $z$
is 0 unless $uy_1/y_2>1$; this inequality is weaker than the
requirement that $u>1$, since, in general $y_1>y_2$.   In this
way, we see that
\begin{eqnarray}
&&I_3\sim T\bigg(P_1(1)P_2(1)+\int_1^{y_2}\bigg(\frac{L}{\log
y_1\log^2 y_2}P_1'\left(\frac{\log \frac{y_1u}{y_2}}{\log
y_1}\right)P_2''\left(\frac{\log u}{\log
y_2}\right)\\
&&\qquad \quad +\frac{2L}{\log^2 y_1 \log
y_2}P_1''\left(\frac{\log \frac{y_1u}{y_2}}{\log
y_1}\right)P_2'\left(\frac{\log u}{\log
y_2}\right) \nonumber\\
&&\qquad \quad+\frac {L^2}{2\log^2 y_1 \log^2y_2}
P_1''\left(\frac{\log \frac{y_1u}{y_2}}{\log
y_1}\right)P_2''\left(\frac{\log u}{\log y_2}\right)\bigg)\log
\frac {y_2}{u} ~du\bigg).\nonumber
\end{eqnarray}
With a change of variables $u=y_2^\eta$, and with
$y_1=T^{\theta_1}$, $y_2=T^{\theta_2}$, we see that
\begin{theorem} Let $I_3$ be as defined in (\ref{eq:I3}). Assuming the ratios conjecture as indicated in (\ref{eq:kratio}), if $P_1$
and $P_2$ are polynomials which vanish at 0 and whose first
derivatives vanish at 0, then for any $\theta>0$ we have
\begin{eqnarray}
&&I_3(0,0,0;P_1,P_2)=
T\bigg(P_1(1)P_2(1)+\int_0^{1}\bigg(\frac{1}{\theta_1}
P_1'\big(1+(1-\eta) \theta_2/\theta_1\big)P_2''(\eta)\\
&& \nonumber\qquad \quad +\frac{2
\theta_2}{\theta_1^2}P_1''\big(1+(1-\eta)\theta_2/\theta_1
\big)P_2'(\eta) \\
&& \nonumber\qquad \quad+\frac {1}{2\theta_1^2}
P_1''\big(1+(1-\eta)\theta_2/\theta_1 \big)P_2''(\eta)
\bigg)(1-\eta)~d\eta+O(1/\log T)\bigg).
\end{eqnarray}
\end{theorem}

\section{Discrete moments of the Riemann zeta
function and its derivatives}\label{sect:discrete}

So far in this paper we've considered integer moments. Another
kind of average which gives useful information about the
distribution of zeros is a discrete moment summing the zeta
function, or its derivatives, at or near the zeros.

In the 80's Gonek [20], assuming the Riemann Hypothesis,  proved,
amongst much more general results, that
\begin{equation}\label{eqn:gonek}
\sum_{1\leq\gamma \leq T} |\zeta'(\rho)|^2=\frac{T}{24 \pi} \log^4
T +O(T\log^3T),
\end{equation}
where $\rho=1/2 +i\gamma$ is a zero of the Riemann zeta function.

Hughes, Keating and O'Connell used the analogy with random matrix
theory to propose the following:
\begin{conjecture}[{Hughes, Keating and O'Connell [23]}] \label{conj:hko}
For $k>-3/2$ and bounded,
\begin{equation}
\sum_{0<\gamma_n\leq
T}|\zeta'(\tfrac{1}{2}+i\gamma_n)|^{2k}\sim\frac{T}{2\pi}\frac{G^2(k+2)}{G(2k+3)}
a(k)(\log T)^{k(k+2)+1}
\end{equation}
as $T\rightarrow\infty$, with
\begin{equation}
a(k)=\prod_p (1-\tfrac{1}{p})^{k^2}
\sum_{m=0}^{\infty}\left(\frac{\Gamma(m+k)}{m!\Gamma(k)}\right)^2p^{-m}.
\end{equation}
Here $G(k)$ is the Barnes $G$-function.
\end{conjecture}

While the above conjecture produces the leading order terms, it
was previously not known how to obtain lower order terms for the
moments of $|\zeta'(\rho)|$. In the first two subsections we
consider the second and fourth moments of $|\zeta'(\rho)|$, and in
the last subsection the second moment of $|\zeta(\rho+a)|$.  Using
the ratios conjecture we show how to obtain all of the lower-order
terms for these averages.

Note that Hughes [22] has conjectured using random matrix theory
\begin{conjecture}
\begin{eqnarray}
&&\sum_{0<\gamma_n\leq T}|\zeta(\rho+2\pi i\alpha
\log^{-1}\tfrac{T}{2\pi}))|^{2k} \sim\frac{T}{2\pi} G^2(k+1) /
G(2k+1) a(k) F_k(2\pi\alpha)  (\log T)^{k^2+1},
\end{eqnarray}
where
\begin{eqnarray}
F_k(2x)  = \tfrac{\pi}{2}( x J_{k+1/2}(x)^2 + x J_{k-1/2}(x)^2 - 2
k J_{k+1/2}(x) J_{k-1/2}(x) ) .
\end{eqnarray}
Here $J_k$ is the usual Bessel function.
\end{conjecture}

\subsection{Second moment of the derivative: $\sum_{1\leq\gamma \leq T}
|\zeta'(\rho)|^2$}\label{sect:2ndderiv}

 In this section we will show how to use
the ratios conjecture to reproduce the result of Gonek
(\ref{eqn:gonek}) and to derive all the lower order terms.

The first step is to write the sum over zeros as a contour
integral
\begin{equation}
\sum_{\gamma< T}|\zeta'(\rho)|^2=\frac{1}{2\pi
i}\int_{C}\frac{\zeta'(z)}{\zeta(z)} \zeta'(z)\zeta'(1-z)dz,
\end{equation}
where the contour $C$ has corners $c$, $c+iT$, $1-c+iT$ and $1-c$,
with $1/2<c<1$. The integrals along the horizontal sides of this
rectangle can be neglected, and so we look first at
\begin{eqnarray}
I_r&=&\frac{1}{2\pi i} \int_{c}^{c+iT}\frac{\zeta'(z)}{\zeta(z)}
\zeta'(z)\zeta'(1-z)dz\\
&=&\frac{1}{2\pi} \int_0^T \frac{\zeta'(c+it)}{\zeta(c+it)}
\zeta'(c+it)\zeta'(1-c-it)dt\nonumber\\
&=&\frac{d}{d\beta}\frac{d}{d\gamma}\frac{d}{d\delta}
\frac{1}{2\pi}\int_0^T\frac{\zeta(c+it+\beta)}{\zeta(c+it)}\nonumber\\
&&\quad\quad\quad \times
\zeta(c+it+\gamma)\zeta(1-c-it+\delta)dt\;
\bigg|_{\beta=\gamma=\delta=0}.\nonumber
\end{eqnarray}

Now we follow the recipe for computing the ratio of zeta functions
in the integrand.  As in Section \ref{sec:intro}, we replace zeta
functions in the numerator by
\begin{equation}
\zeta(s)\sim \sum_{n\leq \sqrt{\tfrac{t}{2\pi}}} \frac{1}{n^s}
+\chi(s)\sum_{n\leq \sqrt{\tfrac{t}{2\pi}}} \frac{1}{n^{1-s}}
\end{equation}
and zeta functions in the denominator by
\begin{equation}
\frac{1}{\zeta(s)} = \sum_{n=1}^{\infty} \frac{\mu(n)}{n^s}.
\end{equation}
Since each of the zeta functions in the numerator is replaced by
two sums, when multiplied out we have a total of eight terms in
the integral.  Considering the term involving the first term from
each approximate functional equation, the resulting sum is
\begin{eqnarray}
&&\sum_{hmn=\ell}
\frac{\mu(h)}{m^{1/2+\beta}n^{1/2+\gamma}h^{1/2}\ell^{1/2+\delta}}
\\
&&\quad=\sum_{h,m,n}\frac{\mu(h)}{m^{1+\beta+\delta}n^{1+\gamma+\delta}
h^{1+\delta}} \nonumber\\
&&\quad=
\frac{\zeta(1+\beta+\delta)\zeta(1+\gamma+\delta)}{\zeta(1+\delta)}.\nonumber
\end{eqnarray}

Of the eight terms in the integrand, only those with the same
number of $\chi$ factors resulting from $\zeta(z)$ as from
$\zeta(1-z)$ will survive.  This means the recipe implies two
further terms, and using (\ref{eqn:chi}), we have
\begin{eqnarray}\label{eqn:IR}
I_r&=& \frac{d}{d\beta} \frac{d}{d\gamma} \frac{d}{d\delta}
\frac{1}{2\pi} \int_0^T \left(\frac{\zeta(1+\beta+\delta)
\zeta(1+\gamma+\delta)} {\zeta(1+\delta)}+
(\frac{t}{2\pi})^{-\beta-\delta}\frac{ \zeta(1-\delta-\beta)
\zeta(1+\gamma-\beta)}{\zeta(1-\beta)}\right.\nonumber\\
&&\quad\quad \left.+(\frac{t}{2\pi})^{-\gamma-\delta}
\frac{\zeta(1+\beta-\gamma) \zeta(1-\delta-\gamma)}
{\zeta(1-\gamma)}\right)(1+O(t^{-\tfrac{1}{2}+\epsilon})) dt
\bigg|_{\beta=\gamma=\delta=0}.
\end{eqnarray}

We now consider the contribution from the other side of the
contour of integration:
\begin{equation}
I_{\ell}= \frac{1}{2\pi i}\int_{1-c+iT}^{1-c}\frac{\zeta'(z)}{
\zeta(z)} \zeta'(z) \zeta'(1-z)dz.
\end{equation}
Replacing $z$ with $1-z$, we have
\begin{equation}
I_{\ell}=\frac{-1}{2\pi i} \int_{c-iT}^{c} \frac{\zeta'(1-z)}
{\zeta(1-z)} \zeta'(1-z)\zeta'(z)dz.
\end{equation}
Differentiating the functional equation gives us
\begin{equation}
\frac{\zeta'}{\zeta}(1-z)=\frac{\chi'}{\chi}(1-z)-\frac{\zeta'}{\zeta}
(z),
\end{equation}
and so
\begin{eqnarray}
I_{\ell}&=&\frac{-1}{2\pi i} \int_{c-iT}^{c}
\frac{\chi'(1-z)}{\chi(1-z)} \zeta'(1-z)\zeta'(z)dz\\
&&\quad\quad+\frac{1}{2\pi
i}\int_{c-iT}^{c}\frac{\zeta'(z)}{\zeta(z)}\zeta'(z)\zeta'(1-z)
dz.\nonumber
\end{eqnarray}
Here the first integral can be shifted over to the
$\tfrac{1}{2}$-line, while the second one is just the complex
conjugate of $I_r$, already calculated in (\ref{eqn:IR}).  In
addition, we can use
\begin{equation}
\chi(s)=(\frac{t}{2\pi})^{1/2-s}e^{it+\pi i/4}(1+O(\tfrac{1}{t}))
\end{equation}
to approximate $\chi'/\chi(1-z)$ with $-\log \tfrac{t}{2\pi}$.
Thus we have
\begin{eqnarray}
\sum_{\gamma< T}|\zeta'(\rho)|^2&=&2I_r+\frac{1}{2\pi}\int_0^T\log
\tfrac{t}{2\pi}\;|\zeta'(1/2+it)|^2 (1+O(t^{-1}))dt\nonumber\\
&=&2I_r+\frac{d}{d\alpha}\frac{d}{d\beta}\frac{1}{2\pi}\int_0^T
\log\tfrac{t}{2\pi}
\;\zeta(1/2+it+\alpha)\zeta(1/2-it+\beta)(1+O(t^{-1}))
dt\bigg|_{\alpha=\beta=0}\\
&=&2I_r+\frac{d}{d\alpha}\frac{d}{d\beta}\frac{1}{2\pi}\int_0^T
\log\tfrac{t}{2\pi}
\;\left(\zeta(1+\alpha+\beta)+(\frac{t}{2\pi})^{-\alpha-\beta}\zeta(1-\alpha-\beta)\right)\nonumber\\
&&\quad\quad\quad(1+O(t^{-1/2+\epsilon}))
dt\bigg|_{\alpha=\beta=0},\nonumber
\end{eqnarray}
where the last line is a further application of the ratios
conjecture (or in this case the simpler moment conjecture [10])
similar to that in Section \ref{sec:intro}.

Using (\ref{eqn:IR}) for $I_r$, it is now necessary to carry out
the differentiation and take the limits as $\alpha, \beta,\gamma$
and $\delta$ tend to zero.  This results in a polynomial in $\log
\tfrac{t}{2\pi}$.  If we write
\begin{equation}
\zeta(1+s)=\frac{1}{s}+\gamma-\gamma_1
s+\frac{\gamma_2}{2!}s^2-\frac{\gamma_3}{3!}s^3\cdots,
\end{equation}
then the final result is
\begin{theorem}
Assuming the ratio conjecture as indicated in (\ref{eqn:IR}), we
have
\begin{eqnarray}
\sum_{\gamma< T}|\zeta'(\rho)|^2&=&\int_0^T \left(\frac{1}{24\pi}
\log^4 \frac{t}{2\pi}
+\frac{\gamma}{3\pi}\log^3\frac{t}{2\pi}+(\frac{\gamma^2}{2\pi}-\frac{\gamma_1}{\pi})
\log^2\frac{t}{2\pi}-(\frac{\gamma^3}{\pi}+\frac{5\gamma
\gamma_1}{\pi} +
\frac{\gamma_2}{2\pi})\log\frac{t}{2\pi}\right.\nonumber\\
&&\quad\quad\quad\left.+\frac{\gamma^4}{\pi} +
\frac{6\gamma^2\gamma_1}{\pi} + \frac{7\gamma_1^2}{\pi}+\frac{4
\gamma\gamma_2}{\pi}+ \frac{5\gamma_3}{3\pi} \right)(1+O(t^{-1/2+\epsilon}))dt\nonumber\\
&=&\frac{T}{24 \pi} \log^4 T +O(T\log^3T).
\end{eqnarray}
\end{theorem}

\begin{remark}
The   leading order term of the above agrees with Gonek's result
(\ref{eqn:gonek}).  It is possible that Gonek's methods could
prove the theorem conditional only on the Riemann Hypothesis.
Also, Pokharel and Rubinstein have checked this numerically.
\end{remark}
\begin{remark}
Since the original version of this paper appeared on the archive,
Milinovich has used Gonek's method to verify all the main terms
above.  He also remarks that this result can probably be obtained
from a theorem of Fujii.
\end{remark}
\subsection{Fourth moment of the derivative:  $\sum_{\gamma<
T}|\zeta'(\rho)|^4$}

Higher moments are more difficult because of complicated
arithmetic contributions.  Unlike the case of the fourth moment of
the zeta function itself,
\begin{eqnarray}\label{eq:4thzeta}
&&\frac{1}{T} \int_0^T|\zeta(\tfrac{1}{2}+it)|^{4}dt\\
&&\qquad = \frac{1}{T}\int_0^T \frac{1}{2\pi^2}
\log^4\tfrac{t}{2\pi} + \frac{8}{\pi^4}
    \left(\gamma{\pi }^{2}-3\zeta'(2)\right) \log^3\tfrac{t}{2\pi}  \nonumber \\&   &\mbox{} +\frac{6}{\pi^6}
    \left(
    -48\gamma\zeta'(2){\pi }^{2}-12\zeta''(2){\pi }^{2}+7\gamma^{2}{\pi }^{4}+
    144\zeta'(2)^{2}-2\,\gamma_1{\pi }^{4}\right) \log^2\tfrac{t}{2\pi}  \nonumber \\&   &\mbox{} +\frac{12}{\pi^8}
    \biggl(
    6\gamma^{3}\pi^{6}-84\gamma^{2}\zeta'(2)\pi^{4}+24\gamma_1\zeta'(2)
    \pi^{4}-1728\zeta'(2)^{3}+
    576\gamma\zeta'(2)^{2}\pi^{2} \nonumber \\&   &\phantom{TTTTT}+288\zeta'(2)\zeta''(2)\pi^{2}-
    8\zeta'''(2)\pi^{4}-10\gamma_1\gamma\pi^{6}-
    \gamma_2\pi^{6}-48\gamma\zeta''(2)\pi^{4}
    \biggr)
    \log\tfrac{t}{2\pi}  \nonumber \\&   &\mbox{} +\frac{4}{\pi^{10}}
    \biggl(
    -12\zeta''''(2){\pi }^{6}+36\gamma_2\zeta'(2){\pi }^{6}+9{\gamma}^{4}{\pi }^{8}+
    21\gamma_1^{2}{\pi }^{8}+432\zeta''(2)^{2}{\pi }^{4} \nonumber \\&   &\phantom{TTTTT}+3456
    \gamma\zeta'(2)\zeta''(2){\pi }^{4}+3024{\gamma}^{2}\zeta'(2)^{2}{\pi }^{4}-
    36{\gamma}^{2}\gamma_1{\pi }^{8}-252{\gamma}^{2}\zeta''(2){\pi }^{6}\nonumber \\&   &\phantom{TTTTT}+3\gamma\gamma_2{\pi}^{8}+
    72\gamma_1\zeta''(2){\pi }^{6}+360\gamma_1\gamma\zeta'(2){\pi }^{6}-216{\gamma}^{
    3}\zeta'(2){\pi }^{6}\nonumber \\&   &\phantom{TTTTT}-864\gamma_1\zeta'(2)^{2}{\pi }^{4}+5\gamma_3{\pi }^{8} +
    576\zeta'(2)\zeta'''(2){\pi }^{4}-20736\gamma\zeta'(2)^{3}{\pi }^{2} \nonumber \\&   &\phantom{TTTTT}-
    15552\zeta''(2)\zeta'(2)^{2}{\pi }^{2}-96\gamma\zeta'''(2){\pi }^{6}+62208
    \zeta'(2)^{4}
    \biggr)dt+o(1)\nonumber\\
    &&\qquad = \frac{1}{12\zeta(2)} \log ^4 T + O(\log^3
    T)\label{eq:zeta2}
\end{eqnarray}
(from the moment conjecture formula of [10] and also implied by
[21]), the fourth moment of the modulus of the derivative involves
arithmetic factors that are more complicated than derivatives of
the zeta function at 2.  In the following we calculate the first
four leading order terms, demonstrating where the first of these
new arithmetic factors appears.

As for the second moment, the first step is to write the sum over
zeros as a contour integral
\begin{equation}
\sum_{\gamma< T}|\zeta'(\rho)|^4=\frac{1}{2\pi
i}\int_{C}\frac{\zeta'(z)}{\zeta(z)}
\zeta'(z)\zeta'(1-z)\zeta'(z)\zeta'(1-z)dz,
\end{equation}
with the contour $C$ running from $c$ to $c+iT$, $1-c+iT$ and
$1-c$. The horizontal integrals don't contribute significantly,
and so we define
\begin{eqnarray}
I_R
&=&\frac{d}{d\alpha}\frac{d}{d\beta}\frac{d}{d\gamma}\frac{d}{d\delta}\frac{d}{d\epsilon}
\frac{1}{2\pi}\int_0^T\frac{\zeta(c+it+\alpha)}{\zeta(c+it)}\nonumber\\
&&\quad\quad\quad \times
\zeta(c+it+\beta)\zeta(c+it+\gamma)\zeta(1-c-it+\delta)\zeta(1-c-it+\epsilon)dt\;
\bigg|_{\alpha=\beta=\gamma=\delta=\epsilon=0}.
\end{eqnarray}

The sum resulting from taking the first half of each approximate
functional equation is
\begin{eqnarray}\label{eq:bigsum}
&&\sum_{m_1m_2m_3h=n_1n_2}
\frac{\mu(h)}{m_1^{1/2+\alpha}m_2^{1/2+\beta}m_3^{1/2+\gamma}h^{1/2}n_1^{1/2+\delta}n_2^{1/2+\epsilon}}
.
\end{eqnarray}
Here we note that if we let $\gamma=0$ we obtain the sum
\begin{eqnarray}
\label{eq:gamma0}
\prod_p\sum_{m_1+m_2+c=n_1+n_2}\frac{1}{p^{(1/2+\alpha)m_1+(1/2+\beta)m_2+(1/2)c
+(1/2+\delta)n_1+(1/2+\epsilon)n_2}}\;\sum_{m_3+h=c} \mu(p^h).
\end{eqnarray}
The final sum is zero unless $c=0$, and thus the whole expression
reduces to the corresponding sum that results from applying the
ratios (or moments) conjecture to the fourth moment of zeta,
$\frac{1}{T}\int_{0}^T
\zeta(1/2+it+\alpha)\zeta(1/2+it+\beta)\zeta(
1/2-it+\delta)\zeta(1/2-it+\epsilon)dt$, which itself produces the
arithmetic contribution $1/\zeta(2+\alpha+\beta+\delta+\epsilon)$
observed as the factor $1/\zeta(2)$ in the leading order term of
(\ref{eq:zeta2}) (see for example [10]).

We keep this in mind as we continue with the sum in
(\ref{eq:bigsum}).  This sum can be written as an Euler product,
and we can pull out the divergent terms in the form of zeta
functions:
\begin{eqnarray}
\label{eq:T}
&&T(\alpha,\beta,\gamma,\delta,\epsilon):=\frac{\zeta(1+\alpha+\delta)
\zeta(1+\alpha+\epsilon)\zeta(1+\beta+\delta)
\zeta(1+\beta+\epsilon)\zeta(1+\gamma+\delta)\zeta(1+\gamma+\epsilon)}{
\zeta(1+\delta)\zeta(1+\epsilon)}\nonumber\\&&\qquad\times \prod_p
\left[
\frac{(1-\tfrac{1}{p^{1+\alpha+\delta}})(1-\tfrac{1}{p^{1+\alpha+\epsilon}})
(1-\tfrac{1}{p^{1+\beta+\delta}})(1-\tfrac{1}{p^{1+\beta+\epsilon}})
(1-\tfrac{1}{p^{1+\gamma+\delta}})(1-\tfrac{1}{p^{1+\gamma+\epsilon}})}
{(1-\tfrac{1}{p^{1+\delta}})(1-\tfrac{1}{p^{1+\epsilon
}})}\right.\\
&&\left.\qquad\times \sum_{m_1+m_2+m_3+h=n_1+n_2}
\frac{\mu(p^h)}{p^{(1/2+\alpha)m_1+(1/2+\beta)m_2+(1/2+\gamma)m_3+h/2+(1/2+\delta)n_1
+(1/2+\epsilon)n_2}}\right].\nonumber
\end{eqnarray}
With the remaining terms from the approximate functional equations
we have
\begin{eqnarray}\label{eq:tenterms}
I_R&=&\frac{d}{d\alpha}\frac{d}{d\beta}\frac{d}{d\gamma}\frac{d}{d\delta}\frac{d}{d\epsilon}
\frac{1}{2\pi}\int_0^T
T(\alpha,\beta,\gamma,\delta,\epsilon)+(\frac{t}{2\pi})^{-\alpha-\delta}
T(-\delta,\beta,\gamma,-\alpha,\epsilon)\nonumber\\
&&+ (\frac{t}{2\pi})^{-\alpha-\epsilon}
T(-\epsilon,\beta,\gamma,\delta,-\alpha)+
(\frac{t}{2\pi})^{-\beta-\delta}
T(\alpha,-\delta,\gamma,-\beta,\epsilon)+(\frac{t}{2\pi})^{-\beta-\epsilon}
T(\alpha,-\epsilon,\gamma,\delta,-\beta)\nonumber\\
&& +(\frac{t}{2\pi})^{-\gamma-\delta}
T(\alpha,\beta,-\delta,-\gamma,\epsilon)+
(\frac{t}{2\pi})^{-\gamma-\epsilon}
T(\alpha,\beta,-\epsilon,\delta,-\gamma)\nonumber\\
&& +
(\frac{t}{2\pi})^{-\alpha-\beta-\delta-\epsilon}T(-\delta,-\epsilon,\gamma,
-\alpha, -\beta)+(\frac{t}{2\pi})^{-\alpha-\gamma-\delta-\epsilon}
T(-\delta,\beta,-\epsilon,-\alpha,-\gamma)\nonumber \\
&& +(\frac{t}{2\pi})^{-\beta-\gamma-\delta-\epsilon}
T(\alpha,-\delta,-\epsilon,-\beta,-\gamma)dt\;
\bigg|_{\alpha=\beta=\gamma=\delta=\epsilon=0}+O(T^{1/2+\epsilon}).
\end{eqnarray}
The most concise way to write the ten terms in the integrand above
is as a contour integral (as described in [10]).  So, we have that
\begin{eqnarray}\label{eq:IRa}
I_R&=&\frac{d}{d\alpha}\frac{d}{d\beta}\frac{d}{d\gamma}\frac{d}{d\delta}\frac{d}{d\epsilon}
e^{\tfrac{\log
t/(2\pi)}{2}(-\alpha-\beta-\gamma-\delta-\epsilon)}\frac{1}{2\pi}\int_0^T
\frac{1}{3!\;2!(2\pi i)^{5}}\oint\cdots\oint e^{\tfrac{\log t/(2\pi)}{2}(z_1+z_2+z_3-z_4-z_5)}\nonumber\\
&&\qquad\qquad\times
\frac{T(z_1,z_2,z_3,-z_4,-z_5)\Delta(z_1,\ldots,z_5)^2}{
\prod_{j=1}^5(z_j-\alpha)(z_j-\beta)(z_j-\gamma)(z_j+\delta)(z_j+\epsilon)}
dz_1\cdots dz_5\; dt+O(T^{1/2+\epsilon}).
\end{eqnarray}

With the formulae
\begin{equation}\label{eq:deriv1}
\frac{d}{d\alpha}\frac{e^{-a\alpha}}
{\prod_{j=1}^n(z_j-\alpha)}\Bigg|_{\alpha=0}=\frac{1}{\prod_{j=1}^n
z_j}\left(\sum_{j=1}^n\frac{1}{z_j}-a\right)
\end{equation}
and
\begin{equation}\label{eq:deriv2}
\frac{d}{d\delta}\frac{e^{-a\delta}}
{\prod_{j=1}^n(z_j+\delta)}\Bigg|_{\delta=0}=\frac{1}{\prod_{j=1}^n
z_j}\left(-\sum_{j=1}^n\frac{1}{z_j}-a\right)
\end{equation}
the differentiation can easily be performed, to yield
\begin{eqnarray}\label{eq:IR}
I_R&=& \frac{1}{2\pi}\int_0^T\frac{1}{3!\;2!(2\pi
i)^{5}}\oint\cdots\oint
\frac{T(z_1,z_2,z_3,-z_4,-z_5)\Delta(z_1,\ldots,z_5)^2}{
\prod_{j=1}^5z_j^5}\nonumber\\
&&\quad\times\left(-\tfrac{\log t/(2\pi)}{2}-\sum_{j=1}^5
\frac{1}{z_j}\right)^2\left(-\tfrac{\log t/(2\pi)}{2}+\sum_{j=1}^5
\frac{1}{z_j}\right)^3\nonumber \\
&&\quad\times e^{\tfrac{\log
t/(2\pi)}{2}(z_1+z_2+z_3-z_4-z_5)}dz_1\cdots
dz_5\;dt+O(T^{1/2+\epsilon}).
\end{eqnarray}

For the contribution from the other side of the contour of
integration we have, following the same method as for the second
moment,
\begin{eqnarray}
I_{L}&= &\frac{1}{2\pi i}\int_{1-c+iT}^{1-c}\frac{\zeta'(z)}{
\zeta(z)} \zeta'(z) \zeta'(1-z) \zeta'(z) \zeta'(1-z)dz\nonumber\\
&=&\frac{1}{2\pi}\int_0^T \bigg(\log
\frac{t}{2\pi}+O(1/(t+1))\bigg)
\;|\zeta'(\tfrac{1}{2} +it)|^4dt+I_R\\
&=&I_R+\frac{d}{d\alpha}\frac{d}{d\beta}\frac{d}{d\gamma}
\frac{d}{d\delta}\frac{1}{2\pi}\int_0^T \bigg(\log
\frac{t}{2\pi}+O(1/(t+1))\bigg)
\;\zeta(1/2+it+\alpha)\zeta(1/2+it+\beta)\nonumber\\
&&\qquad\qquad\zeta(1/2-it+\gamma)\zeta(1/2-it+\delta)
dt\bigg|_{\alpha=\beta=\gamma=\delta=0}.\nonumber
\end{eqnarray}

The most concise way to write the six terms that will result from
applying the ratio conjecture (actually in this case there is no
denominator so it is just a moment conjecture) to the fourth
moment of $\zeta$ in the last line above is as a contour integral
similar to (\ref{eq:IRa}) (as described in [10]). So, we have that

\begin{eqnarray}
I_{L}&=&I_R+\frac{d}{d\alpha}\frac{d}{d\beta}\frac{d}{d\gamma}
\frac{d}{d\delta}e^{\tfrac{\log
t/(2\pi)}{2}(-\alpha-\beta-\gamma-\delta)}\frac{1}{2\pi}\int_0^T
\log\tfrac{t}{2\pi}
\frac{1}{2!\;2!(2\pi i)^{4}}\nonumber\\
&&\quad\times\oint\cdots\oint e^{\tfrac{\log
t/(2\pi)}{2}(z_1+z_2-z_3-z_4)}
\frac{\zeta(1+z_1-z_3)\zeta(1+z_1-z_4) \zeta(1+z_2-z_3)
\zeta(1+z_2-z_4)
}{\zeta(2+z_1+z_2-z_3-z_4)}\\
&&\quad\times\frac{\Delta(z_1,\ldots,z_4)^2}{
\prod_{j=1}^4(z_j-\alpha)(z_j-\beta)(z_j+\gamma)(z_j+\delta)}
dz_1\cdots dz_4+O(T^{1/2+\epsilon}).\nonumber
\end{eqnarray}
Now we use the formulae (\ref{eq:deriv1}) and (\ref{eq:deriv2})
and arrive at
\begin{eqnarray}\label{eq:IL}
I_L&=&I_R+\frac{1}{2\pi}\int_0^T\log\tfrac{t}{2\pi}\frac{1}{2!\;2!(2\pi
i)^{4}}\nonumber\\
&&\qquad\times\oint\cdots\oint
\frac{\zeta(1+z_1-z_3)\zeta(1+z_1-z_4) \zeta(1+z_2-z_3)
\zeta(1+z_2-z_4)
}{\zeta(2+z_1+z_2-z_3-z_4)}\nonumber\\
&&\qquad\times\frac{\Delta(z_1,\ldots,z_4)^2}{
\prod_{j=1}^4z_j^4}\left(-\tfrac{\log t/(2\pi)}{2}-\sum_{j=1}^4
\frac{1}{z_j}\right)^2\left(-\tfrac{\log t/(2\pi)}{2}+\sum_{j=1}^4
\frac{1}{z_j}\right)^2\nonumber \\
&&\qquad\times e^{\tfrac{\log
t/(2\pi)}{2}(z_1+z_2-z_3-z_4)}dz_1\cdots dz_4+O(T^{1/2+\epsilon}).
\end{eqnarray}

We now compute the residues at $z_1=z_2=z_3=z_4=z_5=0$ of the
contour integrals in (\ref{eq:IR}) and (\ref{eq:IL}) (using
Mathematica).  If we write
\begin{equation}
\zeta(1+s)=\frac{1}{s}+\gamma-\gamma_1
s+\frac{\gamma_2}{2!}s^2-\frac{\gamma_3}{3!}s^3\cdots,
\end{equation}
then the final result is
\begin{theorem} Assuming the ratios conjecture as indicated in
(\ref{eq:IRa}), there are constants $C_0,\ldots,C_6$ such that
\begin{eqnarray}\label{eq:Itotal}
\sum_{\gamma< T}|\zeta'(\rho)|^4&=&\frac{1}{2\pi}\int_0^T
\bigg(\frac{\log^9 \tfrac{t}{2\pi}}{8640\zeta(2)}-\frac{(-2\gamma
\zeta(2)+\zeta'(2))\log^8
\tfrac{t}{2\pi}}{480\zeta^2(2)} \nonumber\\
&&\quad+\frac{(7\gamma^2 \zeta^2(2) - 2 \gamma_1\zeta^2(2)  -
    8 \gamma \zeta(2) \zeta'(2) + 4 (\zeta'(2))^2 -
    2  \zeta(2) \zeta''(2))\log^7\tfrac{t}{2\pi}}
    {120\zeta^3(2)}\nonumber\\
    &&\quad+C_6\log^6\tfrac{t}{2\pi} +\cdots+C_0 \bigg)
dt +O(T^{1/2+\epsilon})\nonumber\\
&=&\frac{T}{2 \pi}\frac{ \log^9 T}{8640\zeta(2)} +O(T\log^8T).
\end{eqnarray}
\end{theorem}
This leading term agrees with Conjecture \ref{conj:hko} of Hughes,
Keating and O'Connell [23]; since $G^2(4)/G(7)=8640$ and
$a(2)=1/\zeta(2)$, we see that there is agreement between the
leading order term of (\ref{eq:Itotal}) and the conjecture in the
case that $k=2$.

The first three terms in decreasing powers of $\log T$, shown
in (\ref{eq:Itotal}) contain only arithmetic factors that are
derivatives of the Riemann zeta function evaluated at 2.  This is
not true of the $(\log T)^6$ term.  We can see that this will be
the case by the comment after (\ref{eq:gamma0}).  The arithmetic
factor which forms part of
$T(\alpha,\beta,\gamma,\delta,\epsilon)$ in (\ref{eq:T}),
\begin{eqnarray}
&&A(\alpha,\beta,\gamma,\delta,\epsilon)\nonumber\\&&\qquad=\prod_p
\left[
\frac{(1-\tfrac{1}{p^{1+\alpha+\delta}})(1-\tfrac{1}{p^{1+\alpha+\epsilon}})
(1-\tfrac{1}{p^{1+\beta+\delta}})(1-\tfrac{1}{p^{1+\beta+\epsilon}})
(1-\tfrac{1}{p^{1+\gamma+\delta}})(1-\tfrac{1}{p^{1+\gamma+\epsilon}})}
{(1-\tfrac{1}{p^{1+\delta}})(1-\tfrac{1}{p^{1+\epsilon
}})}\right.\\
&&\left.\qquad\qquad\times \sum_{m_1+m_2+m_3+h=n_1+n_2}
\frac{\mu(p^h)}{p^{(1/2+\alpha)m_1+(1/2+\beta)m_2+(1/2+\gamma)m_3+h/2+(1/2+\delta)n_1
+(1/2+\epsilon)n_2}}\right],\nonumber
\end{eqnarray}
expands as a Taylor series
\begin{eqnarray}
\label{eq:taylorA} &&A(z_1,z_2,z_3,-z_4,-z_5)=A_0 +
 A_1(z_1 + z_2 + z_3 - z_4 -
z_5) \\
&&\quad +   A_{12}(-z_1 z_4 - z_2 z_4 - z_3 z_4 - z_1 z_5
     - z_2 z_5 - z_3 z_5+ z_1 z_2 +
    z_1 z_3 + z_2 z_3 + z_4 z_5)\nonumber \\
    &&\qquad+ \frac{A_{11}}{2}(z_1^2 + z_2^2
    + z_3^2 + z_4^2 + z_5^2)\nonumber \\
    &&\qquad\quad+
    A_{124}(-z_1 z_2 z_4 - z_1 z_3 z_4 - z_2 z_3 z_4 - z_1 z_2 z_5 -
          z_1 z_3 z_5 - z_2 z_3 z_5 + (z_1 + z_2 + z_3) z_4 z_5)
        \nonumber  \\
          &&\qquad\qquad+A_{123}z_1z_2z_3 \nonumber\\ &&\qquad\qquad\quad+
    \frac{A_{112}}{2}(z_4^2(z_1 + z_2 + z_3 - z_5) +
          z_3^2(z_1 + z_2 - z_4 - z_5) + z_2^2(z_1 + z_3 -
          z_4 - z_5) \nonumber\\&&\qquad\qquad\qquad\qquad\qquad\qquad+
          z_1^2(z_2 + z_3 - z_4 - z_5) + (z_1 + z_2 + z_3 -
                z_4) z_5^2)\nonumber \\&&\qquad\qquad\qquad+
    \frac{A_{111}}{6}(z_1^3 + z_2^3 + z_3^3 - z_4^3 -
    z_5^3)+\cdots,\nonumber
    \end{eqnarray}
where $A_j$ is the partial derivative, evaluated at zero, of
$A(z_1,z_2,z_3,z_4,z_5)$ with respect to the $j$th variable.  Note
that $A$ is symmetric amongst the first three variables, and
amongst the final two variables, so for example
$A_{12}=\frac{\partial A(z_1,z_2,z_3,z_4,z_5)}{\partial z_1
\partial z_2}\big|_{z_1=z_2=z_3=z_4=z_5=0} = A_{13}=A_{23}$.  In addition, we noted at
(\ref{eq:gamma0}) that $A(\alpha,\beta,0,\delta,\gamma)$ is just
the same as the arithmetic factor from the fourth moment of zeta,
$\frac{1}{T}\int_{0}^T
\zeta(1/2+it+\alpha)\zeta(1/2+it+\beta)\zeta(
1/2-it+\delta)\zeta(1/2-it+\epsilon)dt$, that is,
$1/\zeta(2+\alpha+\beta+\delta+\epsilon)$.  Therefore all the
partial derivatives of $A$ in (\ref{eq:taylorA}) can be computed
by taking partial derivatives of $1/\zeta(2+z_1+z_2+z_4+z_5)$
except for the derivative $A_{123}=\frac{\partial^3 A}{\partial
z_1\partial z_2\partial z_3}$, which involves all of the first
three variables and gives a contribution that is not expressed as
a derivative of $\zeta(2)$.  This contribution shows up first in
the $\log^6 T$ term.

\subsection{A second moment:
$\sum_{0<\gamma<T}|\zeta(\rho+a)|^2$}

In [18] A. Fujii generalizes work of [20] and proves, under the
assumption of the Riemann Hypothesis, the following theorem.

\begin{theorem}[{Fujii [18]}] \label{theo:fujii}
Assume the Riemann Hypothesis is true.  If $T$ is sufficiently
large and $\alpha$ is a real number such that $|\alpha|\ll \log
T$, then
  \begin{eqnarray}
&&\sum_{1\leq\gamma\leq T} |\zeta(\tfrac{1}{2} +i(\gamma
+2\pi\alpha/\log\tfrac{T}{2\pi}))|^2=\left(1-\left(\frac{\sin\pi\alpha}{\pi\alpha}
\right)^2\right) \frac{T}{2\pi} \log^2T \\
&&\qquad +2\left(-1+\gamma+(1-2\gamma) \frac{\sin 2\pi \alpha}{2\pi \alpha} +\Re \left(\frac{\zeta'}{\zeta}\bigg(1+i\frac{2\pi \alpha}{\log(T/2 \pi)}\bigg)\right)\right)\frac{T}{2\pi}\log \frac{T}{2\pi}  \nonumber \\
&&\qquad +G(T,\alpha) +O(\sqrt{T}\log^3T),  \nonumber
\end{eqnarray}
where $\gamma$ is Euler's constant and $G(T,\alpha)=O(T)$ is
explicitly given.
  \end{theorem}

The ratios conjecture can reproduce this result in a
straightforward way.

First we write the quantity we want to calculate in a form in
which we can apply the ratios conjecture:
\begin{equation}
\sum_{0<\gamma<T} |\zeta(\rho+a)|^2=\lim_{a_1,a_2\rightarrow a}
\frac{1}{2\pi i} \int_C \frac{\zeta'(z)}{\zeta(z)}
\zeta(z+a_1)\zeta(1-z-a_2) dz.
\end{equation}
Now we follow the identical method to that described in Section
\ref{sect:2ndderiv} and so obtain
\begin{eqnarray}
\sum_{0<\gamma<T} |\zeta(\rho+a)|^2&=&\lim_{a_1,a_2\rightarrow
a}\Bigg[ \frac{d}{d\alpha} \frac{1}{2\pi} \int_0^T
\frac{\zeta(1+\alpha-a_2)\zeta(1+a_1-a_2)}{\zeta(1-a_2)}\nonumber\\
&&\qquad\qquad\qquad\qquad+
\left(\frac{t}{2\pi}\right)^{-\alpha+a_2}\frac{\zeta(1+a_2-\alpha)\zeta(1
+a_1-\alpha)}{\zeta(1-\alpha)}\nonumber\\
&&\qquad\qquad\qquad\qquad+
\left(\frac{t}{2\pi}\right)^{-a_1+a_2}\frac{\zeta(1+\alpha-a_1)\zeta(1+a_2-a_1)}
{\zeta(1-a_1)} dt \big|_{\alpha=0} \nonumber\\
&&\qquad\quad+\frac{d}{d\alpha} \frac{1}{2\pi} \int_0^T
\frac{\zeta(1+\alpha+a_1)\zeta(1-a_2+a_1)}{\zeta(1+a_1)}\\
&&\qquad\qquad\qquad\qquad+\left(\frac{t}{2\pi}\right)^{-\alpha-a_1}\frac{\zeta(1-a_1-\alpha)\zeta(1-a_2-\alpha)}
{\zeta(1-\alpha)}\nonumber\\
&&\qquad\qquad\qquad\qquad +\left(\frac{t}{2\pi}\right)^{a_2-a_1}
\frac{\zeta(1+\alpha+a_2)\zeta(1-a_1+a_2)}{\zeta(1+a_2)} dt
\big|_{\alpha=0} \nonumber\\
&&\qquad\quad+\frac{1}{2\pi} \int_0^T \log\tfrac{t}{2\pi} \left(
\zeta(1-a_2+a_1)+\left(\frac{t}{2\pi}\right)^{-a_1+a_2}
\zeta(1-a_1+a_2)\right)dt\Bigg]\nonumber\\
&&\qquad\qquad\qquad\qquad +O(T^{1/2+\epsilon}).\nonumber
\end{eqnarray}
Performing the differentiation and setting $\alpha=0$, we have
\begin{eqnarray}
\sum_{0<\gamma<T} |\zeta(\rho+a)|^2&=&\lim_{a_1,a_2\rightarrow
a}\frac{1}{2\pi}\int_0^T
\frac{\zeta'(1-a_2)\zeta(1+a_1-a_2)}{\zeta(1-a_2)}
-\left(\frac{t}{2\pi}\right)^{a_2} \zeta(1+a_2)\zeta(1+a_1)\nonumber\\
&&\qquad\qquad\qquad\qquad
+\left(\frac{t}{2\pi}\right)^{-a_1+a_2}\frac{\zeta'(1-a_1)\zeta(1+a_2-a_1)}
{\zeta(1-a_1)} \nonumber\\
&&\qquad\qquad+\frac{\zeta'(1+a_1)\zeta(1-a_2+a_1)} {\zeta(1+a_1)}
-\left(\frac{t}{2\pi}\right)^{-a_1}\zeta(1-a_1)\zeta(1-a_2)\\
&&\qquad\qquad\qquad\qquad
+\left(\frac{t}{2\pi}\right)^{a_2-a_1}\frac{\zeta'(1+a_2)\zeta(1-a_1+a_2)}
{\zeta(1+a_2)} \nonumber\\
&&\qquad\qquad+\log \tfrac{t}{2\pi}
\left(\zeta(1-a_2+a_1)+\left(\frac{t}{2\pi}\right)^{-a_1+a_2}\zeta(1-a_1+a_2)\right)dt\nonumber\\
&&\qquad\qquad\qquad\qquad +O(T^{1/2+\epsilon}).\nonumber
\end{eqnarray}
To perform the limit, let $a_1=a$ and $a_2=a+s$.  Then
\begin{equation}
\lim_{s\rightarrow
0}\zeta(1-s)+\left(\frac{t}{2\pi}\right)^s\zeta(1+s) = \log
\tfrac{t}{2\pi} +2\gamma,
\end{equation}
and
\begin{eqnarray}
&&\lim_{s\rightarrow0}
\frac{\zeta'(1-a-s)\zeta(1-s)}{\zeta(1-a-s)}
+\left(\frac{t}{2\pi}\right)^s
\frac{\zeta'(1-a)\zeta(1+s)}{\zeta(1-a)}\\
&&\qquad\qquad =(\log\tfrac{t}{2\pi}+2\gamma)
\frac{\zeta'(1-a)}{\zeta(1-a)} +\frac{\zeta''(1-a)}{\zeta(1-a)}
-\frac{(\zeta'(1-a))^2} {\zeta^2(1-a)}.\nonumber
\end{eqnarray}
Thus, assuming the ratios conjecture as indicated in
(\ref{eqn:IR}), we have
\begin{eqnarray}
\sum_{0<\gamma<T}
|\zeta(\rho+a)|^2&=&\frac{1}{2\pi}\int_0^T(\log\tfrac{t}{2\pi}
+2\gamma)\left(\log\tfrac{t}{2\pi} +\frac{\zeta'(1-a)}{\zeta(1-a)}
+\frac{\zeta'(1+a)}{\zeta(1+a)}\right)\\
&&+ \frac{\zeta''(1-a)}{\zeta(1-a)}
+\frac{\zeta''(1+a)}{\zeta(1+a)}-\left(\frac{\zeta'(1-a)}{\zeta(1-a)}\right)^2
-\left(\frac{\zeta'(1+a)}{\zeta(1+a)}\right)^2\nonumber\\
&& -\left(\frac{t}{2\pi}
\right)^a\zeta(1+a)\zeta(1+a)-\left(\frac{t}{2\pi}\right)^a
\zeta(1-a)\zeta(1-a)\;dt\nonumber\\
&&\qquad\qquad\qquad\qquad +O(T^{1/2+\epsilon}).\nonumber
\end{eqnarray}

This result matches up exactly with Theorem \ref{theo:fujii}; see
in particular the bottom of page 66 in [18].

If we now let $a=2\pi i\alpha\log^{-1}\tfrac{T}{2\pi}$, then for
large $T$ we can use the first few terms of the series
$\zeta(1+s)=\tfrac{1}{s}+\gamma-\gamma_1s+\gamma_2\tfrac{s^2}{2}+\cdots$,
as well as the similar expressions for the derivatives and inverse
of $\zeta(1+s)$, and perform the integration over $t$ to obtain a
more standard expression for the leading terms:
\begin{eqnarray}
&&\sum_{0<\gamma<T}
|\zeta(\rho+2\pi i\alpha \log^{-1}\tfrac{T}{2\pi})|^2\\
&&\qquad=\left(1-\left(\frac{\sin
\pi\alpha}{\pi\alpha}\right)^2\right)
\frac{T}{2\pi}\log^2\tfrac{T}{2\pi}
+\frac{T}{2\pi}\log\tfrac{T}{2\pi}\left(\frac{\sin2\pi\alpha}{\pi\alpha}
-2\gamma\frac{\sin2\pi\alpha}{\pi\alpha}+4\gamma-2\right)\nonumber\\
&&\qquad+\frac{T}{2\pi}(4\gamma\cos
2\pi\alpha-2\cos2\pi\alpha-2\gamma^2\cos2\pi\alpha-4\gamma_1\cos2\pi\alpha+2\gamma^2-4\gamma
+2)+o(T).\nonumber
\end{eqnarray}

\section{Further connections with the literature}\label{sect:conn}
\subsection{Non-vanishing of Dirichlet $L$-functions}

Michel and VanderKam's paper ``Non-vanishing of high derivatives
of Dirichlet's $L$-functions at the central point'' [35] actually
is concerned with non-vanishing of derivatives of
$\Lambda(s,\chi)$, the completed Dirichlet $L$-functions:
\begin{equation}
\Lambda(s,\chi)=\big(\frac{q}{\pi}\big)^{s/2}\Gamma(s/2)L(s,\chi),
\end{equation}
which satisfies
\begin{equation}
\Lambda(s,\chi)=\epsilon_{\chi} \Lambda(1-s,\overline{\chi}).
\end{equation}
  They use mollifying techniques to give a lower
bound for the frequency of $\Lambda^{(k)}(1/2,\chi)\ne 0$ as
$\chi$ ranges over primitive characters modulo $q$. They find (a)
that using a mollifier with two pieces
\begin{eqnarray}
M_k(1/2, \chi):=\sum_{m\le y} \frac{\mu(m)\left(\chi(m)+(-1)^k
\overline{\epsilon}_{\chi}(\chi)\overline{\chi}(m)\right)P\left(\frac{\log
\frac y m}{\log y}\right)}{m^{1/2}}
\end{eqnarray}
is more efficient than the conventional mollifier
\begin{eqnarray}
M^*_k(1/2,\chi):=\sum_{m\le y}
\frac{\mu(m)\chi(m)P\left(\frac{\log \frac y m}{\log
y}\right)}{m^{1/2}};
\end{eqnarray}
in fact for $\Lambda(1/2,\chi)$ (no derivatives) they find that
asymptotically 1/2 do not vanish, improving work of Iwaniec and
Sarnak [25] who had 1/3 in place of 1/2 here. They also prove (b)
that when mollifying the high derivatives of $\Lambda$, the
proportion of non-vanishing of the $k$th derivative can be shown
to approach 2/3 as $k\to \infty$, and that to obtain this result
it is critical to use the general mollifier.

In this section we use the ratios conjecture to reproduce the
asymptotic formulae of [35].  In addition we indicate how one can
show that the proportion of non-vanishing for $L^{(k)}(1/2,\chi)$
does approach 100\% as $k\to\infty$ for this family.

Rather than work with the Dirichlet $L$-functions, we find it
convenient to work with the Riemann zeta-function in $t$ aspect;
they are both unitary families, so the results will be identical.
For the analogue of the $\Lambda^{(k)}(1/2,\chi)$ function we will
use the function
\begin{eqnarray}
\chi(s)^{1/2} \mathcal Z^{(k)}(s),
\end{eqnarray}
where \begin{eqnarray} \mathcal Z(s)=\chi(s)^{-1/2} \zeta(s)
\end{eqnarray}
is a complex analogue of Hardy's $Z(t)$ function.  This is
appropriate because the $\mathcal{Z}$-function associated with
$L(s,\chi)$ is
\begin{equation}
\mathcal{Z}(s,\chi)=\epsilon_{\chi}^{-1/2}\Lambda(s,\chi),
\end{equation}
so that
\begin{equation}
\Lambda^{(k)}(s,\chi)=\epsilon_{\chi}^{1/2}\mathcal{Z}^{(k)}(s,\chi).
\end{equation}
Note that $\chi$ is used in two different roles here; recall that
$\chi(s)$ is the factor from the functional equation of the zeta
function, see (\ref{eq:approxfunczeta}), and it plays the role of
$\epsilon_{\chi}$.

The analogue of the quantity considered in [35] (in Section 7 of
that paper) is
\begin{eqnarray}\label{eq:twopiece}
\frac{\left(\int_0^T \chi(s)^{1/2} \mathcal
Z^{(k)}(s)M(s)~dt\right)^2} {\int_0^T \mathcal Z^{(k)}(s)\mathcal
Z^{(k)}(1-s)M(s)M(1-s)~dt},
\end{eqnarray}
for a two-piece mollifier.  However, here we will illustrate the
calculation with the conventional mollifier
\begin{equation}\label{eq:MSP}
M(s)=M(s,P)=\sum_{m\leq y} \mu(m) P(\log (y/m)/\log y)/m^s.
\end{equation} The object is to choose $P$ is such a way that this ratio
is minimized.

 Since
\begin{eqnarray}
\chi(s)=\left(\frac t {2\pi}\right)^{1/2-s}(1+O(1/t))
\end{eqnarray}
for $t>1$, we have
\begin{eqnarray}
\frac d {ds} \left(\chi(s)^{-1/2}\right)=\frac{\ell}{2}
\chi(s)^{-1/2}(1+O(1/t)),
\end{eqnarray}
with $\ell=\log \tfrac{t}{2\pi}$.  Hence,
\begin{eqnarray}
\chi(s)^{1/2} \mathcal Z^{(k)}(s)=\left(\frac{d}{d\alpha}\right)^k
e^{\alpha \ell/2} \zeta(s+\alpha)(1+O(1/t))\bigg|_{\alpha=0}.
\end{eqnarray}
Thus the integral in the numerator can be evaluated by considering
\begin{eqnarray}
\left(\frac {d}{d\alpha}\right)^k \int_0^T e^{\alpha \ell /2}
\mathcal \zeta(s+\alpha)M(s,P)~dt\bigg|_{\alpha=0}\sim \left(\frac
{d}{d\alpha}\right)^k T e^{\alpha L/2}
P(1)\bigg|_{\alpha=0}=TP(1)2^{-k}L^k
\end{eqnarray}
where $L=\log T$. (The ratios conjecture gives $\int_0^T
\frac{\zeta(s+\alpha)}{\zeta(s+w)}~dt \sim T$.) The denominator
is, by (\ref{eqn:mollunit}) and after rescaling $\alpha = a/L,
\beta=b/L$,
\begin{eqnarray}&&
\quad\sim T L^{2k} \left(\frac {d}{da}\right)^k\left(\frac
{d}{db}\right)^k  e^{(a+b)/2}\nonumber\\
&&\qquad \qquad \times\left(P(1)^2+\frac{1}{\theta}
\frac{d}{dw}\frac{d}{dz} e^{-a\theta w -b \theta  z
}\int_0^1\int_0^1 e^{-(a+b)u} P (w+r)P
(z+r)~dr~du\bigg|_{w=z=0 \atop a=b=0}\right)\nonumber\\
&& \quad = TL^{2k}\left( \frac{P(1)^2}{2^{2k}} \right.\nonumber\\
&&\qquad\qquad\left.+\frac 1\theta\frac
d{dw}\frac{d}{dz}\int_0^1\int_0^1 P(r+w)(1/2-u-\theta w)^k
P(r+z)(1/2-u-\theta z)^k ~dr~du\right)\bigg|_{w=z=0}\nonumber\\
&&\quad=TL^{2k}\bigg( \frac{P(1)^2}{2^{2k}} + \frac
1\theta\int_0^1\int_0^1 \bigg(P'(r)(1/2-u )^k-k\theta P(r)(1/2-u
)^{k-1}\bigg)^2 ~dr~du\bigg)\nonumber\\
&&\quad=\frac{TL^{2k}}{2^{2k}}\bigg(  P(1)^2  + \frac
1\theta\int_0^1
 \frac{P'(r)^2}{2k+1}~dr + 4\theta
 \int_0^1\frac{k^2P(r)^2}{2k-1}~dr\bigg).
 \end{eqnarray} This corresponds to the evaluation of $\mathcal{Q}_1$, accomplished in equation (17)
 in [35].  Note that $\Delta=2\theta$.

Thus, the ratio (\ref{eq:twopiece})  is
\begin{eqnarray}\label{eq:ratio}
\frac{P(1)^2}{P(1)^2  + \frac 1\theta\int_0^1
 \frac{P'(r)^2}{2k+1}~dr + 4\theta
 \int_0^1\frac{k^2P(r)^2}{2k-1}~dr}.
 \end{eqnarray}
 If $k=0$ we take $P(r)=r$ and $\theta=1/2$ and deduce that at
 least 1/3 of $L$-functions do not vanish at 1/2.  This is the result of Iwaniec and Sarnak [25]. For large
$k$ if we take $P(r)=r^k$ we see that this ratio is
\begin{eqnarray}
\frac{1}{1+\left(4\theta+\frac{1}{\theta}\right)\frac{k^2}{4k^2-1}}\sim\frac{1}{1+\left(\theta+\frac
{1}{4\theta}\right)}
\end{eqnarray}
which is 1/2 when $\theta=1/2$.

In general if $A,B>0$, the minimum of $A\int_0^1
P'(x)^2~dx+B\int_0^1 P(x)^2~dx$ over smooth functions $P$
satisfying $P(0)=0$ and $P(1)=1$ is $AP'(1)$ and is achieved by
$P(x)=(\sinh \sqrt{B/A} x)/(\sinh\sqrt{B/A})$. So the optimal
choice for (\ref{eq:ratio}) is
\begin{eqnarray}
P(r)=\frac{\sinh(\Lambda r)}{\sinh \Lambda}, \qquad
\Lambda=2\theta k \sqrt{\frac{2k+1}{2k-1}}
\end{eqnarray}
as in [35]; however, this still gives that the ratio is $1/2+
O(1/k^2)$.

Next we explain the use of the two part mollifier in the case that
$k=0$. For this, we consider a mollifier of the form
\begin{eqnarray}
M(s,P,a):=\sum_{n\le y} \mu(n)P\left(\frac {\log \frac y n}{\log
y}\right)(n^{-s}+a\chi(1-s)n^{s-1})
\end{eqnarray}
and we want to maximize the ratio
\begin{eqnarray}
\frac{\left(\int_0^T \zeta(s) M(s,P,a)~dt\right)^2}{\int_0^T
\zeta(s)\zeta(1-s)M(s,P,a)M(1-s,P,a)~dt}.
\end{eqnarray}
The key things to observe here are that, with $M(s,P)$ as in
(\ref{eq:MSP}),
\begin{eqnarray}
\int_0^T \zeta(s)\chi(1-s)M(1-s,P)~dt=\int_0^T
\zeta(1-s)M(1-s,P)~dt \sim T P(1)
\end{eqnarray}
and \begin{eqnarray} &&\int_0^T \zeta(s)\zeta(1-s)
M(1-s,P)\chi(1-s)M(1-s,P)~dt\\&&\qquad=\int_0^T
\zeta(1-s)^2M(1-s,P)^2~dt\nonumber\\
&&\qquad \sim P(1)^2 T.\nonumber
\end{eqnarray}
Thus, the ratio is
\begin{eqnarray}
\sim T
\frac{(1+a)^2P(1)^2}{(1+a^2)\left(P(1)^2+\frac{1}{\theta}\int_0^1
P'(t)^2~dt\right) +2 aP(1)^2}.
\end{eqnarray}
The optimal choices here are $P(r)=r$ and $a=1$; for $\theta=1/2$
this gives a ratio of 1/2 as claimed in [35].

To handle high derivatives of $L(s,\chi)$ at $s=1/2$ we consider,
by analogy, high derivatives of $\zeta(s)$. The trick is to insert
a factor of $\chi(s)$ and to ask about the non-vanishing of
\begin{eqnarray}\chi(s)\zeta^{(k)}(1-s) .
\end{eqnarray}
Thus, we want to maximize the ratio of
\begin{eqnarray}
\frac{\left(\int_0^T
\chi(s)\zeta^{(k)}(1-s)M(s,P)~ds\right)^2}{\int_0^T
|\chi(s)\zeta^{(k)}(1-s)M(s,P)|^2~dt}.
\end{eqnarray}
Now,
\begin{eqnarray}
\zeta^{(k)}(1-s)&=&\left(\frac{d}{ds}\right)^k\chi(1-s)\zeta(s)=\sum_{j=0}^k
\left({k \atop j}\right)
(-1)^j\chi^{(j)}(1-s) \zeta^{(k-j)}( s)\\
&=&  \chi(1-s) \left(\frac{d}{d\alpha}\right)^k
e^{\alpha\ell}\zeta(s+\alpha)\bigg|_{\alpha=0}.
\end{eqnarray}
Thus, the numerator is
\begin{eqnarray} \sim\ (TP(1)L^k)^2.
\end{eqnarray}
The denominator is evaluated optimally in [14] and is
$T|P(1)|^2L^{2k}(1-O(1/k^2))$. Applying this method to $L(s,\chi)$
it can be deduced that there is a constant $C>0$ such that the
proportion of $L^{(k)}(1/2,\chi)$ which vanish is smaller than
$C/k^2$.

\begin{remark}
Michel and VanderKam [35] give a nice explanation at the end of
Section 2 of why one cannot expect to do better than 1/2
non-vanishing of $\Lambda^{(k)}(1/2,\chi)$ using a conventional
mollifier and 2/3 using a two-piece mollifier.  The reason relies
on the symmetry of the approximate functional equation for
$\Lambda^{(k)}(1/2,\chi)$ and the uniform distribution of
$\epsilon_{\chi}$.  That we can get a proportion of non-vanishing
of $L^{(k)}(1/2,\chi)$ approaching 1 as $k\to \infty$ does not
contradict their argument because of the lack of symmetry of the
approximate functional equation for $L^{(k)}(1/2,\chi)$.  For
applications to bounding multiplicities of central zeros,
information about non-vanishing of $L^{(k)}(1/2,\chi)$ is equally
as good as for $\Lambda^{(k)}(1/2,\chi)$.
\end{remark}

\subsection{Non-vanishing of automorphic $L$-functions}

The main theorem of Kowalski, Michel, and VanderKam's paper
``Non-vanishing of high derivatives of automorphic  $L$-functions
at the center of the critical strip'' [31] is a mollification of
the second moment of weight 2 primitive cusp forms of a prime
level. The formula achieved is slightly different than the result
we mention above (Theorem \ref{theo:mollorthog}) for mollifying
the second moment in an orthogonal family.  The reason for this is
that they use a slightly different mollifier. Instead of choosing
a smoothed sum of the coefficients of the inverses of the
Dirichlet series in question they choose a mollifier of the shape
\begin{eqnarray}
\sum_{m\le y}\frac{\lambda_f(n)\mu(n)P\left(\frac{\log \frac y
n}{\log y} \right)}{\psi(n)n^{1/2}},
\end{eqnarray}
where $\lambda_f(n)$ are the coefficients of the $L$-function
which is to be mollified and where $\psi(n)=\prod_{p\mid n}
\left(1+\frac{1}{p}\right).$  The analogue of our Theorem
\ref{theo:mollorthog} has the right side replaced by
\begin{eqnarray}
\qquad\frac{1}{\theta^2}\left((Q(1)P'(1)+\theta
Q'(1)P(1))^2+\frac{1}{\theta}\int_0^1\int_0^1
(P''(x)Q(y)-\theta^2P(x)Q''(y))^2~dx~dy\right).
\end{eqnarray}
This result, which is not deducible from our ratios conjecture,
was reported in the paper of Conrey and Farmer [8] as the general
result one would obtain from mollifying a second moment in an
orthogonal family. We wish to correct that statement and replace
it with the statement of Theorem \ref{theo:mollorthog}.
\subsection{Non-vanishing of quadratic $L$-functions}   The papers of Soundararajan and
Conrey \& Soundararajan deal with non-vanishing of Dirichlet
$L$-functions for real quadratic characters, at the central point
and on the real axis. The results of Theorem
\ref{theo:mollsymplec} are consistent with the results of these
papers [40,15].

\section{Conclusion}
The purpose of this paper was to illustrate the use of the ratios
conjectures by deriving from them a number of important results
from the theory of $L$-functions.  The variety of applications is
by no mean exhausted by what we have presented.  Other
calculations that might be valuable include lower order terms in
moments of $S(t)$, $\log |\zeta(1/2+it)|$ and  $S(t+h)-S(t)$. For
the second moment of $S(t)$ the lower order terms have already
been computed by Tsz Ho Chan [6], while  lower order terms of the
second moment of $S(t+h)-S(t)$ have been considered in [1].
Precise evaluations of $n$-level correlations  might be combined
to obtain the secondary terms in the nearest neighbour spacing
distribution for the zeros of the Riemann zeta function. In [39],
the leading term in the $n$-correlation function is calculated for
a restricted space of test functions, but for essentially any
$L$-function. Ratios conjectures could also be used to evaluate
possible schemes to improve lower bounds for proportions of zeros
on the critical line.

\end{document}